\documentclass[mnsc,nonblindrev]{informs4} 
\OneAndAHalfSpacedXI

\usepackage{subcaption}
\usepackage{mathrsfs}
\usepackage{amsmath,bm}
\usepackage{graphicx}
\usepackage{array}
\usepackage{tabularx}
\usepackage{url}
\definecolor{DarkBlue}{rgb}{0,0.08,0.45}
\usepackage[backref = false, bookmarks, colorlinks = true, plainpages = false, citecolor = DarkBlue , urlcolor = DarkBlue, filecolor = DarkBlue, linkcolor = DarkBlue]{hyperref}
\usepackage{newunicodechar}
\usepackage[utf8]{inputenc}
\usepackage{natbib}
 \bibpunct[, ]{(}{)}{,}{a}{}{,}%
 \def\bibfont{\small}%
 \def\bibsep{\smallskipamount}%

 \def\bd{\boldsymbol{d}}
 \def\bx{\boldsymbol{x}}
 \def\bW{\boldsymbol{W}}

 \def\bbeta{\boldsymbol{\beta}}

\TheoremsNumberedThrough     
\ECRepeatTheorems
\graphicspath{{plots/}}
\EquationsNumberedThrough    

\begin{document}
\newcommand{\abs}[1]{\left|  #1 \right| }
\newcommand{\brak}[1]{\left(#1\right)}    
\newcommand{\crl}[1]{\left\{#1\right\}}   
\newcommand{\edg}[1]{\left[#1\right]}     
\newcommand{\norm}[1]{\|#1\|}
\newcommand{\floor}[1]{\lfloor #1 \rfloor}

\newcommand{\cA}{{\mathcal A}}
\newcommand{\cB}{{\mathcal B}}
\newcommand{\cD}{{\mathcal D}}
\newcommand{\cF}{{\mathcal F}}
\newcommand{\cG}{{\mathcal G}}
\newcommand{\cH}{{\mathcal H}}
\newcommand{\cK}{{\mathcal K}}
\newcommand{\cL}{{\mathcal L}}
\newcommand{\cM}{{\mathcal M}}
\newcommand{\cR}{{\mathcal R}}
\newcommand{\cS}{{\mathcal S}}
\newcommand{\cT}{{\mathcal T}}
\newcommand{\cX}{{\mathcal X}}

\newcommand{\bb}{{\mathbf b}}
\newcommand{\bp}{{\mathbf p}}
\newcommand{\bc}{{\mathbf c}}
\newcommand{\bg}{{\mathbf g}}
\newcommand{\bl}{{\mathbf l}}
\newcommand{\bbf}{{\mathbf f}}
\newcommand{\by}{{\mathbf y}}

\newcommand{\bA}{{\mathbf A}}
\newcommand{\bB}{{\mathbf B}}
\newcommand{\bC}{{\mathbf C}}
\newcommand{\bD}{{\mathbf D}}
\newcommand{\bG}{{\mathbf G}}
\newcommand{\bL}{{\mathbf L}}
\newcommand{\bS}{{\mathbf S}}
\newcommand{\bQ}{{\mathbf Q}}
\newcommand{\bU}{{\mathbf U}}
\newcommand{\bV}{{\mathbf V}}
\newcommand{\bX}{{\mathbf X}}
\newcommand{\bZ}{{\mathbf Z}}
\newcommand{\bF}{{\mathbf F}}
\newcommand{\ba}{{\mathbf a}}
\newcommand{\bz}{{\mathbf z}}

\newcommand{\bmu}{{\boldsymbol \mu}}
\newcommand{\balpha}{{\boldsymbol \alpha}}

\newcommand{\btheta}{\boldsymbol{\theta}}
\newcommand{\bsigma}{\boldsymbol{\sigma}}
\newcommand{\bnu}{\boldsymbol{\nu}}
\newcommand{\bSigma}{\boldsymbol{\Sigma}}

\newcommand{\F}{\mathbb{F}}
\newcommand{\p}{\mathbb{P}}
\newcommand{\Q}{\mathbb{Q}}

\newcommand{\R}{\mathbb{R}}
\newcommand{\q}{\mathbb{Q}}

\newcommand{\tr}{{\rm tr}}

\newcommand{\id}{{\mathbbm 1}}

\newcommand{\expect}{\mathbb{E}}


\RUNAUTHOR{Han, Hu, and Shen}

\RUNTITLE{
Deep Neural Newsvendor
}

\TITLE{\Large Deep Neural Newsvendor}
\vspace{-1.6cm}

\ARTICLEAUTHORS{%
\AUTHOR{Jinhui Han}
\AFF{Rotman School of Management, University of Toronto,  Toronto, Ontario, Canada M5S 3E6, \EMAIL{jinhuizjcu@gmail.com}} 
\AUTHOR{Ming Hu}
\AFF{Rotman School of Management, University of Toronto,  Toronto, Ontario, Canada M5S 3E6, \EMAIL{ming.hu@rotman.utoronto.ca}}
\AUTHOR{Guohao Shen}
\AFF{Department of Applied Mathematics, Hong Kong Polytechnic University, Hung Hom, Kowloon, Hong Kong, \EMAIL{guohao.shen@polyu.edu.hk}}
} 

\vspace{-1.6cm}

\ABSTRACT{%
We consider a data-driven newsvendor problem, where one has access to past demand data and the associated feature information. We solve the problem by estimating the target quantile function using a deep neural network (DNN). The remarkable representational power of DNN allows our framework to incorporate or approximate various extant data-driven models. We provide theoretical guarantees in terms of excess risk bounds for the DNN solution characterized by the network structure and sample size in a non-asymptotic manner, which justify the applicability of DNNs in the relevant contexts. Specifically, the convergence rate of the excess risk bound with respect to the sample size increases in the smoothness of the target quantile function but decreases in the dimension of feature variables. This rate can be further accelerated when the target function possesses a composite structure. 
In particular, our theoretical framework can be extended to accommodate the data-dependent scenarios, where the data-generating process could be time-dependent but not necessarily identical over time. 
Building on our theoretical results, we provide further managerial insights and practical guidance through simulation studies. 
Finally, we apply the DNN method to a real-world dataset obtained from a food supermarket. Our numerical experiments demonstrate that (1) the DNN method consistently outperforms other alternatives across a wide range of cost parameters, and (2) it exhibits good performance when the sample size is either very large or relatively limited.
}%



\maketitle

%



\section{Introduction}\label{sec:introduction}

The newsvendor problem manages to find the optimal inventory level to serve random demand in balancing the overage or underage costs if the inventory is higher or lower than the realized demand. The accurate prediction of demand (distribution) is essential in newsvendor decision-making. In particular, the demand may be influenced by a variety of factors, such as customer demographics, seasonality, and general economic indexes. Thanks to the big data era, the availability of historical demand data and observations on associated features renders data-driven demand prediction viable and effective. Motivated by this viewpoint, we investigate the newsvendor problem in a data-driven context where the observed feature variables are deployed to better predict unknown demand.

Specifically, we propose to solve the data-driven newsvendor problem using a one-step, distribution-free, nonparametric deep neural network (DNN) method. In recent decades, DNN has achieved tremendous success beyond the computer science community in various fields, such as medical science, engineering, and finance. However, its applications in the context of operations management (OM) are still relatively limited compared to its popularity in other areas, especially for exploring the theoretical properties that support its efficacy in solving OM problems \citep{feng2018research}. In our setting, the inherent structure of the newsvendor problem allows us to express its solution as a conditional quantile function of demand-related features, which can be conveniently approximated by DNNs based on their remarkable representation power. In addition, the DNN method also enjoys the advantages of handling high-dimensional data over other alternatives such as the parametric two-step, estimate-then-optimize method.

Our data-driven inventory decision is made by feeding the observed feature information as inputs into a fitted DNN with the optimized inventory decision under the empirical risk minimization as an output, a process analogous to the high-dimensional quantile regression. For our nonparametric DNN model, we can properly design the network structure, such as the network depth and width, that essentially determines the model complexity and approximation power.

It is natural to ask whether the black-box DNN can produce the right solution we desire. We give a positive answer to this question by providing theoretical guarantees on bounding the excess risk of the DNN solution in terms of the network structure and sample size in a non-asymptotic manner. Here, the excess risk refers to the discrepancy in the expected newsvendor loss between any decision and the optimal one when knowing the demand distribution. In particular, we discover that the excess risk can be separated into two independent components. The first component is the \emph{stochastic error} due to random realizations in the observed data, while the second one is the \emph{approximation error} capturing the distance between the target quantile function and the function class represented by DNNs. If the given DNN can directly represent the target quantile function, then only the stochastic error is present, which encompasses parametric (piecewise) linear models as special cases. In such a situation, further increasing the network complexity is unnecessary, as doing so would only amplify the stochastic error bound. 
In general situations, the approximation error decreases with network complexity. Therefore, choosing appropriate network configurations to strike a balance between stochastic and approximation errors becomes crucial.

We analyze both the high-probability and expected excess risk bounds. Specifically, for the former, the stochastic error bound scales as $1/\sqrt{n}$ up to a logarithm factor of the sample size $n$, and as $\sqrt{\log(1/\delta)}$ for $\delta$, where $1-\delta$ denotes the probabilistic accuracy of the high-probability bound. After balancing between the stochastic and approximation errors, the excess risk scales as $n^{-\frac{\beta}{(2\beta+p)}}$ for $n$ and also as $\sqrt{\log(1/\delta)}$ for $\delta$, where $\beta$ measures the smoothness of the underlying target function and $p$ is the dimension of feature variables. This implies that the convergence rate is faster when the target function is smoother or when the dimension of features is lower. More importantly, we \emph{explicitly} characterize the impact of each model parameter, such as the network width, depth, and cost parameters, on the risk bound. {\color{black}Building on this, practitioners can explore suitable network designs while gaining insight into how the properties of the target function influence the required sample size and network complexity. Detailed guidance on this is provided in Section \ref{sec:simulation_implementation}. 
We further illustrate that the DNN method can achieve the minimax optimal rate in terms of the expected excess risk.} This is done by showing that the excess risk can be bounded both from above and below by polynomial terms with identical orders in the sample size up to logarithmic factors. The obtained convergence rate of the expected excess risk bound is sharper than that of the high-probability bound due to technical reasons. Nevertheless, they are both state-of-the-art when respectively compared with the existing results \citep{ban2019big,schmidt2020nonparametric}. 
With the above results, we justify the applicability of DNNs in our data-driven newsvendor problem and anticipate that DNN models hold significant promise in addressing other OM problems.

On the technical side, we resort to the Lipschitz property of the Newsvendor loss function to decompose the excess risk into stochastic and approximation errors. The analysis of the stochastic error relies on empirical process theory \citep{anthony1999neural,bartlett2019nearly}, where concentration inequalities play a critical role in connecting the stochastic error to complexity measures of DNN function classes. We further exploit the results in \cite{jiao2023deep} regarding approximation theory for DNNs to establish non-asymptotic error bounds. The lower bound for the expected excess risk can be established by formulating an equivalent multiple hypothesis testing problem and applying Fano's inequality \citep{scarlett2019introductory}. In particular, by leveraging the new approximation theory developed in \cite{jiao2023deep}, all these bounds do not require special distributional conditions or network restrictions, differing themselves from those found in the literature (e.g., \citealt{schmidt2020nonparametric}).


We further consider two practically critical extensions that enhance our theoretical contributions and the applicability of our framework. First, we provide the excess risk bound of the DNN solution when the data is not independent and identically distributed (i.i.d.). This extension is relevant in practical applications as we usually deal with time series data where the feature information and past demand are collected over a certain timeline. Without any restriction on the underlying data-generating stochastic process, such as stationarity, we reveal that a similar excess risk bound as in the i.i.d. case still holds, with an additional term measuring the intrinsic discrepancy of the underlying data-generating process. As the data-generating process exhibits greater stationarity, this discrepancy diminishes, and the excess risk bound converges to the counterpart in the i.i.d. case. Our findings demonstrate that even for the weakly dependent and non-stationary data, the DNN method can still yield satisfactory results. This observation extends the existing results in the DNN theory.

Second, we explore the scenarios where the curse of dimensionality can be alleviated when using the DNN method. Notice that the convergence rate of the excess risk bound is inversely proportional to the dimension of feature variables $p$. When $p$ is sufficiently large, the convergence can be extremely slow, requiring a substantial amount of data to obtain a good prediction. Nevertheless, we find that if the target function possesses a favorable composite structure, the convergence rate is affected by the intrinsic dimension of the composite function rather than the ambient dimension $p$. Specifically, when a generalized linear model is assumed for the target function, the intrinsic dimension is one, and we can establish an excess risk bound that scales as $n^{-\frac{\beta}{2\beta+1}}$, which is considerably faster than the original rate, particularly when $p$ is large.

{\color{black} Related to our theoretical results, we conduct simulation studies to provide practical insights for implementation. A central argument of our theory is the need to carefully balance the stochastic and approximation errors with opposite behaviors with respect to network complexity, given a fixed amount of data. Echoing this, we observe that the excess risk does not decrease monotonically with either an increasing network depth or width; indeed, beyond a certain threshold, the excess risk turns to rise if we use a wider or deeper network. In addition, we find that a range of network configurations can achieve comparably good performance, which not only indicates the flexibility of the DNN method but also aligns with our theory regarding selecting proper network designs to attain risk minimization. In practice, we recommend starting with a sufficiently large network width and gradually increasing the depth through fine-tuning to find a good candidate network. Of course, the specific network design will depend on various factors such as the dimensionality and smoothness of the target function. We also examine the convergence of excess risk in the sample size, which is found to be polynomial with statistical significance. Consistent with our error analysis, the convergence rate is slower than that of the stochastic error of order $1$, while the difference is attributed to the presence of the approximation error.}

Finally, using a unique dataset from our collaborated food supermarket in China, we conduct a case study to assess the practical performance of the DNN method and compare it with several typical data-driven approaches in the literature. To address the possible concern that the DNN method may require massive data to produce a satisfactory result, we first apply it to a subset of the full dataset where the sample size is relatively small and the number of feature variables is also moderate. We then evaluate the DNN method's performance on the complete dataset. Both sets of experiments yield consistent results,  demonstrating that the DNN method outperforms other alternatives across a wide range of critical levels $\rho$ (defined in \eqref{newsvendor_quantile}). For instance, when $\rho=0.65$, the best alternative (kernel optimization) still incurs a newsvendor loss that is 11.94\% greater than the DNN method. Furthermore, thanks to the well-established computational modules, such as \emph{Adam} \citep{kingma2014adam}, the training and execution processes of DNNs are convenient.


\section{Literature Review}\label{sec:literature}

Echoing the proposals put forth in  \cite{feng2018research,feng2022developing}, our work contributes to the advancement of the theory and practice of OM problems using state-of-the-art data-driven methods. In particular, this paper 
relates to three streams of the existing literature. 

The first stream applies advanced machine learning techniques in OM. For example, \cite{oroojlooyjadid2020applying} experiment using DNN algorithms to simultaneously optimize order quantities for a list of products based on features. \cite{qi2023practical} propose a deep learning framework to address the multi-period inventory replenishment problem. \cite{gijsbrechts2022can} and \cite{oroojlooyjadid2022deep} explore the application of deep reinforcement learning in solving inventory problems. \cite{ye2023deep} combine deep learning and doubly robust estimation to conduct causal inference for large-scale experiments. \cite{chan2022machine} employ a machine learning model to estimate objective values for out-of-sample clients based on in-sample data and solve real-world cycling infrastructure planning problems. More recently, \cite{gabel2022product} and \cite{aouad2022representing} have also incorporated neural networks into product choice modeling and assortment scheduling.

Of course, we are not among the first to embrace machine learning advances such as DNN in the OM community, but we exploit the special structure of the newsvendor problem to relate the optimization problem to the statistical quantile regression. To this end, we are able to establish theoretical guarantees that justify the DNN approach. Notably, we explicitly explain the impact of relevant parameters, such as network structures and the feature dimension, on the potential generalization error. Such detailed error analysis regarding complex neural networks is usually absent in this stream of literature, where it is typically limited to relatively simple parametric models. Moreover, we apply the DNN method on a unique real-world dataset and confirm its practical effectiveness and superiority over other alternative data-driven approaches, further complimenting previous empirical experiments as in \cite{oroojlooyjadid2020applying}.

Our work also closely relates to the second stream of literature on data-driven prescriptive analytics in inventory management. To name a few, \cite{levi2007provably} describe a sampling-based method to reach an inventory decision that achieves near-optimal performance, while \cite{levi2015data} improve the analytical bounds therein. \cite{lin2022data} further relax their assumptions to prove upper bounds for the general and worse-case regrets. \cite{qi2022distributionally} examine the conditional quantile prediction when the data is not identically distributed. 
\cite{zhang2023optimal} study a distributionally robust optimal policy that can generalize unseen feature values well. Other commonly used approaches in this field also include stochastic gradient descent algorithms and bandit controls, see, e.g., \cite{shi2016nonparametric,zhang2018perishable,simchi2022bypassing,simchi2023phase}. 
In a general optimization framework beyond inventory management, \cite{bertsimas2020predictive} adopt a nonparametric method to approximate the underlying target distribution through a weighted empirical distribution and then solve a conditional stochastic optimization problem. \cite{kallus2023stochastic} construct random forests based on a learning scheme to solve for the optimal decision and illustrate their effectiveness. \cite{bertsimas2022data} instead use the reproducing kernel Hilbert space to find the decision from covariates.

In this stream of literature, the closest work to ours is \cite{ban2019big}, which solves a similar data-driven feature-based newsvendor problem using linear and kernel decision rules. 
{\color{black} Their linear policy, which is convenient to implement, accommodates both basic and nonlinearly transformed features, making it highly adaptable. Meanwhile, the developed kernel method, despite being nonlinear, offers an elegant solution that can be expressed analytically and solved efficiently via a simple ranking algorithm. The authors also provide finite-sample performance bounds for out-of-sample costs, depending on the dimensionality and sample size. In contrast, we employ a completely different DNN method to address the problem. Inspired by the seminal work of  \cite{ban2019big}, we establish non-asymptotic performance bounds for the excess risk, which accounts for the approximation error introduced when the target function lies outside the optimization function class. In addition to high-probability bounds, we verify that the derived rate is minimax optimal for the expected excess risk—an aspect not covered by \cite{ban2019big}.
}


Our DNN solution also differs from the two-step predict(or estimate)-then-optimize method (see, e.g., \citealt{hu2022fast,chen2022using,perakis2022optimizing}) in that the DNN model directly outputs the optimal decisions. In fact, as noticed in \cite{liyanage2005practical} and \cite{ban2019big}, the two-step approach may lead to amplified errors if the first-step estimation model is misspecified. More recently, \cite{siegel2021profit,siegel2023data} identify the statistical estimation error for the expected profit in data-driven newsvendor models and explore how to correct it asymptotically. {\color{black}Nevertheless, the DNN method offers an additional advantage over other approaches by mitigating model misspecification errors through the use of a sufficiently large network.} Moreover, the general optimization framework laid out in some previously mentioned papers, such as \cite{kallus2023stochastic} and \cite{chen2022using}, is not specifically tailored to newsvendor analysis, which leaves much room for a more delicate treatment by exploiting the newsvendor problem's special structure. Our theoretical and empirical results complement those related papers in this regard.

Third, our paper contributes to statistical learning and high-dimensional quantile regression theories regarding DNNs. On the one hand, DNNs are typical statistical learning vehicles that study how to summarize and generalize useful information from empirical processes \citep{anthony1999neural,mohri2018foundations,bartlett2019nearly}. On the other hand, due to the similarity in the loss function between the newsvendor problem and quantile regression, our work also contributes to the studies in the latter using DNNs.  Recently, there have been many relevant papers dedicated to establishing generalization bounds for DNNs in the context of least square or quantile regression (e.g., \citealt{schmidt2020nonparametric,farrell2021deep,padilla2022quantile}). Compared to these papers, we explicitly characterize the prefactor of excess risk bounds as a polynomial function of the dimensionality $p$, significantly improving upon the extant results with an exponential prefactor. 
{\color{black}We also do not require the related restrictions on the network structure, such as the parameters should be bounded by a prescribed constant and the whole network has to be sparse in a certain level \citep{schmidt2020nonparametric,padilla2022quantile}. 
Finally, our paper also differs from \cite{jiao2023deep} and \cite{shen2021deep} in that we focus primarily on high-probability and expectation bounds with fast convergence rates corresponding to the newsvendor loss using different techniques. Specifically, \cite{jiao2023deep} examine the deep least square regression problem, where the square loss function benefits from the variance-bounding trick, linking the variance of the loss difference to its expectation. 
However, this approach does not directly apply to the newsvendor loss, requiring us to use conditional expectations to derive new generalization bounds for the empirical process (see Theorem \ref{thm:stochastic_error2}).  \cite{shen2021deep} investigate a relevant quantile regression problem, but their primary technical contribution lies in dealing with possibly unbounded response variables through truncation techniques, where the final rate depends on moment conditions. In contrast to both papers, we further validate that the obtained rate for the expected excess risk bound is indeed minimax optimal, using Fano's inequality \citep{scarlett2019introductory}. 
While those two papers focus on dimensionality challenges of DNN estimation with i.i.d. data, we are mainly interested in newsvendor-related issues, such as comparisons with the classical SAA or linear models, handling dependent data, and providing practical guidance alongside real-world applications.}

\section{Deep Neural Newsvendor with Feature Data}\label{sec:problem_formulation}

\subsection{Feature-Based Newsvendor Problem}

We consider a data-driven newsvendor problem where the firm needs to decide the inventory level of perishable goods to serve future random demand. Specifically, the demand distribution is unknown to the firm, but the historical demand realizations $\bd(n)=(d_1,d_2,\ldots,d_n)^\top$ are accessible. Throughout the paper, we denote vectors or matrices using bold symbols. Besides the demand data, we assume that the firm is also able to collect feature data $\bx(n)=(\bx_1,\bx_2,\ldots,\bx_n)^\top$ consisting of observable (and not necessarily independent) features along with each demand realization. These features may include, among others, macroeconomic data such as GDP and CPI, weather conditions, and even location details \citep{ban2019big}. In practice, it is reasonable to imagine that the demand is influenced by various factors, and we are interested in using these factors to predict the demand. In this paper, we assume that the features have already been carefully and properly selected based on their significant relevance to the demand, and their data are available prior to decision-making. The variable selection problem falls beyond the scope of the current study.

Let $D\in\mathbb{R}$ and $\bX\in\mathbb{R}^p$ denote the random demand and feature vector, respectively. The classical newsvendor problem aims to solve for the optimal ordering quantity that minimizes the expected newsvendor cost loss expressed as
\begin{equation}\label{newsvendor_cost_function}
    C(q):=\mathbb{E}_D[b(D-q)^++h(q-D)^+],
\end{equation}
where the expectation $\mathbb{E}_D$ is taken with respect to $D$, $q$ is the inventory decision, $b$ (resp., $h$) represents the unit underage (resp., overage) cost for the unsatisfied demand (resp., unsold inventory). Let $F(\cdot)$ denote the cumulative distribution function (CDF) of the random demand $D$. In an oracle case with full information of $F(\cdot)$, the optimal quantity is then
\begin{equation}\label{newsvendor_quantile}
    q^*=\inf\{q\geq 0: F(q)\geq \rho\}, \quad \text{where } \rho:=\frac{b}{b+h}.
\end{equation}
Obviously, the above oracle newsvendor solution cannot be directly obtained in practice since we do not know $F(\cdot)$ in advance. Data-driven methods then come into play by using the empirical distribution of random samples to estimate the unknown $F(\cdot)$, and a typical example is the celebrated SAA approach (see \citealt{shapiro2009lectures} for more details).  

Alternatively, the feature-based newsvendor problem aims to find an optimal measurable mapping $f^*: \mathbb{R}^p\rightarrow \mathbb{R}$ such that, based on the observed features $\bX=\bx$, the optimal inventory decision is determined accordingly by $f^*(\bx)$. In other words, 
\begin{equation*}
    f^*:=
    \underset{f}{\argmin}\;  \mathbb{E}_{D,\bX}[b(D-f(\bX))^++h(f(\bX)-D)^+],
\end{equation*}
where the expectation $\mathbb{E}_{D,\bX}$ is taken with respect to the joint distribution of $D$ and $\bX$. Equivalently, for any given $\bX=\bx$, the optimal decision can be established with respect to the conditional distribution $D|\bX$ as follows
\begin{equation}\label{feature_newsvendor}
    f^*(\bx):=
    \underset{f(\bx): f{\rm\ measurable}}{\argmin}  \mathbb{E}_{D\mid \bX}[b(D-f(\bX))^++h(f(\bX)-D)^+\mid \bX=\bx].
\end{equation}
Similar to \eqref{newsvendor_quantile}, the oracle solution to \eqref{feature_newsvendor}, when knowing the joint distribution of $\bZ:=(\bX,D)$, is exactly the conditional quantile function given by
\begin{equation}\label{true_conditional_quantile}
    f_{\rho}(\bx):=\inf\{q\geq 0: F(q\vert \bX=\bx)\geq \rho\},
\end{equation}
where $F(\cdot\vert \bx)$ is the conditional CDF of $D$ given $\bX=\bx$.

With only finite random samples $\bS_n=\{(\bx_i,d_i)\}_{i=1}^n$ from the joint distribution $\bZ$ in hand, we then consider the data-driven, feature-based newsvendor problem by minimizing the empirical risk $\mathcal{R}^{\rho}_n(\cdot)$ to obtain the empirical risk minimizer (ERM) $\hat{f}_{n}$ by:
\begin{equation}\label{data_driven_newsvendor}
    \hat{f}_{n}\in \underset{f\in\mathcal{F}_n}{\argmin} \mathcal{R}^{\rho}_n(f):= \underset{f\in\mathcal{F}_n}{\argmin} \frac{1}{n}\sum_{i=1}^n\bigg[b(d_i-f(\bx_i))^++h(f(\bx_i)-d_i)^+\bigg],
\end{equation}
where $\mathcal{F}_n$ is a suitable function class, which may depend on the sample size $n$. 
In this paper, we choose $\mathcal{F}_n$ to be a class of functions represented by DNNs. In addition to exploring the practical performance of DNNs for our problem, we are particularly interested in establishing theoretical guarantees so that the DNN tool does not simply work as a black box. Our nonparametric DNN solution integrates the forecasting and optimization steps into one, unlike the separated estimation and optimization (SEO) approach. Moreover, we do not require prior knowledge or assumptions about the underlying demand distribution.

\subsection{Newsvendor Solution Using DNNs}\label{subsec:DNN_model}

\begin{figure}[!ht]
	\centering
 \includegraphics[scale=0.636]{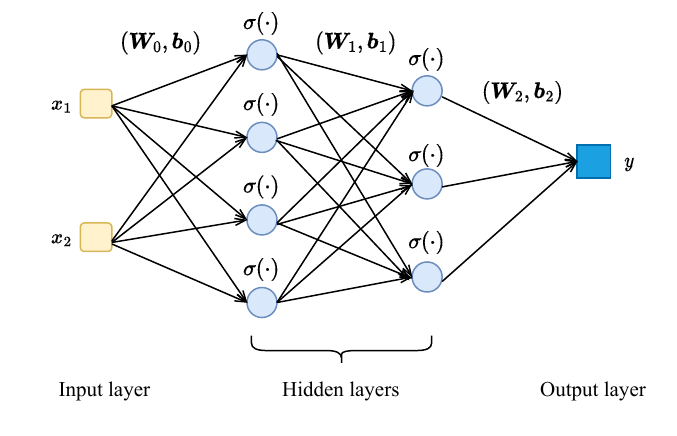}	
	\caption{\color{black}An illustrative fully connected feedforward neural network with two hidden layers.}
		\label{fig:NN_structure}
\end{figure}

Typical multi-layer perceptrons (MLPs), the commonly used feedforward neural networks, consist of several key components: the number of layers that decides the neural network's depth, the number of neurons in each layer that determines the network width, the connections between neurons specifying how data is processed through the network, and the activation functions in each neuron that introduce nonlinearity. Figure \ref{fig:NN_structure} illustrates the architecture of a two-hidden-layer network that maps the input $(x_1,x_2)$ to the output $y$. 

Mathematically, we can also represent the MLP architecture in terms of the composition of a series of functions as follows:
\begin{equation}\label{network_functional_composition}
    f(\bx;\bW,\bb)=\mathcal{L}_\mathcal{D}\circ \sigma \circ \mathcal{L}_{\mathcal{D}-1}\circ \sigma \circ \cdots \circ \sigma \circ \mathcal{L}_1\circ\sigma \circ\mathcal{L}_0(\bx), \quad \bx\in\mathbb{R}^p,
\end{equation}
where $\mathcal{D}$ is the number of hidden layers. 
Throughout the paper, we consider the standard component-wise rectified linear unit (ReLU) activation function $\sigma(x)=\max\{0,x\}$ and 
$\mathcal{L}_i(\bx)=\boldsymbol{W_i}\bx+\bb_i,$ for $i=0,1,\ldots,\mathcal{D},$ in which $\boldsymbol{W_i}\in\mathbb{R}^{w_{i+1}\times w_{i}}, \bb_i\in\mathbb{R}^{w_{i+1}}$, and $w_i$ is the number of neurons (width) of the $i$-th layer. Consider the network shown in Figure \ref{fig:NN_structure} as an example. The input $(x_1,x_2)$ first undergoes a linear transformation with weights $\bW_0$ and bias $\bb_0$, 
resulting in a 4-dimensional vector. The ReLU activation function is then applied component-wise in the first hidden layer. This process is repeated for the subsequent hidden layer. Finally, the output $y$ is obtained by linearly combining the activated values from the last hidden layer with weights $\bW_2$ and bias $\bb_2$.

{\color{black}
For a given neural network, its width, denoted by $\mathcal{W}$, is defined as the maximum number of neurons in any hidden layer, i.e., $\mathcal{W}:=\max\{w_1,\ldots,w_{\mathcal{D}}\}$; its depth $\mathcal{D}$ refers to the number of hidden layers; its size $\mathcal{S}$ is the total number of parameters (weights and biases) involved, given by $\mathcal{S}:=\sum_{i=0}^{\mathcal{D}}w_{i+1}(w_i+1)$; and the number of neurons, denoted by $\mathcal{U}$, counts the neurons across all hidden layers, i.e., $\mathcal{U}:=\sum_{i=1}^{\mathcal{D}}w_i$. We then define the optimization function class $\mathcal{F}_n$ as $\mathcal{F}_{\mathcal{D},\mathcal{W},\mathcal{U},\mathcal{S},\mathcal{B}}$, which represents a class of feedforward neural network mappings $f_{\phi}:\mathbb{R}^p\rightarrow \mathbb{R}$, parameterized by $\phi$, with the depth $\mathcal{D}$, width $\mathcal{W}$, size $\mathcal{S}$, and number of neurons $\mathcal{U}$, such that $f_\phi$ satisfies $\|f_\phi\|_{\infty}\leq\mathcal{B}$ for some positive constant $\mathcal{B}$, where $\|\cdot\|_\infty$ denotes the sup-norm. Notably, we do not require all layers to have the same width, as imposed in \cite{chen2022using}, nor do we need to stipulate that the network is sparse or that the parameters are uniformly bounded by 1 (see, e.g., \citealt{schmidt2020nonparametric,padilla2022quantile}). In what follows, we write $\mathcal{F}_{\mathcal{D},\mathcal{W},\mathcal{U},\mathcal{S},\mathcal{B}}$ as $\mathcal{F}_{D\!N\!N}$ for brevity. Additionally, we define a DNN as any neural network with $\mathcal{D}\geq 1$ hidden layers.

Due to the special structure of ReLU networks in \eqref{network_functional_composition}, their outputs are piecewise linear functions. Although alternative activation functions, such as Sigmoid and Tanh, are available, ReLU networks enjoy advantages in terms of mathematical simplicity with nice theoretical properties, such as the universal approximation property \citep{hornik1991approximation}. {\color{black}Within the context of the considered newsvendor problem, we can further demonstrate that ReLU networks attain the minimax optimal convergence rate with respect to the sample size (see Theorems \ref{thm:upper_expected_bound} and \ref{thm:lower_expected_bound}). This result is achieved because the newsvendor loss is both Lipschitz continuous and convex in the inventory decision, and the optimal decision amounts to estimating a specific quantile of a conditional demand distribution.} 
Meanwhile, they have gained wide popularity in deep learning applications due to their computational efficiency, ease of training, and strong generalization performance (see, e.g., \citealt{nair2010rectified, krizhevsky2012imagenet}). This success can be attributed, in part, to practical factors such as effective initialization strategies \citep{glorot2010understanding,he2015delving}, which help mitigate issues like vanishing or exploding gradients in deep architectures. Given their theoretical strengths and empirical effectiveness, we focus on ReLU networks as our optimization function space.

}

Our feature-based, data-driven deep neural newsvendor problem aims to find a measurable function $\hat{f}_{D\!N\!N}\in\mathcal{F}_{D\!N\!N}$ such that 
\begin{equation}\label{DNN_solution_definition}
    \hat{f}_{D\!N\!N}:= \underset{f\in\mathcal{F}_{D\!N\!N}}{\argmin}\,   \mathcal{R}^{\rho}_n(f),
\end{equation}
where $\mathcal{R}_n^{\rho}(f)$ is the empirical newsvendor risk defined in \eqref{data_driven_newsvendor}. To evaluate the theoretical performance of the deep neural newsvendor solution, we employ the oracle solution $f_{\rho}$ in \eqref{true_conditional_quantile} as a benchmark and quantify the \emph{excess risk} bounds in terms of the sample size and network parameters in the next section. Specifically, the \emph{excess risk} is defined as the difference in the expected newsvendor loss between $\hat{f}_{D\!N\!N}$ and $f_{\rho}$ based on the joint distribution of $\bZ=(\bX,D)$, i.e.,
\begin{equation*}
    \mathcal{R}^{\rho}(\hat{f}_{D\!N\!N})-\mathcal{R}^{\rho}(f_{\rho}):=\mathbb{E}_{\bZ}[b(D-\hat{f}_{D\!N\!N}(\bX))^++h(\hat{f}_{D\!N\!N}(\bX)-D)^+]-\mathbb{E}_{\bZ}[b(D-f_{\rho}(\bX))^++h(f_{\rho}(\bX)-D)^+].
\end{equation*}
We notice that the excess risk remains a random variable, contingent upon the random samples $\bS_n$ and through the fitting of $\hat{f}_{D\!N\!N}$.

\section{Theoretical Guarantees: Non-Asymptotic Excess Risk Bounds}\label{sec:theoretical_guarantees}

{\color{black}Machine learning methods using neural networks are often perceived as black boxes, where outputs are produced by feeding covariates into trained networks. They have demonstrated remarkable practical performance in various tasks, e.g., image classification \citep{krizhevsky2012imagenet,he2015delving}, game intelligence \citep{silver2016mastering}, and speech recognition \citep{hinton2012deep}, given that the data amount is ``adequate." Despite these successes, several questions naturally arise: How does the performance improve as more data is available? How does the network structure, such as the network depth and width, influence prediction accuracy? Most importantly, can DNNs solve the target problem with theoretical guarantees rather than relying solely on trial and error?} In this section, we address these questions in the context of the data-driven newsvendor problem by examining the theoretical non-asymptotic excess risk bounds for DNN solutions in terms of the sample size and network structure. 

\subsection{High-Probability Bound for Excess Risk}\label{subsec:probability_bound}

In general, the excess risk arises from two different sources: the \emph{stochastic error}, resulting from the random realizations of observed data, and the \emph{approximation error}, caused by the inability of the DNN with a given network structure to exactly represent the target quantile function. The former risk is inherent in data-driven problems since we aim to apply the learned decision rule from finite observed data to any future random realization. The latter one emerges when seeking the optimal rule within a prescribed set $\mathcal{F}_{D\!N\!N}$, rather than the general function space. However, it is expected that with more complex network structures, DNNs can represent a broader class of functions, reducing this error to a negligible level. We formally state these intuitions in the following lemma.

\begin{lemma}\label{lemma:risk_decomposition}
    The excess risk of the ERM $\hat{f}_{D\!N\!N}$ satisfies
    \begin{equation*}
        \mathcal{R}^{\rho}(\hat{f}_{D\!N\!N})-\mathcal{R}^{\rho}(f_{\rho})\leq \underbrace{2\underset{f\in\mathcal{F}_{D\!N\!N}}{\sup}|\mathcal{R}^{\rho}(f)-\mathcal{R}_n^{\rho}(f)|}_\text{stochastic error}+\underbrace{\underset{f\in\mathcal{F}_{D\!N\!N}}{\inf}\left\{\mathcal{R}^{\rho}(f)-\mathcal{R}^{\rho}(f_{\rho})\right\}}_\text{approximation error}.
    \end{equation*}
\end{lemma}
As elucidated in Sections \ref{subsubsec:stochastic_error} and \ref{subsec:approximation_balanced_errors}, on the one hand, a more complex network can result in a smaller approximation error because of the growing representation power; on the other hand, it will amplify the stochastic error due to the increased number of model parameters needing to be estimated. 
We first separately analyze each error and then strike a balance between them to achieve an optimal error rate with respect to the sample size.

\subsubsection{Stochastic Error: When the Fixed DNN Can Represent $f_{\rho}$.}\label{subsubsec:stochastic_error}

According to Lemma \ref{lemma:risk_decomposition}, if the given DNN can represent $f_{\rho}$, i.e., $f_{\rho}\in\mathcal{F}_{D\!N\!N}$, the approximation error term vanishes, and only the stochastic error remains when using this given DNN to solve for the newsvendor decision. {\color{black}However, it is important to note that this scenario is an ideal one, and we typically do not know the exact form of the underlying conditional quantile function $f_{\rho}$ or whether it belongs to the representation class of a certain DNN. 
Nevertheless, we will address the general case later in Theorem \ref{thm:excess_risk}, knowing only the target function's smoothness rather than its specific functional form.}

We make the following assumption regarding the underlying distributions.

\begin{assumption}\label{assumption_1}
 The random demand $D$ has a bounded compact support, which is assumed to be $[0,\Bar{D}]$ for some $\Bar{D}>0$ without loss of generality.
\end{assumption}

Assumption \ref{assumption_1} is very mild and reasonable since the demand data should be finite in values in practice. We can derive a stochastic error bound based on this assumption.

\begin{theorem}[{\sc Stochastic Error Bound}]\label{thm:stochastic_error}
    Under Assumption \ref{assumption_1}, for any $\delta\in(0,1)$ and $n\geq C\cdot \mathcal{S}\mathcal{D}\log(\mathcal{S})$ with a large enough $C>0$, with probability at least $1-\delta$ over the random draw of $\bS_n=\{(\bx_i,d_i)\}_{i=1}^n$, where $(\bx_i,d_i)$ are i.i.d. samples from the unknown joint distribution $(\bX,D)$, we have
    \begin{equation*}
        \underset{f\in\mathcal{F}_{D\!N\!N}}{\sup}|\mathcal{R}^{\rho}(f)-\mathcal{R}_n^{\rho}(f)|\leq \sqrt{2}(b\bar{D}+h\mathcal{B})\bigg(C_1\sqrt{\frac{\mathcal{S}\mathcal{D}\log(\mathcal{S})\log(n)}{n}}+\sqrt{\frac{\log(1/\delta)}{n}}\bigg),
    \end{equation*}
    where $C_1>0$ is an independent universal constant.
\end{theorem}

{\color{black}Theorem \ref{thm:stochastic_error} is established using empirical process theory \citep{anthony1999neural,bartlett2019nearly}, and the error bound does not require restrictions on the DNN structure.  
As a direct consequence of Lemma \ref{lemma:risk_decomposition}, we obtain the excess risk bound under the assumption that $f_{\rho}\in\mathcal{F}_{D\!N\!N}$.}

{\color{black}\begin{corollary}\label{corollary:stochastic_error}
    If $f_\rho\in\mathcal{F}_{D\!N\!N}$, under conditions of Theorem \ref{thm:stochastic_error}, with probability at least $1-\delta$ over the random draw of $\bS_n$,
    \begin{equation*}
        \mathcal{R}^{\rho}(\hat{f}_{D\!N\!N})-\mathcal{R}^{\rho}(f_{\rho})\leq 2\sqrt{2}(b\bar{D}+h\mathcal{B})\bigg(C_1\sqrt{\frac{\mathcal{S}\mathcal{D}\log(\mathcal{S})\log(n)}{n}}+\sqrt{\frac{\log(1/\delta)}{n}}\bigg).
    \end{equation*}
\end{corollary}}

When $f_{\rho}\in\mathcal{F}_{D\!N\!N}$, Theorem \ref{thm:stochastic_error}, and in particular, Corollary \ref{corollary:stochastic_error}, describe the additional cost of the in-sample learned inventory policy for the out-of-sample newsvendor problem compared to the optimal oracle solution. It is reasonable to imagine that with more collected samples, the DNN can better learn the underlying distributional structure so that its performance gets better. Specifically, we notice that with a fixed DNN and given cost parameters $b$ and $h$, the stochastic error (excess risk) scales in the order of $O(\sqrt{\log(n)/n})$.
Meanwhile, the stochastic error bound increases with the network complexity. {\color{black}In other words, if the target function can be exactly represented by a fixed neural network (e.g., assuming it is linear as in Example \ref{example:linear_model}), there is no need for a more complex and redundant configuration with additional width or depth. Of course, this assumption is not guaranteed in practice (hence referred to as the ideal case), but we will get rid of this impractical condition for the general case in the next subsection.}

\begin{example}[Linear Demand Model]\label{example:linear_model}
    We consider the linear demand model with random errors
    $D|(\bX=\bx)=\alpha+\bbeta^\top\bx+\epsilon,$ 
    where $\alpha\in\mathbb{R}, \bbeta\in\mathbb{R}^p$, and $\epsilon\sim F_{\epsilon}$ is independent of the random feature vector $\bX$, with its $\rho$-th quantile being zero. Therefore, the conditional $\rho$-th quantile function is given by
    $
        f_\rho(\bx)=\alpha+\bbeta^\top\bx
    $, which implies that the optimal oracle inventory policy is simply a linear combination of the feature variables. In this case, the following lemma shows that a one-layer ReLU neural network suffices to represent it.
    {\color{black}\begin{lemma}\label{lemma:linear_functions_nn}
        The linear function $f_\rho(\bx)=\alpha+\bbeta^\top\bx$ can be represented by the one-layer ReLU neural network architecture \eqref{network_functional_composition} via
        $f_{\rho}(\bx)=\alpha+\bbeta^\top\bx=\bW_1\sigma(\bW_0\bx+\bb_0)+b_1,$ 
        where $\bW_0=(\bbeta,-\bbeta)^\top\in\mathbb{R}^{2\times p}$, $\bb_0=(\alpha,-\alpha)^\top\in\mathbb{R}^{2}$, $\bW_1^\top=(1,-1)^\top\in\mathbb{R}^{2}$, and $b_1=0$.  \end{lemma}

        \begin{figure}[!ht]
	\centering
 \includegraphics[scale=0.396]{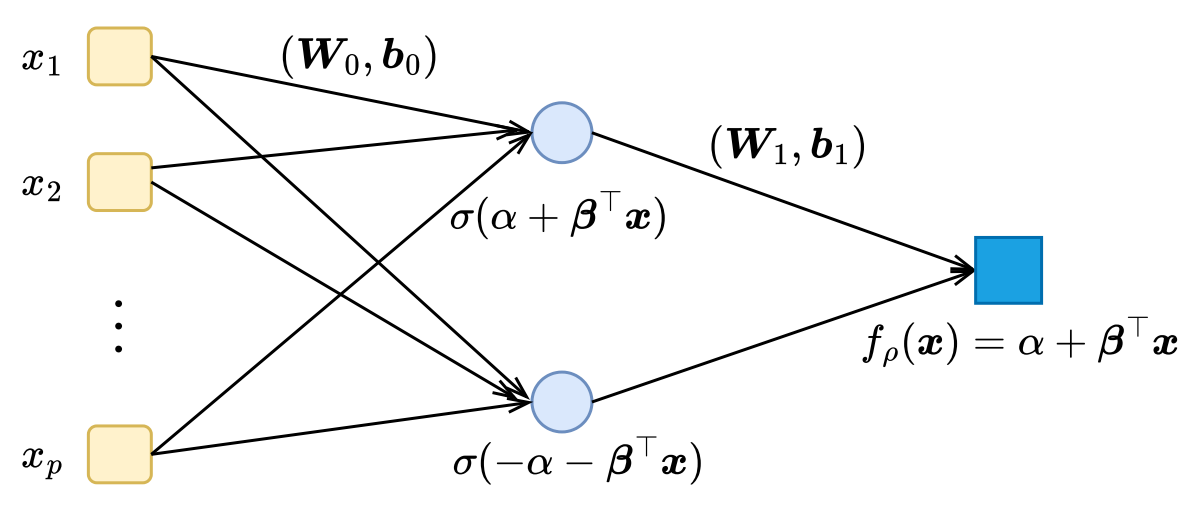}	
	\caption{Illustration of a one-layer ReLU neural network representing linear functions.}
		\label{fig:linear_structure}
\end{figure}

  Figure \ref{fig:linear_structure} illustrates the one-layer ReLU neural network. In this parametric setting, we can directly apply Corollary \ref{corollary:stochastic_error} to obtain an excess risk bound that scales in the order of $O(\sqrt{p\log(n)/n})$, given that $\mathcal{S}=\sum_{i=0}^1w_{i+1}(w_i+1)=2p+5$ for the network structure shown in Figure \ref{fig:linear_structure}. For this special case, our bound is comparable, up to a logarithm factor, to the counterparts in \cite{ban2019big}, whether using linear or kernel decision rules. Nevertheless, beyond parametric models such as in Example \ref{example:linear_model}, we are able to provide theoretical guarantees for general target functions without requiring a priori specification of their explicit functional forms.}
    \end{example}

\subsubsection{Approximation Error and Balanced Excess Risk  Bounds.}\label{subsec:approximation_balanced_errors}

The approximation error quantifies the ability of any DNN to approximate the target function, regardless of the sample size. Indeed, it is determined by the smoothness and dimensionality of the target function. For instance, if the target function is simple, say a linear function, as in Example \ref{example:linear_model}, the approximation error could be extremely small. To encompass a wide range of scenarios, we here consider the general case where the target function $f_{\rho}$ belongs to a H\"older function class.

\begin{definition}[H\"older Function Class]\label{def:holder_class}
    For $\beta,B_0>0$, and $\mathcal{X}\in\mathbb{R}^{p}$, the H\"older function class $\mathcal{H}^\beta(\mathcal{X},B_0)$ is defined by
    \begin{equation*}
        \mathcal{H}^\beta(\mathcal{X},B_0):=\bigg\{f:\mathcal{X}\rightarrow\mathbb{R}, \underset{\|\balpha\|_1\leq s}{\max}\|\partial^\balpha f\|_{\infty}\leq B_0, \text{ and } \underset{\|\balpha\|_1=s}{\max}\underset{x\neq y}{\sup}\frac{|\partial^\balpha f(x)-\partial ^\balpha f(y)|}{\|x-y\|_2^r}\leq B_0\bigg\},
    \end{equation*}
    where $s:=\lfloor\beta\rfloor$ denotes the largest integer that is strictly smaller than $\beta$ and $r=\beta-\lfloor\beta\rfloor\in(0,1]$; $\balpha:=(\alpha_1,\ldots,\alpha_p)^\top$ is a vector of non-negative integers and $\|\balpha\|_1:=\sum_{i=1}^p\alpha_i$; and $\partial^\balpha=\partial^{\alpha_1}\cdots\partial^{\alpha_p}$ denotes the partial derivative operators.
\end{definition}

{\color{black}By definition, a function in $\mathcal{H}^\beta(\mathcal{X},B_0)$ is $\lfloor\beta\rfloor$-times differentiable over $\mathcal{X}$, with all derivatives up to order $\lfloor\beta\rfloor$ uniformly bounded by $B_0$, and its $\lfloor\beta\rfloor$-th derivative is H\"older continuous with exponent $\beta-\lfloor\beta\rfloor$. For example, when $\beta=1$, this class includes Lipschitz continuous functions, and when $\beta=2$, it includes functions with continuous first-order derivatives.}
We need the following assumptions to analytically analyze the approximation error bound in terms of the network structure.

\begin{assumption}\label{assumption_2}
(i) The support $\mathcal{X}$ of the feature vector is a compact set in $\mathbb{R}^p$, and without loss of generality, we assume that $\mathcal{X}:=[0,1]^p$. The marginal probability measure of $\bX$ is absolutely continuous with respect to the Lebesgue measure; 
 (ii) The conditional quantile function $f_{\rho}$ belongs to the H\"older class $\mathcal{H}^{\beta}([0,1]^p,B_0)$ for a given $\beta>0$ and a finite constant $B_0>0$ {\color{black}with $B_0\leq \mathcal{B}$}.
\end{assumption}

{\color{black}Assumption \ref{assumption_2}(i) imposes a mild boundedness condition on the domain of the feature data, while Assumption \ref{assumption_2}(ii) introduces a basic smoothness condition on the target function. Although the smoothness parameter may not be directly inferred in practice, this condition is rather reasonable and general because it encompasses a large subset of continuous functions. For example, by conservatively setting $\beta=1$, Assumption \ref{assumption_2}(ii) implies that the target quantile function is Lipschitz continuous, a condition that already covers many frequently used demand-covariate models, such as Example \ref{example:linear_model}, where we can choose $B_0$ to be the largest absolute value of the elements in the coefficient vector. 
} 

\begin{proposition}\label{prop:approximation_error}
     Let $\mathcal{F}_{D\!N\!N}$ be the ReLU neural network mappings with the  width and depth respectively specified as
    {\color{black}\begin{equation}\label{choice_of_depth_width}
        \mathcal{W}=38(\lfloor\beta\rfloor+1)^2p^{\lfloor\beta\rfloor+1}N\lceil\log_2(8N)\rceil \text{ and }
        \mathcal{D}=21(\lfloor\beta\rfloor+1)^2M\lceil\log_2(8M)\rceil+2,
    \end{equation}
    }where $\lceil a\rceil$ denotes the smallest integer that is no less than $a$.
    Under Assumption \ref{assumption_2} and for any $N, M\in\mathbb{N}_+$, we have
    \begin{equation*}
        \underset{f\in\mathcal{F}_{D\!N\!N}}{\inf}\left[\mathcal{R}^{\rho}(f)-\mathcal{R}^{\rho}(f_{\rho})\right]\leq 18\max\{b,h\}B_0 (\lfloor\beta\rfloor+1)^2p^{\lfloor\beta\rfloor+(\beta\vee 1)/2}(NM)^{-2\beta/p},
    \end{equation*}
    where $a\vee b:=\max\{a,b\}$.
\end{proposition}

As opposed to Theorem \ref{thm:stochastic_error}, Proposition \ref{prop:approximation_error} demonstrates that the approximation error bound is decreasing in the network complexity (measured by its width $\mathcal{W}$ and depth $\mathcal{D}$ via $N$ and $M$). This observation is quite intuitive since a network with a larger width or depth has greater approximation power. Therefore, a careful network structure design is necessary to strike a balance between the stochastic and approximation errors since the excess risk is a sum of them, as shown in Lemma \ref{lemma:risk_decomposition}.
Moreover, the feature dimension $p$ also affects the approximation error, and the approximation to the target function becomes tougher if the number of covariates grows. Finally, the smoothness parameter $\beta$ influences the approximation error because ReLU neural networks essentially use piecewise linear functions to approximate the target function. 

{\color{black}Although the network depth and width can vary by tuning the free parameters $N$ and $M$, we always treat the uniform upper bound $\mathcal{B}$ for the network mappings as a fixed positive constant. This can be enforced by truncating the network output to ensure that $\|f_\phi\|_\infty\leq\mathcal{B}$, even as the network size grows. Technically, this truncation can be implemented using the following composite ReLU operations: $y=-\max\{\mathcal{B}-\max\{x,0\},0\}+\max\{\mathcal{B}-\max\{-x,0\},0\}$, which ensures that for any $x\in\mathbb{R}$, the output $y$ is constrained within $[-\mathcal{B},\mathcal{B}]$.}
{\color{black}Consequently, to achieve an optimal convergence rate of the excess risk bound in the sample size, we can balance the two errors in Theorem \ref{thm:stochastic_error} and Proposition \ref{prop:approximation_error} by appropriately relating $M$, $N$, or both, to the sample size $n$. Using any of these approaches can result in the same optimal convergence rate but would incur different total numbers of parameters in the network.} 

\begin{theorem}[{\sc High-Probability Bound for Excess Risk}]\label{thm:excess_risk}
    {\color{black}Let the network width and depth be defined according to \eqref{choice_of_depth_width} with $NM=\lfloor n^{\frac{p}{2p+4\beta}}\rfloor$.
    Under Assumptions \ref{assumption_1} and \ref{assumption_2}, for all $n\geq C$ with a large enough $C>0$,} 
    with probability at least $1-\delta$ over the random draw of $\bS_n$,  we have
    \begin{equation}\label{excess_risk_bound}
        \mathcal{R}^{\rho}(\hat{f}_{D\!N\!N})-\mathcal{R}^{\rho}(f_{\rho})\leq 2\sqrt{2}(b\bar{D}+h\mathcal{B})\bigg(C_1(\lfloor\beta\rfloor+1)^4p^{\lfloor\beta\rfloor+1}(\log(n))^2n^{-\frac{\beta}{2\beta+p}}+\sqrt{\frac{\log(1/\delta)}{n}}\bigg),
    \end{equation}
    where $C_1>0$ is an independent universal constant.
\end{theorem}

\begin{remark}
    If the excess risk has some nice local quadratic structure (see Assumption \ref{assumption_3}), the approximation rate in Proposition \ref{prop:approximation_error} can be further accelerated (Proposition \ref{prop:approximation_error_new}). As a result, the excess risk bound can also be slightly improved to $O(n^{-\frac{2\beta}{4\beta+p}})$ (Theorem \ref{thm:excess_risk_e_companion}).
\end{remark}

{\color{black}It is worth highlighting that the excess risk bound in Theorem \ref{thm:excess_risk} is established without requiring prior knowledge of the functional form of $f_\rho$ and applies to a large class of target functions.  
Nevertheless, DNNs can still provide satisfactory solutions with explicit theoretical guarantees in this general setting. The DNN method itself is nonparametric, meaning that we simply need to train a DNN with sufficient data and a properly chosen network structure. In fact, Theorem \ref{thm:excess_risk} suggests that a wide range of network configurations can achieve a designed precision, and this flexibility is also observed numerically (see Section \ref{sec:simulation_implementation}).
In some sense, we ``unwrap" the black box of DNNs in the data-driven newsvendor context by demonstrating that the DNN solution can be reliable in theory and should be considered a valuable candidate in the problem-solving toolbox, particularly for complex scenarios where parametric methods fall short.}

The excess risk bound holds without requiring specific distributional conditions or network restrictions, exhibiting a polynomial dependence on the feature dimension $p$.
It scales as $O(n^{-\frac{\beta}{2\beta+p}})$ up to a logarithmic factor, which comes from the first error term in \eqref{excess_risk_bound}. This indicates that a smoother target function (i.e., $\beta$ is larger) leads to a faster convergence rate. Again, this is because ReLU neural networks are essentially piecewise linear functions, and they excel at approximating smooth functions. However, we also notice that the convergence rate decreases with the feature dimension, so in a high-dimensional scenario, the convergence of the error bound could be very slow. This is a commonly encountered problem in high-dimensional statistics, often referred to as the ``curse of dimensionality." 
To mitigate this issue in the current setting, we propose focusing on the intrinsic dimension of the prediction problem in Section \ref{subsec:dimensionality_issues}. By doing so, the convergence order no longer depends on the feature dimension $p$ but on the much smaller intrinsic dimension if the target function has a composite structure. We also remark that the second error term in \eqref{excess_risk_bound} is a result of the concentration inequalities that bound the deviation between sample values and expected values, which is common in measuring the performance of sample-based policies. 

\subsection{Sharper Bound in Expectation: An Optimal Rate}\label{subsec:lower_bound}

In this section, we show that the DNN method can attain the minimax convergence rate in that the expected excess risk can be bounded both from above and below by the same polynomial term in the sample size, up to logarithmic factors. In this sense, the derived rate is optimal.

{\color{black}\begin{assumption}[Local Quadratic Condition of Excess Risk]
	\label{assumption_3}
	There exist some constants $c^0_\rho=c^0_\rho(\rho,D,\bX)>0$ and $\delta^0_\rho=\delta^0_\rho(\rho,D,\bX)>0$ such that
	$\mathcal{R}^\rho(f) -\mathcal{R}^\rho(f_\rho)\leq c^0_\rho \Vert f-f_\rho \Vert^2_{L^2(\nu)},$
	for any $f$ satisfying $\Vert f-f_\rho \Vert_{L^\infty(\mathcal{X}^0)}\leq \delta^0_\rho$, where $\mathcal{X}^0$ is any subset of $\mathcal{X}$ with $\mathbb{P}(\bX\in\mathcal{X}^0)=\mathbb{P}(\bX\in\mathcal{X})$.
\end{assumption}}

\begin{theorem}[{\sc Upper Bound for Expected Excess Risk}]\label{thm:upper_expected_bound}
    Under Assumptions \ref{assumption_1}, \ref{assumption_2}, and \ref{assumption_3}, with the same network structure as in Theorem \ref{thm:excess_risk}, we have
    \begin{equation}\label{bound_expectation_upper}
        \mathbb{E}\left[\mathcal{R}^{\rho}(\hat{f}_{D\!N\!N})-\mathcal{R}^{\rho}(f_{\rho})\right]
         \leq C(\log(n))^4n^{-\frac{2\beta}{p+2\beta}},
    \end{equation}
    for $n$ large enough, where $C$ is a universal constant depending on model parameters.
\end{theorem}

In fact, the technical assumption of a local quadratic structure for the excess risk is minor, which is satisfied when the density of the conditional distribution $D|\bX$ has an upper bound near $f_\rho(\bX)$. Notably, the bound on the expected excess risk in Theorem \ref{thm:upper_expected_bound} is tighter than the high-probability bound in Theorem \ref{thm:excess_risk} in that the convergence rate with respect to the sample size almost doubles. On the one hand, this is not surprising because one may not be able to directly recover the high-probability bound \eqref{excess_risk_bound}, additive in the logarithm of $\delta$, from \eqref{bound_expectation_upper} via the Markov inequality. The additive form in \eqref{excess_risk_bound} with respect to the logarithm of $\delta$ is actually more common in statistical learning theory \citep{mohri2018foundations}. We also keep our main focus on the high-probability bound, as shown in Section \ref{subsec:probability_bound}, to facilitate the comparisons with the existing findings. On the other hand, the proofs of these two bounds are totally different, where some technical gaps may exist to improve the convergence rate in the high-probability bound. While this is out of the scope of our paper, we leave the explorations to future research. Finally, compared with the state-of-the-art results in the literature (e.g., \citealt{padilla2022quantile}), our result does not require restrictive assumptions on the network parameters, such as boundedness and sparsity, which is theoretically stronger and can be practically more relevant. 

Meanwhile, we are able to establish a parallel lower bound for the expected excess risk using any sample-based policy.

\begin{theorem}[{\sc Lower Bound for Expected Excess Risk}]\label{thm:lower_expected_bound}
    Under Assumptions \ref{assumption_1} and \ref{assumption_2} and additional assumptions that (i) the density function of the conditional distribution $D|\bX$ is Lipschitz continuous and uniformly lower bounded by $\kappa>0$; (ii) the marginal density function of $\bX$ has a finite upper bound over $\mathcal{X}$, we have
    \begin{align}\label{expected_risk_lower_bound}
	\inf_{\hat{f}_n}\sup_{f_\rho\in \mathcal{H}^{\beta}([0,1]^p,B_0)} \mathbb{E} \Big[ \mathcal{R}^\rho(\hat{f}_n)-\mathcal{R}^\rho(f_\rho)\Big] &\ge c n^{-\frac{2\beta}{p+2\beta}},
\end{align}
where $c$ is some positive constant and the infimum is taken over all possible estimators based on the random samples $\bS_n=\{(\bx_i,d_i)\}_{i=1}^n$ from the joint distribution $(\bX, D)$ with the $\rho$-th quantile of the conditional distribution $D|\bX$ being $f_\rho(\bX)$.
\end{theorem}

The lower bound in Theorem \ref{thm:lower_expected_bound} matches the upper bound in Theorem \ref{thm:upper_expected_bound}, up to certain logarithm factors, indicating that the rate $n^{-\frac{2\beta}{p+2\beta}}$ is indeed minimax optimal. It further demonstrates the effectiveness of the DNN method in dealing with the data-driven, feature-based newsvendor problem. To our knowledge, our optimal rate result is also new to the literature on DNN theories and is not restricted to special distributions and network structures (unlike, e.g.,  \citealt{schmidt2020nonparametric}).

\section{Extensions and Further Discussions}\label{sec:extensions}

We further consider two practical extensions within the general framework established earlier. The first extension focuses on accommodating a learning environment with dependent data. We then explore how to mitigate the curse of dimensionality when the target function exhibits a composite structure.
Each subsection presents its background and conditions, and they stand alone without necessarily being connected.

\subsection{Learning from Dependent Data}\label{subsec:dependent_data}

In practice, the sample data is usually a time series where the feature information $\bx_i$ and demand $d_i$ are collected over a certain timeline. In such cases, imposing inter-independence assumptions on the data may be challenging or even unreasonable, as time series data commonly exhibits various dependence structures. Such a dependency can arise from auto-regressive effects, memory effects, or momentum patterns, where each observation at a given time point may be influenced by past observations. Consequently, it becomes necessary to extend the previously established theory when dealing with potentially dependent data.

In the following, we consider a general time series prediction problem where $\bS_T=\{\bZ_i:=(\bX_i, D_i)\}_{i=1}^T$ comes from some stochastic process. We note that no specific assumptions are imposed on the data-generating stochastic process so it could be non-stationary and non-mixing. One can imagine the sample data to be the chronological features and demand over the past several weeks, months, or quarters. The objective is to predict the inventory decision in the next period $(T+1)$ in order to minimize the expected newsvendor loss.

Given the potential dependence of the data, we are particularly interested in the risk conditioning on the past realizations of the stochastic process $\bS_T$. 
With a little abuse of notation, we define the \emph{path-dependent} conditioned risk as
$$\mathcal{R}^\rho(f)=\mathbb{E}[b(D_{T+1}-f(\bX_{T+1}))^++h(f(\bX_{T+1})-D_{T+1})^+\mid \bZ_1,\ldots, \bZ_T]$$
 for any function $f:\mathbb{R}^p\to\mathbb{R}$, 
which differs from the averaged version $\mathbb{E}[b(D_{T+1}-f(\bX_{T+1}))^++h(f(\bX_{T+1})-D_{T+1})^+]$ in that the latter averages the newsvendor loss over all possible historical paths. Meanwhile, the path-dependent risk is considered more reasonable due to the fact that we have already observed a feature-demand time series. Similarly, we define
\begin{equation*}
    f_{\rho}(\bx):=\argmin_{f} \mathcal{R}^\rho(f)=\inf\{q\geq0: F(q\vert \bX_{T+1}=\bx,\bZ_1,\ldots,\bZ_T)\geq \rho\},
\end{equation*}
which is the $\rho$-th quantile of the conditional CDF of $D_{T+1}$ given $\bX_{T+1}=\bx$ and past observations $(\bZ_1,\ldots,\bZ_{T})$.
Correspondingly, the empirical risk is given by 
\begin{equation}\label{empirical_risk_dependent}
    \mathcal{R}^\rho_T(f)=\frac{1}{T}\sum_{i=1}^T[\{b(D_{i}-f(\bx_{i}))^++h(f(\bx_{i})-D_{i})^+\}],
\end{equation}
and the minimizer within a given DNN function class is denoted by $\hat{f}_{D\!N\!N}:= \underset{f\in\mathcal{F}_{D\!N\!N}}{\argmin}\,  \mathcal{R}^{\rho}_T(f)$.
 
Again, we intend to provide a non-asymptotic upper bound for the excess risk $\mathcal{R}^\rho(\hat{f}_{D\!N\!N})-\mathcal{R}^\rho(f_\rho)$ under this new dependent-data scenario.

\begin{theorem}[{\sc Excess Risk Bound for Dependent Data}]\label{thm:excess_risk_dependent}
    Suppose the same network structure as in Theorem \ref{thm:excess_risk}.
    Under Assumptions \ref{assumption_1} and \ref{assumption_2} and for $T$ large enough, with probability at least $1-\delta$ over the random draw of $\bS_T$, where the elements are possibly dependent, we have
    \begin{equation}\label{excess_risk_bound_dependent}
        \begin{aligned}
   \mathcal{R}^\rho(\hat{f}_{D\!N\!N})-\mathcal{R}^\rho(f_\rho) &\le 2\sqrt{2}(b\bar{D}+h\mathcal{B})\bigg(C(\lfloor\beta\rfloor+1)^4p^{\lfloor\beta\rfloor+1}(\log(T))^2T^{-\frac{\beta}{2\beta+p}}+\sqrt{\frac{\log(1/\delta)}{T}}\bigg) +2\Delta,
   \end{aligned}
    \end{equation}
    where $C>0$ is a universal constant and 
        $\Delta:=\sup_{L\in\mathcal{F}_L}\Big(\mathbb{E}[L(\bZ_{T+1})\mid \bZ_1,\ldots,\bZ_T]-\sum_{t=1}^T\frac{1}{T}\mathbb{E}[L(\bZ_{t})\mid \bZ_1,\ldots,\bZ_{t-1}]\Big),$ 
    in which $L(\bz):=b(d-f(\bx))^++h(f(\bx)-d)^+$ and $\mathcal{F}_L:=\{L(f(\bx),d):f\in\mathcal{F}_{D\!N\!N}\}$.
\end{theorem}

A few comments are in order. First, we notice that the main difference between the dependent-data excess risk bound in Theorem \ref{thm:excess_risk_dependent} and the counterpart in Section \ref{sec:theoretical_guarantees} lies in the term $\Delta$. By definition, it evaluates the discrepancy between the target distribution and the distribution of past observations. Specifically, it can also serve as a natural measure of non-stationarity for the stochastic process $\bS_t$, which is tailored to adapt to the loss function $L$ and the hypothesis set $\mathcal{F}_{D\!N\!N}$ \citep{kuznetsov2015learning}. As the network size increases, the discrepancy also increases. 
In particular, when the sample data is i.i.d., it is apparent that $\Delta=0$ and this excess risk bound reduces to the one established in Section \ref{sec:theoretical_guarantees}. In general cases, it links to the familiar measures such as 
the total variation distance $(\Vert\cdot\Vert_{TV})$ and relative entropy $(D_{KL}(\cdot\Vert\cdot))$ between the conditional distribution of $\bZ_{T+1}$ given $\bZ_1,\ldots,\bZ_T$ (denoted by $\mathbb{P}(\cdot\vert \bZ_1,\ldots,\bZ_T)$) and the mixture distribution of sample margins $\sum_{t=1}^T\frac{1}{T}\mathbb{P}(\cdot\vert \bZ_1,\ldots,\bZ_{t-1})$ with an application of the Pinsker's inequality,
{\small\begin{align*}
    \Delta &\leq C\bigg\Vert \mathbb{P}(\cdot\vert \bZ_1,\ldots,\bZ_T)-\sum_{t=1}^T\frac{1}{T}\mathbb{P}(\cdot\vert \bZ_1,\ldots,\bZ_{t-1})\bigg\Vert_{TV}\leq \sqrt{2D_{KL}\bigg(\mathbb{P}(\cdot\vert \bZ_1,\ldots,\bZ_T)\bigg\Vert\sum_{t=1}^T\frac{1}{T}\mathbb{P}(\cdot\vert \bZ_1,\ldots,\bZ_{t-1})\bigg)},
\end{align*}
}provided that $L(\bz)\leq C$ holds uniformly for a positive constant $C$. From the inequality above, we can see more clearly that as the target distribution and the distribution of samples get closer, measured by a smaller total variation distance or relative entropy, the discrepancy measure $\Delta$ becomes smaller, and the stochastic process is also more stationary. The excess risk bound in Theorem \ref{thm:excess_risk_dependent} grows with $\Delta$, implying that when the time series becomes more non-stationary, it is more challenging to make predictions upon past observations because they may be less informative about future outcomes.

Second, there is additional flexibility in adjusting the weights of the historical data in the prediction. Recall that we directly compute the empirical risk as the average cost for all past observations in \eqref{empirical_risk_dependent}. This is a natural treatment for i.i.d. data, where each sample contributes equally. However, in a dependent-data scenario, it may be the case that observations closer to the prediction date or those occurring on the same day of the week are more informative and carry greater predictive power. Therefore, it is no longer appropriate to treat every past observation equally. Instead, assigning different weights to them can be beneficial. For example, we can assign higher weights to closer data points or samples from the same day of the week. Importantly, even with a flexible weighting scheme, the excess risk bound in Theorem \ref{thm:excess_risk_dependent} still holds, with a corresponding modification on the definition of the discrepancy measure $\Delta$ \citep{kuznetsov2015learning}.

\subsection{Mitigating Curse of Dimensionality}\label{subsec:dimensionality_issues}

Recall that the excess risk bound in Theorem \ref{thm:excess_risk} scales as $O(n^{-\frac{\beta}{2\beta+p}})$ up to a logarithmic factor, where $n$ is the sample size, $\beta$ represents the smoothness of the target function, and $p$ is the dimension of the feature variables. While the smoothness parameter is endogenously given, the dimensionality strongly influences the convergence rate so that a large dimension parameter $p$ could result in an extremely slow convergence rate, commonly referred to as the \emph{curse of dimensionality} in data science. A slow convergence rate thus requires many more samples so as to attain certain theoretical accuracy, which is usually impractical in real-world scenarios.
We here discuss how to mitigate this issue by maintaining a desirable convergence rate even with a very large $p$.

Generally speaking, there are two main approaches in the literature that impose assumptions either on the target function or the covariate distribution to achieve effective dimension reduction \citep{jiao2023deep}. For the former, it can be shown that if the target function exhibits a suitable composite structure, the effective dimension parameter appearing in the convergence rate could be much smaller than the nominal dimension $p$. In what follows, we use an example with a generalized linear demand (GLM) model to illustrate this.

\begin{example}[{\sc Generalized Linear Demand}]\label{example:GLM}
    We assume that the target function $f_{\rho}$ is a generalized linear function in the feature covariates:
    $f_{\rho}(\bx)=g(\btheta^\top \bx)$ for  $\btheta,\bx\in\mathbb{R}^{p+1},$ 
    where the first element of $\bx$ is $1$ without loss of generality and $g:\mathbb{R}\rightarrow\mathbb{R}$ is a univariate function. {\color{black}This model includes the commonly used MNL model and linear demand model as special cases, and it has been widely adopted in the literature dealing with contextual information (e.g., \citealt{ban2019big,miao2022context}).} With this special structure, we are able to show that the effective dimension in the convergence order is indeed only 1, indicating a faster rate without being plagued by the curse of dimensionality.
 
    \begin{theorem}[{\sc Excess Risk Bound for GLM}]\label{thm:excess_risk_GLM}
        {\color{black}Let the network width and depth be respectively $\tilde{\mathcal{W}}=\max\{\mathcal{W}, 2p\}$ and $\tilde{\mathcal{D}}=\mathcal{D}+3$, where $\mathcal{W}$ and $\mathcal{D}$ are defined according to \eqref{choice_of_depth_width} with $MN=\lfloor n^{\frac{1}{4\beta+2}}\rfloor$.
    Under Assumptions \ref{assumption_1} and \ref{assumption_2}(i) and the additional assumption that the target quantile function follows $f_{\rho}(\bx)=g(\btheta^\top \bx)$, where $g\in\mathcal{H}^{\beta}([a,c],B_0)$ for a given $\beta>0$, a finite constant $B_0>0$, and $-\infty<a<c<\infty$,
    with probability at least $1-\delta$ over the random draw of $\bS_n$ (which consists of i.i.d. samples from the joint distribution of $(\bX,D)$), we have}
    \begin{equation*}\label{excess_risk_bound_GLM}
        \mathcal{R}^{\rho}(\hat{f}_{D\!N\!N})-\mathcal{R}^{\rho}(f_{\rho})\leq 2\sqrt{2}(b\bar{D}+h\mathcal{B})\bigg(C(\lfloor\beta\rfloor+1)^4(\log(n))^2n^{-\frac{\beta}{2\beta+1}}+\sqrt{\frac{\log(1/\delta)}{n}}\bigg)
    \end{equation*}
    for $n$ large enough, where $C>0$ is a universal constant.
    \end{theorem}
\end{example}

Theorem \ref{thm:excess_risk_GLM} demonstrates that the excess risk bound now scales as $O(n^{-\frac{\beta}{2\beta+1}})$ under the GLM. {\color{black}This bound is obviously tighter, offering a much faster convergence rate, particularly when the feature dimension is large. This improvement is achieved under the additional GLM assumption, along with corresponding modifications to the proposed network design outlined in Theorem \ref{thm:excess_risk}.}
The rationale behind this is quite straightforward: a composite function behaves similarly to the structure of a DNN in \eqref{network_functional_composition}, which means that we can approximate each layer in the composition by a neural network and then stack them together. Since a linear function can be perfectly represented by a one-layer network as shown in Lemma \ref{lemma:linear_functions_nn}, we only need to consider the approximation of the univariate function $g$ in the outer layer, and this leads to an excess risk bound with an effective dimension of one. Theoretically, this result can be further extended to any general composite function with finite composition layers using similar reasoning, and interested readers are referred to \cite{shen2021deep} for more details.

{\color{black}

\section{Simulation Study and Implementation Guidance}\label{sec:simulation_implementation}

In this section, we conduct simulation experiments to examine the implementation of the DNN method, while the next section will demonstrate its effectiveness using real-world data. Our goal is to provide practical insights associated with our developed theory. The key insights are threefold. First, our numerical experiments reveal that the error of the DNN solution does not decrease monotonically with the network depth or width for a given amount of data. This finding confirms the trade-off between stochastic and approximation errors, which is central to our theoretical error analysis, and also highlights the importance of balancing these two error sources. 
Second, in terms of network design, it is more efficient and effective to use a deeper yet wide enough network during the tuning process to identify the optimal architecture. In fact, a wide range of network configurations can achieve comparably good performance.
Third, we numerically check the convergence rate of the DNN solution's excess risk, which is statistically significant and can be explained by our theoretical predictions.
The remainder of this section will elaborate on these points.

In our numerical experiments, we consider the data-generating mechanism represented by $D=f(\bX)+\epsilon$, where $\epsilon$ follows a standard normal distribution. Therefore, the target $\rho$-th conditional quantile function is $f_{\rho}(\bx)=f(\bx)+\Phi^{-1}(\rho)$, with $\Phi(\cdot)$ denoting the CDF of the standard normal distribution. Specifically, the following scenarios are examined:
\begin{itemize}
    \item[(a)] Univariate H\"older continuous (UHC) case: $f(x)=2x^{1/2}$; 
    \item[(b)] Multivariate logistic (ML) model: $f(\bx)=\frac{2\exp(\btheta^\top\bx)}{1+\exp(\btheta^\top\bx)}$, where $\btheta:=(4,-2,2,-1)^\top$;
    \item[(c)] Multivariate additive (MA) model: $f(\bx)=\exp(x_1-0.5)+2(x_2+x_3-1)^2+|x_4-0.5|$.
\end{itemize}  

These examples are carefully selected to encompass various feature dimensions and smoothness levels of the target function. 
The feature covariates $\bX$ are assumed to be uniformly distributed within the unit interval $[0,1]$ for the univariate case, and within $[0,1]^p$ for multivariate cases. To facilitate illustration, we present the numerical performance at three critical quantile levels: $\rho=0.25,0.5,0.75$, where $b+h$ is assumed to be $1$.

For each case mentioned above, we generate training data $(\bx_i^{train},d_i^{train})_{i=1}^n$ of size $n$ to train the empirical risk minimizer at $\rho\in\{0.25,0.5,0.75\}$, denoted by $\hat{f}_{n,\rho}$. To estimate the excess risk of the obtained DNN solution, we also generate testing data $(\bx_t^{test}, d_t^{test})_{t=1}^T$ of size $T$ from the same distribution as the training data. The testing error is then calculated in terms of the difference in average newsvendor loss between $\hat{f}_{n,\rho}$ and $f_\rho$ represented by
{\small\begin{align*}
    \frac{1}{T}\sum_{i=1}^T\bigg[b\left(d_i^{test}-f_{\rho}(\bx_i^{test})\right)^++h\left(f_\rho(\bx_i^{test})-d_i^{test}\right)^+\bigg]-\frac{1}{T}\sum_{i=1}^T\bigg[b\left(d_i^{test}-\hat{f}_{n,\rho}(\bx_i^{test})\right)^++h\left(\hat{f}_{n,\rho}(\bx_i^{test})-d_i^{test}\right)^+\bigg]. 
\end{align*}
}To ensure that this testing error approximates $\hat{f}_{n,\rho}$'s excess risk, we set $T=10^5$. Additionally, we report the mean of this error over $100$ replications to ensure robustness.

It is important to note that we desire to provide implementation guidance and insights that correspond to the theoretical findings established in this paper. However, this process presents certain challenges. While we have analyzed the stochastic and approximation errors associated with the DNN method, practical implementations also introduce an additional optimization error arising from the numerical solution of \eqref{DNN_solution_definition}. Any practical solution based on finite samples will inevitably involve all three types of errors. 
Our theoretical results have tackled only the first two by assuming that $\hat{f}_{n,\rho}$ precisely achieves the minimum empirical newsvendor loss within the given DNN function class. 
In practice, the optimization error can result from a variety of factors, including the choice of optimizer, initialization, and stopping criteria, making a unified theoretical analysis for this error particularly challenging.
For this reason, we strive to minimize the impact of the optimization error by many trials to ensure that the numerical results align as closely as possible with our theoretical findings, though some level of error may be unavoidable. Despite these challenges, we have identified useful practical guidelines, which will be discussed in detail later.

A detailed discussion of the optimization error, including those factors mentioned above contributing to it, lies beyond the scope of the current paper, and a vast literature has discussed the relevant issues (see, e.g., \citealt{glorot2010understanding}, \citealt{he2015delving}, \citealt{adcock2021gap}). In particular, we refer to \cite{adcock2021gap} for a comprehensive computational framework and extensive numerical illustrations for implementing DNNs in practical situations, and we largely follow their recommendations in our numerical experiments. 
Indeed, we have checked in our setting via many trials that their suggestions effectively achieve a tolerant optimization error, so we do not reiterate the related numerical experiments here. For example, the \emph{Adam} optimizer with an adaptively decaying learning rate can contribute to smooth and stable numerical outcomes, and network parameters can be initialized using symmetric uniform or normal distributions with small variances \citep{adcock2021gap}.

\subsection{Trade-Off in Network Size: Bigger Isn't Always Better}\label{subsec:network_size}

As demonstrated in Section \ref{sec:theoretical_guarantees}, there is a delicate balance between stochastic and approximation errors with respect to the network architecture. When dealing with limited data, a key observation from Theorem \ref{thm:stochastic_error} and Proposition \ref{prop:approximation_error} is that stochastic error increases with network complexity, while approximation error decreases with it. In other words, 
our theory suggests that an appropriately chosen network architecture is essential. The architecture must provide sufficient representational capacity while keeping stochastic error at a moderate level. A network that is too deep or too wide may not automatically yield optimal outcomes, let alone the high computational cost incurred.

Figure \ref{fig:fixed_depth_4} supports our claim. Using a sample size of 256, we fix the network depth and observe how the excess risk behaves as the network width increases. In the UHC and MA cases, the excess risk initially decreases. Linking to our developed theory, this reflects a dominance of approximation error over stochastic error. However, beyond a certain threshold, the excess risk stops decreasing and switches to rising as the width grows. This also aligns with our theory because the network at that point has already achieved an adequate approximation of the target function, causing the stochastic error to dominate. 
Interestingly, in the ML case, the excess risk keeps increasing with the network width. We hypothesize that this may be due to the high smoothness of the logistic model, where the approximation error is already minimal (negligible) even with a small network, leaving the stochastic error to dominate throughout. 
This phenomenon is also robust across different sample sizes (see also Figures \ref{fig:excess_risk_surface}, \ref{fig:excess_risk_surface_samples_256}, and \ref{fig:excess_risk_surface_samples_512}).
Figure \ref{fig:fixed_width_256} again witnesses this trade-off, this time in terms of the network depth with a fixed width.

\begin{figure}[!ht]
	\centering
 \includegraphics[scale=0.336]{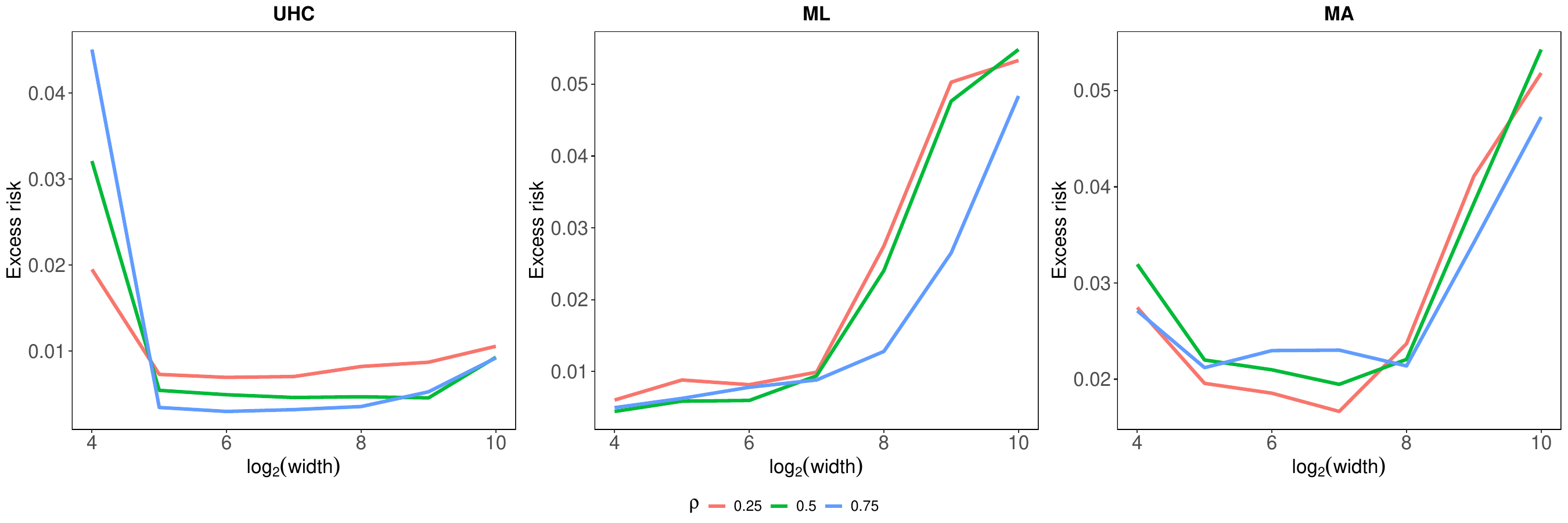}	
	\caption{Excess risk against network width with a fixed depth of 4 when $n=256$.}
		\label{fig:fixed_depth_4}
\end{figure}
\begin{figure}[ht]
	\centering
 \includegraphics[scale=0.336]{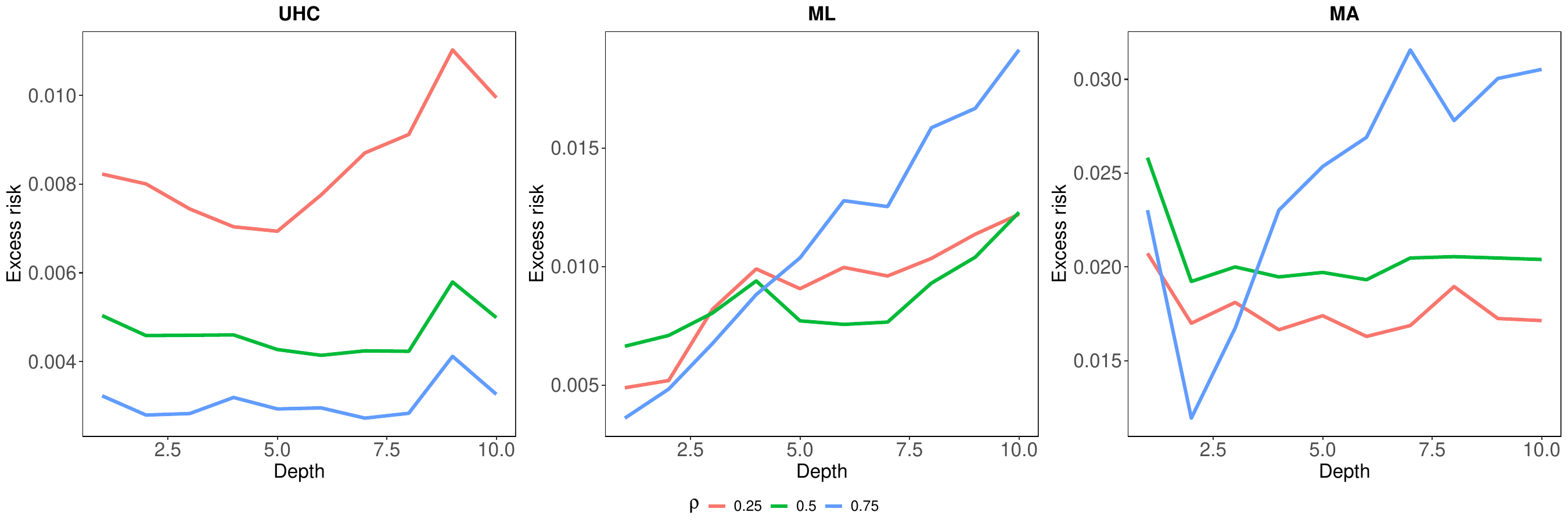}	
	\caption{Excess risk against network depth with a fixed width of 128 when $n=256$.}
		\label{fig:fixed_width_256}
\end{figure}

\subsection{Depth vs. Width: The Efficiency Edge of Deeper Networks}\label{subsec:numerical_depth_width}

Recall that in Theorems \ref{thm:excess_risk} and \ref{thm:upper_expected_bound}, we achieved an optimal balance between stochastic and approximation errors by specifying the network architecture that grows polynomially in the sample size. Intuitively, more data provides additional room for reducing the optimization error by enlarging network complexity, while keeping the stochastic error under control. 
We also have great flexibility in applying various combinations of network configurations. For example, one could set a fixed width and an expanding depth by setting $N=C$ for some $C\in\mathbb{N}_+$ and $M=\lfloor n^{\frac{p}{2p+4\beta0}}/C\rfloor$. Notably, the network size $\mathcal{S}=\sum_{i=0}^{\mathcal{D}}w_{i+1}(w_i+1)$ grows quadratically in the width but linearly in the depth. Such a special configuration can potentially save the number of parameters and also the computational cost. Moreover, as noted in \cite{adcock2021gap}, we observe that deep networks tend to exhibit better performance than shallower ones, provided they are sufficiently wide. We thus recommend using a network with a sufficiently large width and experimenting with different numbers of hidden layers to search for the best configuration.

We remind that Theorems \ref{thm:excess_risk} and \ref{thm:upper_expected_bound} hold uniformly with respect to target functions within a broad H\"older function class, representing somehow a ``worst-case" scenario. In practice, network configurations may deviate from those suggested in the theorems while still achieving good performance. For instance, as shown in Example \ref{example:linear_model}, linear functions can be effectively represented by a one-layer network with only two neurons in the hidden layer. For very smooth functions, such as linear ones, shallow networks with a moderate width may suffice. Nevertheless, our theoretical results deliver a critical message that optimal performance can be attained across a great range of network architectures, underscoring the flexibility of the DNN method. 

Figure \ref{fig:excess_risk_surface} validates this implication because plenty of different network configurations yield comparably good performance, evidenced by the relatively flat region corresponding to (nearly) minimum excess risk. This observation is persistent across different model setups, critical levels, and sample sizes (see also Figures \ref{fig:excess_risk_surface_samples_256}-\ref{fig:excess_risk_surface_samples_1024_0.75}). 
When comparing performances across models with varying dimensions, regularities, and sample sizes, the following general insights emerge: (1) for smoother target functions, a relatively small network can already do a decent job, and overly complex architectures may backfire, (2) as dimensionality increases, larger sample sizes are typically required to achieve a satisfactory precision, and (3) deeper or wider networks may be necessary when the target function exhibits greater nonlinearities and lower smoothness. These insights offer practical guidance, complementing the earlier recommendation. Indeed, by following the provided guidelines, practical experimentation can conveniently achieve good performance by tuning the network based on rough ideas of the target function's regularity and dimensionality.

\begin{figure}[!ht]
	\centering
 \includegraphics[scale=0.416]{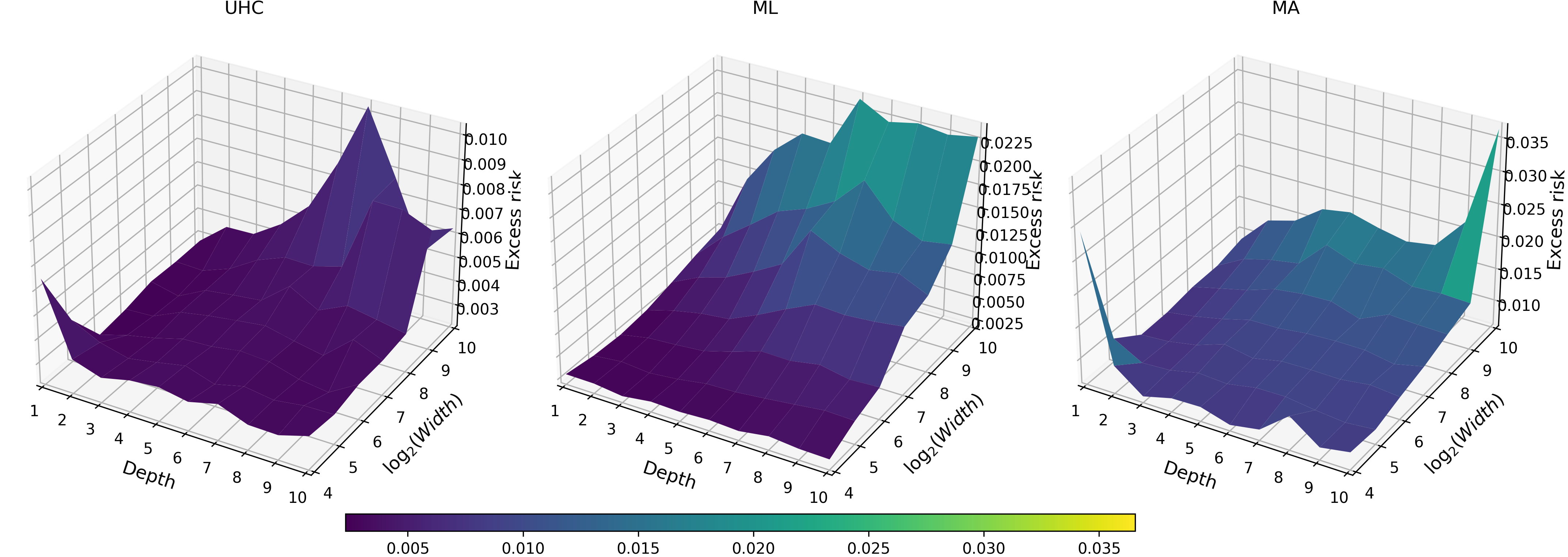}	
	\caption{Excess risk against network depth and width when $n=1024$ and $\rho=0.25$.}
		\label{fig:excess_risk_surface}
\end{figure}

\subsection{Understanding Convergence Rates in Sample Size}\label{subsec:understanding_convergence_rates}

One of our main theoretical results is the optimal convergence rate  $2\beta/(p+2\beta)$ of the expected excess risk in the sample size (up to logarithm factors), under properly designed network architectures. However, 
it is hard or even impossible to directly witness this rate in practical examples because of the presence of optimization errors, the somehow ``worst-case" nature of the upper bound, and the fact that the lower bound becomes numerically significant only asymptotically (as $c$ in \eqref{expected_risk_lower_bound} can be taken to be very small). Despite these challenges, we conduct a series of experiments with exponentially increasing sample sizes to investigate whether a clear convergent pattern shows up and how the rate links to our theory. To mitigate the impact of optimization errors, we explore a wide range of network configurations---varying combinations of the depth and width---for each sample size, and then select the minimum excess risk among them as the optimally balanced error for that sample size. This treatment is reasonable because the practically optimal design may deviate from the conservative plan outlined in Theorem \ref{thm:upper_expected_bound}, and it is stable because a variety of network designs can yield comparatively good performances, as found in Section \ref{subsec:numerical_depth_width}.

Figure \ref{fig:excess_risk_convergence_rate} presents the numerical convergence results. We can see clear polynomial convergence phenomena as overall, the linear regressions show a statistically significant fit to the log-transformed data points with very small $p$-values. As expected, the theoretically optimal rates are not directly visible. Instead, the observed rates are generally below $1$ (disregarding minor random perturbations in the ML case), which is the convergence rate of the stochastic error for a fixed network design, as established in Theorem \ref{thm:stochastic_error2}. Our findings are consistent with the theory in that the convergence rate of excess risk reflects a compromise between approximation and stochastic errors. The extent of the inferiority to the stochastic error rate depends on the degree of approximation error, which in turn is influenced by the regularity and dimensionality of the target function. Interestingly, the convergence rate approaches that of the stochastic error when the target function is very smooth, as the approximation error becomes nearly negligible (which is also consistent with our conjecture in Section \ref{subsec:network_size}); see the ML case in Figure \ref{fig:excess_risk_convergence_rate}. It deviates from $1$ to a smaller value when the smoothness is insufficient and the dimension effect kicks in, reflecting the compromise introduced by the approximation error (see the UHC and MA cases).

\begin{figure}[!ht]
	\centering
 \includegraphics[scale=0.286]{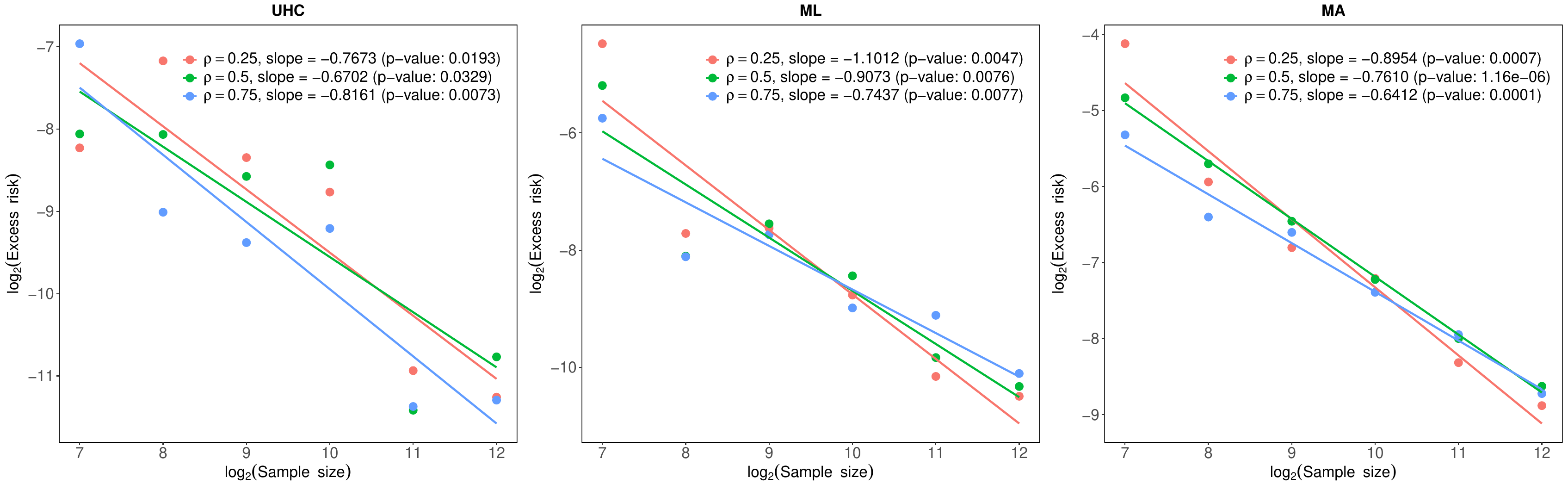}	
	\caption{Excess risk against sample size under proper network architectures.}
		\label{fig:excess_risk_convergence_rate}
\end{figure}

}

\section{Case Study: Inventory Management for a Food Supermarket Chain}\label{sec:case_study}

In this section, we showcase the practical application of DNNs for inventory decision-making. Specifically, we collaborate with a food company in China to collect unique sales data along with various features. Based on the empirical data, we have conducted comprehensive investigations, including testing the practical performance of DNNs and comparing it against many other data-driven benchmarks available in the literature.

\subsection{Data Description and Candidate Methods}\label{sec_7.1}

\begin{figure}[!ht]
	\centering
 \includegraphics[scale=0.366]{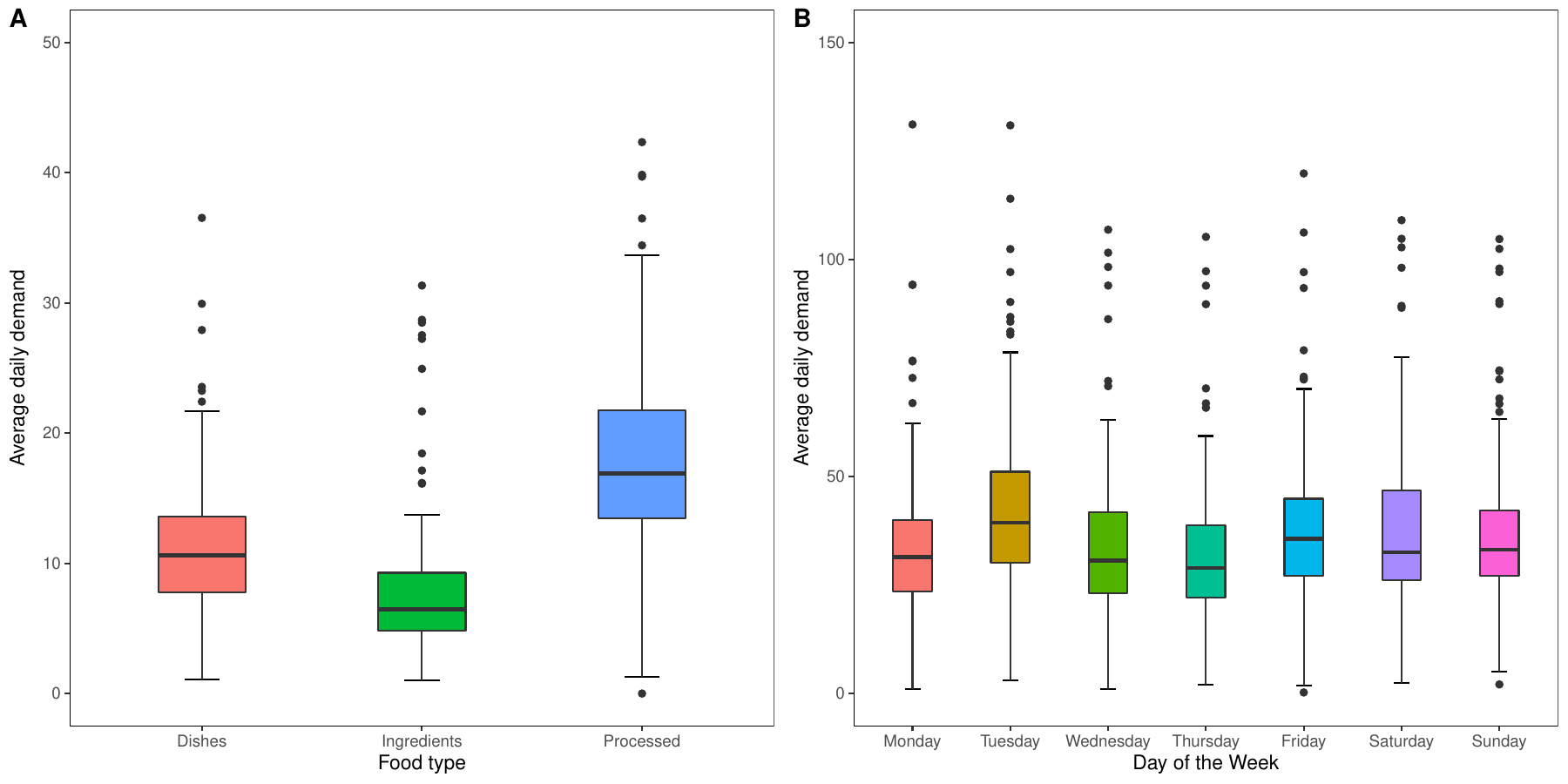}	
	\caption{Box plots of average daily demand for 131 stores versus food type and days of the week.}
		\label{fig:summary_plots}
\end{figure}

The dataset comprises daily sales records from our collaborated food company throughout the entire year of 2013. The firm operates 131 stores, and each store offers three types of food: processed food, dishes, and food ingredients. While there are plenty of choices within each food type, we focus on analyzing the demand for the general food category. Given the significant variation in food choices due to availability and seasonality, we aggregate the sales within each category to represent the actual demand for that category. Figure \ref{fig:summary_plots} provides the summary statistics for the average daily demand of stores, categorized by food types and days of the week. Each box plot depicts the demand variability among stores. By comparing the boxes in both panels, we can observe clear dependencies of demand on these two categorical variables. For instance, the overall demand for processed food is higher than the other two types, and the food demand is relatively higher on average on Tuesdays. {\color{black}To facilitate our numerical experiments, we convert the categorical features into binary vectors where the length equals the number of categories, with exactly one element set to 1 and the rest to 0.} 

We also consider an additional set of features in our experiments, which consist of 14 past demands and the associated \emph{operational statistics}. The so-called \emph{operational statistics} refer to the sample average and the differences between consecutive pairs of order statistics of these 14 past demands. Their relevance is theoretically justified in \cite{liyanage2005practical}, and their effectiveness has been observed in \cite{ban2019big}. In particular, it is important to note that when incorporating past demands as features, one should be cautious about potential time-series dependence in practice. Nevertheless, our analysis in the dependent-data scenario, as discussed in Section \ref{subsec:dependent_data}, allows us to incorporate this practical consideration within our framework seamlessly. 

We construct a set of feature variables to predict the demand. There are, in total, 54,780 demand records (for different food types and stores) along with their associated features. Following the common practice,
we {\color{black}chronologically} allocate 70\% of the data for training and validation purposes, while the remaining 30\% is reserved for testing. Specifically, the training and validation data are used to train and fine-tune the model and hyperparameters (if any) for each data-driven method under consideration. Once calibrated, these models are then applied to predict demand in the test set, and the corresponding losses are calculated to compare the performance of different methods. To calculate the newsvendor losses, we need to set values for the unit underage and overage costs. Following the standard treatment (e.g., \citealt{ban2020model}), we choose the average food price of CNY 18.01 as the underage cost $b$ without loss of generality. 
Since the overage cost $h$ is hard to determine in practice, depending on various factors such as the production cost and salvaging treatment, we compute and present the newsvendor losses for various $\rho:=b/(b+h)$, each corresponding to a unique $h$. These ratios implicitly relate to the product's profitability,
so we focus on the case where $\rho>0.5$ to reflect the practical situation.

The following set of data-driven methods are implemented on our dataset.
\begin{enumerate}
    \item DNN method. {\color{black}Given our dataset with relatively simple features, we adopt a three-layer neural network with 512 neurons per layer (and justify this choice in Section \ref{sec_supp:justifincation_network_design}). 
   The DNN model is built in Python using \emph{Pytorch}, and the \emph{Adam} optimization algorithm is applied with an initial learning rate of $0.001$. The default coefficients $\beta=(0.9,0.99)$ are used for computing running averages of gradients and their squares (see \citealt{kingma2014adam} for more details). We also adopt a standard early stopping rule to save computational time.}
    \item SAA method. We use the empirical quantile of the training data to give predictions for the demand in the test data.
    \item Linear empirical risk minimization (LERM) method. We implement the linear decision rules proposed in \cite{ban2019big}. Specifically, we investigate the performance of LERM under three different settings: one without any regularization and the remaining two with either $\ell_1$ or $\ell_2$ regularization.
    \item Kernel optimization (KO) method. We utilize the nonparametric Gaussian kernel-based smoothing method described in \cite{ban2019big}, which demonstrated superior performance over other alternatives on the emergency room data considered therein. 
    \item Reproducing Kernel Hilbert Space (RKHS) method. We conduct simulations using the RKHS approach proposed in \citet{bertsimas2022data}. For the simulation setting, we use a radial basis function kernel $K(x_1,x_2)=\exp(-\gamma\Vert x_1-x_2\Vert^2)$ with $\gamma=10^{-6}$ and a tuning parameter $\lambda=10^{-4}$ for regularization. More detailed descriptions can be found in \citet{bertsimas2022data}.
    
    \item Separated estimation and optimization (SEO) method. 
    A typical treatment in the literature is to assume a normal demand distribution, allowing for mean and variance estimation through linear regression with respect to feature variables. The inventory decision is then determined by the $\rho$-th quantile of the normal distribution.
\end{enumerate}

\subsection{Preliminary Experiments on a Small Dataset}

We first carry out an experiment on a subset of the full dataset to demonstrate the effectiveness of the DNN method in dealing with small datasets. Specifically, we focus on the sales record of one selected store that provides three types of food.
This subset consists of 990 daily sales records for all food types and includes a total of 38 feature variables {\color{black}($p=38$)}, such as the day of the week, past demands, and operational statistics, as previously introduced. Following the procedure outlined in the previous section, we test the performance of each considered method on this small dataset.
\begin{table}[!ht]
				\caption{\color{black}Newsvendor loss (in CNY \yen) for the selected store across different methods. Each cell displays the average newsvendor loss on the first line, followed by the 95\% interval on the second line.}
			\centering\color{black}
		{\def\arraystretch{1.66}  
			\begin{tabular*}{\textwidth}{@{\extracolsep{\fill}}lcccccccc}
				\hline
				\hline
				$\rho$   & 0.6 & 0.65 & 0.7 &0.75  & 0.8  & 0.85 & 0.9 & 0.95 \\ \hline
				\parbox[c]{1cm}{SAA}  	&\parbox[c]{1.46cm}{\centering\scriptsize 119.2\\(0, 479.1)}	&\parbox[c]{1.46cm}{\centering\scriptsize 108.0\\(0, 443.0)}	&\parbox[c]{1.46cm}{\centering\scriptsize 97.4\\(0, 425.0)}	&\parbox[c]{1.476cm}{\centering\scriptsize 86.2\\(6.0, 389.0)}	&\parbox[c]{1.476cm}{\centering\scriptsize 74.3\\(4.5, 353.0)}	&\parbox[c]{1.476cm}{\centering\scriptsize 60.7\\(3.2, 299.0)}	&\parbox[c]{1.476cm}{\centering\scriptsize 45.5\\(3.6, 194.5)} & \parbox[c]{1.46cm}{\centering\scriptsize 27.0\\(4.2, 47.2)}
                     \\ \hline
                    \parbox[c]{1cm}{SEO}  	&\parbox[c]{1.476cm}{\centering\scriptsize 90.8\\(3.0, 276.3)}	&\parbox[c]{1.476cm}{\centering\scriptsize 83.3\\(6.1, 243.4)}	&\parbox[c]{1.476cm}{\centering\scriptsize 75.2\\(4.6, 217.2)}	&\parbox[c]{1.476cm}{\centering\scriptsize 66.5\\(3.7, 189.6)}	&\parbox[c]{1.476cm}{\centering\scriptsize 51.1\\(3.7, 158.9)}	&\parbox[c]{1.476cm}{\centering\scriptsize 46.9\\(3.7, 123.2)} & \parbox[c]{1.46cm}{\centering\scriptsize 35.1\\(4.6, 78.2)}& \parbox[c]{1.476cm}{\centering\scriptsize 	21.5\\(3.1, 29.5)}
                \\ \hline 
				\parbox[c]{1cm}{LERM}  	&\parbox[c]{1.476cm}{\centering\scriptsize 90.5\\(2.1, 314.0)}	&\parbox[c]{1.476cm}{\centering\scriptsize 82.2\\(3.1, 298.3)}	&\parbox[c]{1.476cm}{\centering\scriptsize 71.9\\(4.0, 258.4)}	&\parbox[c]{1.476cm}{\centering\scriptsize 62.8\\(3.0, 238.9)}	&\parbox[c]{1.476cm}{\centering\scriptsize 53.1\\(1.7, 228.2)}	&\parbox[c]{1.476cm}{\centering\scriptsize 44.5\\(1.7, 175.7)} & \parbox[c]{1.476cm}{\centering\scriptsize 35.8\\(2.1, 148.7)}& \parbox[c]{1.476cm}{\centering\scriptsize 	23.5\\(2.0, 57.9)}
                 \\ \hline  
    			{\small LERM-$\ell_1$} 	&\parbox[c]{1.476cm}{\centering\scriptsize 86.4\\(3.9, 317.8)}	&\parbox[c]{1.476cm}{\centering\scriptsize 77.4\\(2.9, 284.1)}	&\parbox[c]{1.476cm}{\centering\scriptsize 68.0\\(2.4, 259.2)}	&\parbox[c]{1.476cm}{\centering\scriptsize 60.0\\(2.4, 223.5)}	&\parbox[c]{1.476cm}{\centering\scriptsize 50.9\\(2.3, 209.7)}	&\parbox[c]{1.476cm}{\centering\scriptsize 41.1\\(2.1, 159.9)} & \parbox[c]{1.476cm}{\centering\scriptsize 32.3\\(3.5, 123.5)}& \parbox[c]{1.476cm}{\centering\scriptsize 	21.4\\(3.8, 35.0)}\\ \hline  
    			{\small LERM-$\ell_2$} 	&\parbox[c]{1.476cm}{\centering\scriptsize 99.3\\(3.1, 293.8)}	&\parbox[c]{1.476cm}{\centering\scriptsize 80.5\\(2.8, 242.1)}	&\parbox[c]{1.476cm}{\centering\scriptsize 70.7\\(2.9, 214.8)}	&\parbox[c]{1.476cm}{\centering\scriptsize 77.4\\(2.6, 255.4)}	&\parbox[c]{1.476cm}{\centering\scriptsize 55.8\\(1.3, 197.7)}	&\parbox[c]{1.476cm}{\centering\scriptsize 50.1\\(1.5, 197.6)} & \parbox[c]{1.476cm}{\centering\scriptsize 49.5\\(4.4, 154.4)}& \parbox[c]{1.476cm}{\centering\scriptsize 	17.8\\(1.8, 51.5)}
 \\ \hline 
				\parbox[c]{1cm}{KO}  	&\parbox[c]{1.476cm}{\centering\scriptsize 90.5\\(2.1, 314.0)}	&\parbox[c]{1.476cm}{\centering\scriptsize 82.2\\(3.1, 298.3)}	&\parbox[c]{1.476cm}{\centering\scriptsize 71.8\\(3.8, 257.4)}	&\parbox[c]{1.476cm}{\centering\scriptsize 62.8\\(3.0, 238.9)}	&\parbox[c]{1.476cm}{\centering\scriptsize 53.1\\(1.8, 228.2)}	&\parbox[c]{1.476cm}{\centering\scriptsize 44.5\\(1.7, 175.6)} & \parbox[c]{1.476cm}{\centering\scriptsize 35.8\\(2.1, 148.7)}& \parbox[c]{1.476cm}{\centering\scriptsize 	23.5\\(2.0, 57.9)}
 \\ \hline  
                \parbox[c]{1cm}{RKHS}  	&\parbox[c]{1.476cm}{\centering\scriptsize 119.3\\(0.5, 478.8)}	&\parbox[c]{1.476cm}{\centering\scriptsize 108.1\\(0.4, 443.2)}	&\parbox[c]{1.476cm}{\centering\scriptsize 97.6\\(0.4, 424.6)}	&\parbox[c]{1.476cm}{\centering\scriptsize 86.4\\(6.0, 389.1)}	&\parbox[c]{1.476cm}{\centering\scriptsize 74.5\\(4.9, 351.5)}	&\parbox[c]{1.476cm}{\centering\scriptsize 60.9\\(3.5, 297.7)} & \parbox[c]{1.476cm}{\centering\scriptsize 45.5\\(4.2, 206.5)}& \parbox[c]{1.476cm}{\centering\scriptsize 	26.9\\(3.7, 64.4)}
 \\ \hline 
				\parbox[c]{1cm}{DNN}  	&\parbox[c]{1.476cm}{\centering\scriptsize 91.6\\(2.5, 283.3)}	&\parbox[c]{1.476cm}{\centering\scriptsize 77.6\\(4.6, 242.6)}	&\parbox[c]{1.476cm}{\centering\scriptsize 68.3\\(3.3, 205.2)}	&\parbox[c]{1.476cm}{\centering\scriptsize 59.5\\(4.5, 176.1)}	&\parbox[c]{1.476cm}{\centering\scriptsize 46.7\\(1.9, 136.2)}	&\parbox[c]{1.476cm}{\centering\scriptsize 38.6\\(1.7, 118.5)} & \parbox[c]{1.46cm}{\centering\scriptsize 24.8\\(0.7, 88.8)}& \parbox[c]{1.476cm}{\centering\scriptsize 	12.9\\(0.9, 41.3)} \\
				\hline
				\hline
			\end{tabular*}
			\label{table:sub-dataset:absolute}
		}
\end{table}

{\color{black}Table \ref{table:sub-dataset:absolute} presents the average newsvendor loss along with the 95\% empirical quantile interval, defined by the 0.025 and 0.975 quantiles of the losses, for each trained model on the common test dataset. We observe that the DNN method performs well across all critical levels and achieves the lowest average loss for critical levels above 0.7.  This suggests that it remains effective even with limited data. In particular, as the critical level $\rho$ increases, the overage penalty becomes less significant with a fixed underage cost, which reasonably explains the general improvement in prediction performance as the critical level grows.

From Figure \ref{fig:relative_loss_small_dataset}, we can further see the relative advantage of the DNN method over others, which tends to expand as the critical level increases. A similar pattern is also observed in the empirical experiments in \cite{oroojlooyjadid2020applying}.
This highlights the DNN method's ability to predict extreme distributional quantiles. In comparison, the DNN method shows an advantage in reducing predictive newsvendor loss, outperforming several recent competitive methods, and improving upon the empirical SAA rule that is commonly used in practice.

}

\subsection{Comparison of Results on the Entire Dataset}

We proceed to conduct numerical experiments on the entire dataset, where predictive models corresponding to different methods are built using feature information such as store ID and past demands {\color{black}($p=169$)}. The results are presented in Table \ref{table:full-dataset:absolute}. We clarify that these trained predictive models are essentially different from the ones in the preliminary experiments of the previous subsection. Therefore, the loss values in Table  \ref{table:full-dataset:absolute} are not necessarily smaller than their counterparts in Table \ref{table:sub-dataset:absolute} despite the larger sample size. 

\begin{table}[!ht]
				\caption{\color{black}Newsvendor loss (in CNY \yen) across different methods for the entire dataset. Each cell displays the average newsvendor loss on the first line, followed by the 95\% interval on the second line.}
			\centering\color{black}
		{\def\arraystretch{1.66}  
			\begin{tabular*}{\textwidth}{@{\extracolsep{\fill}}lcccccccc}
				\hline
				\hline
				$\rho$   & 0.6 & 0.65 & 0.7 &0.75  & 0.8  & 0.85 & 0.9 & 0.95 \\ \hline
				\parbox[c]{1cm}{SAA}  	&\parbox[c]{1.46cm}{\centering\scriptsize 160.6\\(0, 972.5)}	&\parbox[c]{1.46cm}{\centering\scriptsize 150.4\\(0, 936.5)}	&\parbox[c]{1.46cm}{\centering\scriptsize 141.3\\(0, 918.5)}	&\parbox[c]{1.476cm}{\centering\scriptsize 130.4\\(6.0, 882.5)}	&\parbox[c]{1.476cm}{\centering\scriptsize 117.9\\(4.5, 828.5)}	&\parbox[c]{1.476cm}{\centering\scriptsize 103.7\\(3.2, 756.4)}	&\parbox[c]{1.476cm}{\centering\scriptsize 86.4\\(6.0, 648.4)} & \parbox[c]{1.476cm}{\centering\scriptsize 61.8\\(7.6, 396.2)}
                     \\ \hline
                    \parbox[c]{1cm}{SEO}  	&\parbox[c]{1.476cm}{\centering\scriptsize 120.7\\(5.4, 539.7)}	&\parbox[c]{1.476cm}{\centering\scriptsize 115.1\\(6.5, 502.4)}	&\parbox[c]{1.476cm}{\centering\scriptsize 108.1\\(7.9, 455.6)}	&\parbox[c]{1.476cm}{\centering\scriptsize 99.3\\(9.1, 412.2)}	&\parbox[c]{1.476cm}{\centering\scriptsize 88.7\\(9.9, 366.8)}	&\parbox[c]{1.66cm}{\centering\scriptsize 76.0\\(10.0, 309.5)} & \parbox[c]{1.476cm}{\centering\scriptsize 60.7\\(9.4, 245.1)}& \parbox[c]{1.476cm}{\centering\scriptsize 	41.6\\(7.8, 150.2)}
                \\ \hline 
				\parbox[c]{0.86cm}{LERM}  	&\parbox[c]{1.476cm}{\centering\scriptsize 115.3\\(2.4, 540.7)}	&\parbox[c]{1.476cm}{\centering\scriptsize 117.0\\(6.0, 481.4)}	&\parbox[c]{1.66cm}{\centering\scriptsize 122.7\\(25.0, 430.4)}	&\parbox[c]{1.476cm}{\centering\scriptsize 77.1\\(2.9, 294.9)}	&\parbox[c]{1.476cm}{\centering\scriptsize 60.9\\(3.4, 243.4)}	&\parbox[c]{1.66cm}{\centering\scriptsize 80.9\\(12.2, 232.3)} & \parbox[c]{1.476cm}{\centering\scriptsize 31.4\\(1.3, 124.9)}& \parbox[c]{1.476cm}{\centering\scriptsize 	27.8\\(1.7, 79.9)}
                 \\ \hline  
    			{\small LERM-$\ell_1$} 	&\parbox[c]{1.66cm}{\centering\scriptsize 185.0\\(28.4, 654.5)}	&\parbox[c]{1.66cm}{\centering\scriptsize 131.3\\(13.4, 511.4)}	&\parbox[c]{1.476cm}{\centering\scriptsize 114.3\\(3.2, 409.9)}	&\parbox[c]{1.476cm}{\centering\scriptsize 108.2\\(2.7, 365.3)}	&\parbox[c]{1.476cm}{\centering\scriptsize 53.3\\(1.4, 247.6)}	&\parbox[c]{1.66cm}{\centering\scriptsize 102.8\\(47.8, 272.7)} & \parbox[c]{1.476cm}{\centering\scriptsize 30.3\\(1.0, 132.4)}& \parbox[c]{1.476cm}{\centering\scriptsize 	26.9\\(1.4, 98.7)}\\ \hline  
    			{\small LERM-$\ell_2$} 	&\parbox[c]{1.66cm}{\centering\scriptsize 166.9\\ (20.5, 619.0)}	&\parbox[c]{1.66cm}{\centering\scriptsize 142.5\\(18.5, 519.3)}	&\parbox[c]{1.66cm}{\centering\scriptsize 119.8\\(17.3, 429.2)}	&\parbox[c]{1.66cm}{\centering\scriptsize 91.1\\(14.5, 331.1)}	&\parbox[c]{1.476cm}{\centering\scriptsize 59.6\\(2.0, 247.8)}	&\parbox[c]{1.476cm}{\centering\scriptsize 49.5\\(2.2, 200.7)} & \parbox[c]{1.66cm}{\centering\scriptsize 56.9\\(23.3, 156.7)}& \parbox[c]{1.476cm}{\centering\scriptsize 	25.3\\(1.2, 111.3)}
 \\ \hline 
				\parbox[c]{1cm}{KO}  	&\parbox[c]{1.476cm}{\centering\scriptsize 106.3\\(2.5, 579.8)}	&\parbox[c]{1.476cm}{\centering\scriptsize 97.5\\(2.6, 535.3)}	&\parbox[c]{1.476cm}{\centering\scriptsize 88.5\\(2.5, 488.0)}	&\parbox[c]{1.476cm}{\centering\scriptsize 79.2\\(2.4, 449.5)}	&\parbox[c]{1.476cm}{\centering\scriptsize 69.6\\(2.3, 400.8)}	&\parbox[c]{1.476cm}{\centering\scriptsize 59.2\\(2.3, 339.6)} & \parbox[c]{1.476cm}{\centering\scriptsize 47.4\\(2.3, 252.1)}& \parbox[c]{1.476cm}{\centering\scriptsize 	33.3\\(2.5, 135.7)}
 \\ \hline  
                \parbox[c]{1cm}{RKHS}  &\parbox[c]{1.476cm}{\centering\scriptsize 824.2\\(578, 1707)}	&\parbox[c]{1.476cm}{\centering\scriptsize 281.2\\(45, 496)}	&\parbox[c]{1.66cm}{\centering\scriptsize 3249.4\\(2790, 3363)}	&\parbox[c]{1.476cm}{\centering\scriptsize 467.6\\(220, 1354)}	&\parbox[c]{1.476cm}{\centering\scriptsize 325.0\\(77, 1211)}	&\parbox[c]{1.476cm}{\centering\scriptsize 1117.0\\(933, 1162)} & \parbox[c]{1.476cm}{\centering\scriptsize 241.2\\(140, 261)}& \parbox[c]{1.476cm}{\centering\scriptsize 	62.5\\(7.5, 416.1)}
 \\ \hline 
				\parbox[c]{1cm}{DNN}  &\parbox[c]{1.476cm}{\centering\scriptsize 102.6\\(2.5, 537.7)}	&\parbox[c]{1.476cm}{\centering\scriptsize 87.1\\(2.1, 441.7)}	&\parbox[c]{1.476cm}{\centering\scriptsize 75.3\\(1.6, 371.9)}	&\parbox[c]{1.476cm}{\centering\scriptsize 61.2\\(1.7, 305.9)}	&\parbox[c]{1.476cm}{\centering\scriptsize 49.1\\(1.4, 242.1)}	&\parbox[c]{1.476cm}{\centering\scriptsize 36.9\\(1.0, 181.5)} & \parbox[c]{1.476cm}{\centering\scriptsize 24.9\\(0.7, 123.8)}& \parbox[c]{1.476cm}{\centering\scriptsize 	12.8\\(0.4, 66.2)} \\
				\hline
				\hline
			\end{tabular*}
			\label{table:full-dataset:absolute}
		}
\end{table}


On the entire dataset, the DNN method consistently outperforms other data-driven models for all considered critical levels. {\color{black}The KO method performs closest to the DNN model among the alternative methods when $\rho$ is relatively small, while the LERM algorithms tend to surpass it as $\rho$ is large (see also Figure \ref{fig:relative_loss_full_dataset}).
Additionally, the RKHS method turns out to have poor performance when both the number of features and sample size are large, indicating that its applicability to large datasets requires caution and further exploration, which is beyond the scope of this paper.

Finally, it is worth noting that the training time of the DNN model is quite short, typically around 15 seconds when employing a standard early stopping rule. Furthermore, the execution time for a trained DNN model is approximately 0.01 seconds. All these findings highlight the convenience and efficiency of applying the DNN model rapidly (see Table \ref{table:training-execution_times_full} for further comparison details). }

\section{Conclusion}\label{sec:conclusion}

In this paper, we propose to solve the feature-based newsvendor problem using a data-driven DNN method. 
We provide analytical non-asymptotic excess risk bounds in terms of the sample size and network configuration. Specifically, the excess risk bound can be separated into stochastic and approximation errors, offering a trade-off when selecting the network structure to achieve the optimal convergence rate. 
{\color{black}The DNN method can effectively mitigate model misspecification errors by achieving a vanishing excess risk when sufficient data is available. Moreover, as reflected in the trade-off between two types of errors, an excessively wide or deep network does not automatically lead to optimal performance. Instead, starting with a sufficiently wide network and adjusting the depth to find the optimal configuration is found to be quite effective. In general, the DNN method is flexible, allowing a range of network designs to perform well, though a specific good design depends on the available data and the properties of the target function.}
In addition, it can handle dependent data, providing reliable solutions even in the presence of weak non-stationarity. The convergence of the excess risk bound is also accelerated when the underlying target function exhibits a composite structure. Besides theoretical characterizations, we implement the DNN method on a unique real-world dataset. The numerical experiments showcase its good performance in terms of accurate prediction and efficient and fast training and execution.
Therefore, our work justifies the applicability of DNNs to the data-driven newsvendor problem and contributes to the advancement of theory and practice in OM problems by utilizing state-of-the-art DNN methods. {\color{black}Further discussions on some related topics can be found in Section \ref{sec_supp:further_discussion} in the online appendix.}

%
%




\setlength{\bibsep}{0pt}
{\small\bibliographystyle{informs2014} 
\bibliography{manuscript.bib} 
}



\ECSwitch


\ECHead{\centering Online Appendix to\\ ``Deep Neural Newsvendor"} 

 
\section{Auxiliary Results}\label{ec_sec:auxiliary}

To prove the non-asymptotic error bounds in Section \ref{sec:theoretical_guarantees}, we first present some preliminary concepts and theorems that will be used later regarding complexity theory (Chapter 3 of \citealt{mohri2018foundations}) and network learning (Chapters 11 and 12 of \citealt{anthony1999neural}).

\begin{definition}[{\sc Covering Number}]
    Let $\mathcal{F}$ be a class of functions from $\mathcal{X}^p$ to $\mathbb{R}$.
    For any given sequence $\bx=(x_1,\ldots,x_p)\in\mathcal{X}^p$, let
    \begin{equation*}
       \mathcal{F}|_{\bx}:= \{(f(x_1),\ldots,f(x_p)): f\in\mathcal{F}\}\subset \mathbb{R}^p.
    \end{equation*}
    For a positive constant $\delta$, the covering number $\mathcal{N}(\delta,\mathcal{F}|_{\bx},\|\cdot\|)$ is the smallest positive integer $N$ such that there exist  $f_1,\ldots,f_N\in\mathcal{F}|_{\bx}$ and  $\mathcal{F}|_{\bx}\subset \bigcup_{i=1}^N\{f: \|f-f_i\|\leq \delta\}$, i.e., $\mathcal{F}|_{\bx}$ is fully covered by $N$ spherical balls centered at $f_i$ with radius $\delta$ under the norm $\|\cdot\|$. Moreover, the uniform covering number $\mathcal{N}_n(\delta,\mathcal{F},\|\cdot\|)$ is defined to be the maximum covering number over all  $\bx$ in $\mathcal{X}^p$, i.e.,
    \begin{equation*}
        \mathcal{N}_n(\delta,\mathcal{F},\|\cdot\|):=\max\{\mathcal{N}(\delta,\mathcal{F}|_{\bx},\|\cdot\|): \bx\in\mathcal{X}^p\}.
    \end{equation*}
\end{definition}
\begin{definition}[Pseudo-Dimension]\label{def:pseudo_dimension}
    The pseudo-dimension of $\mathcal{F}$, denoted as Pdim$(\mathcal{F})$, is defined as the largest integer $m$ for which there exist $(\bx_1,y_1),\ldots,(\bx_m,y_m)\in\mathcal{X}^p\times\mathbb{R}$ such that for all $\boldsymbol{\eta}\in\{0,1\}^m$, there exists $f\in\mathcal{F}$ so that the following two arguments are equivalent:
    \begin{equation*}
        f(\bx_i)>y_i  \iff \eta_i, \quad \text{ for } i=1,\ldots,m.
    \end{equation*}
\end{definition}
Both the covering number and pseudo-dimension quantify the complexity of a given function class from different aspects. Indeed, the following results are available to establish the connections between them.
\begin{theorem}[{\sc Theorem 12.2 in \citealt{anthony1999neural}}]\label{thm:covering_number}
    If $\mathcal{F}$ is a set of real functions from a domain to a bounded interval $[0,B]$ and the pseudo-dimension of $\mathcal{F}$ is $\text{Pdim}(\mathcal{F})$, then
    \begin{equation*}
        \mathcal{N}_n(\delta,\mathcal{F},\|\cdot\|_{\infty})\leq\sum_{i=1}^{\text{Pdim}(\mathcal{F})}\binom{n}{i}\bigg(\frac{B}{\delta}\bigg)^i<\bigg(\frac{enB}{\delta \cdot \text{Pdim}(\mathcal{F})}\bigg)^{\text{Pdim}(\mathcal{F})},
    \end{equation*}
    for $n\geq \text{Pdim}(\mathcal{F})$.
\end{theorem}


\section{Supplementary Materials in Section \ref{sec:theoretical_guarantees}}

\subsection{Section \ref{subsec:probability_bound}: High-Probability Bound}

\proof{Proof of Lemma \ref{lemma:risk_decomposition}.}
Define the ``best" estimator $f^*$ in the function class $\mathcal{F}_{D\!N\!N}$ as 
\begin{equation*}
    f^*=\argmin_{f\in\mathcal{F}_{D\!N\!N}}\mathcal{R}^{\rho}(f).
\end{equation*}
We have
\begin{align*}
    &\mathcal{R}^{\rho}(\hat{f}_{D\!N\!N})-\mathcal{R}^{\rho}(f_{\rho})\\
    =&\left[\mathcal{R}^{\rho}(\hat{f}_{D\!N\!N})-\mathcal{R}_n^{\rho}(\hat{f}_{D\!N\!N}) \right]
    +\left[\mathcal{R}_n^{\rho}(\hat{f}_{D\!N\!N})-\mathcal{R}_n^{\rho}(f^*)\right]
    +\left[\mathcal{R}_n^{\rho}(f^*)-\mathcal{R}^{\rho}(f^*)\right]
    +\left[\mathcal{R}^{\rho}(f^*)-\mathcal{R}^{\rho}(f_{\rho})\right]\\
    \leq &\left[\mathcal{R}^{\rho}(\hat{f}_{D\!N\!N})-\mathcal{R}_n^{\rho}(\hat{f}_{D\!N\!N}) \right]
    +\left[\mathcal{R}_n^{\rho}(f^*)-\mathcal{R}^{\rho}(f^*)\right]
    +\left[\mathcal{R}^{\rho}(f^*)-\mathcal{R}^{\rho}(f_{\rho})\right]\\
    \leq& 2\sup_{f\in\mathcal{F_{D\!N\!N}}} \vert \mathcal{R}^{\rho}(f)-\mathcal{R}_n^{\rho}(f) \vert +\mathcal{R}^{\rho}(f^*)-\mathcal{R}^{\rho}(f_{\rho})\\
    =& 2 \sup_{f\in\mathcal{F_{D\!N\!N}}} \vert \mathcal{R}^{\rho}(f)-\mathcal{R}_n^{\rho}(f) \vert +\inf_{f\in\mathcal{F}_{D\!N\!N}}\left[\mathcal{R}^{\rho}(f)-\mathcal{R}^{\rho}(f_{\rho})\right],
\end{align*}
where the first inequality follows from the definition of $\hat{f}_{D\!N\!N}$ as the minimizer of $\mathcal{R}^\rho_n(f)$ in $\mathcal{F}_{D\!N\!N}$, the second inequality holds due to the fact that both $\hat{f}_{D\!N\!N}$ and $f^*$ belong to the function class $\mathcal{F}_{D\!N\!N}$, and the last equality is valid by the definition of $f^*$. 
\Halmos

\proof{Proof of Theorem \ref{thm:stochastic_error}.}
    By definition, we know that the function space $\mathcal{F}_{D\!N\!N}$ contains networks with the parameter $\phi$, depth $\mathcal{D}$, width $\mathcal{W}$, size $\mathcal{S}$, and number of neurons $\mathcal{U}$ so that $\|f_{\phi}\|_{\infty}\leq\mathcal{B}$. Since the  demand variable has a bounded domain as specified in Assumption \ref{assumption_1}, the newsvendor loss is also bounded:
    \begin{align*}
        L(f(\bx),d):=b(d-f(\bx))^++h(f(\bx)-d)^+\leq b\Bar{D}+h\mathcal{B},
    \end{align*}
    for any $f\in\mathcal{F}_{D\!N\!N}$. We can further obtain that
    {\color{black}\begin{equation}\label{EC.2}
        |L(f(\bX), D)-\mathbb{E} L(f(\bX), D)\vert \leq 2(b\Bar{D}+h\mathcal{B}).
    \end{equation}
    }Meanwhile, it is obvious that the newsvendor loss function $L(y, y')$ is $\max\{b,h\}$-Lipschitz in its both two arguments.

    For any given $\epsilon>0$, let $f_1,f_2,\ldots,f_{\mathcal{N}}$ be the anchor points of an $\epsilon$-covering for the function class $\mathcal{F}_{D\!N\!N}$, and we denote $\mathcal{N}:=\mathcal{N}_n(\epsilon,\mathcal{F}_{D\!N\!N},\Vert\cdot\Vert_\infty)$ as the covering number of $\mathcal{F}_{D\!N\!N}$ with radius $\epsilon$ under the norm $\Vert\cdot\Vert_\infty$. By definition, for any $f\in\mathcal{F}_{D\!N\!N}$, there exists an anchor $f_j$ for $j\in\{1,\ldots,\mathcal{N}\}$ such that 
    $\Vert f_j-f\Vert_\infty\le \epsilon$. The Lipschitz property of $L(y,y')$ then implies $\vert L(f(x),d)-L(f_j(x),d)\vert\le \max\{b,h\}\epsilon$, based on which an application of triangular inequality gives
    \begin{align}
        \vert\mathcal{R}^\rho(f)-\mathcal{R}^\rho_n(f)\vert&\leq \vert\mathcal{R}^\rho(f_j)-\mathcal{R}^\rho_n(f_j)\vert +\vert\mathcal{R}^\rho(f)-\mathcal{R}^\rho(f_j)\vert+\vert\mathcal{R}^\rho_n(f_j)-\mathcal{R}^\rho_n(f)\vert \notag\\
        &= \vert\mathcal{R}^\rho(f_j)-\mathcal{R}^\rho_n(f_j)\vert \notag\\&\quad+|\mathbb{E}_{\bZ}[L(f(\bX),D)-L(f_j(\bX),D)]|+\bigg|\frac{1}{n}\sum_{i=1}^n[L(f(\bX_i),D_i)-L(f_j(\bX_i),D_i)]\bigg|\notag\\
        &\le \vert\mathcal{R}^\rho(f_j)-\mathcal{R}^\rho_n(f_j)\vert+2\max\{b,h\}\epsilon \label{EC.1}.
    \end{align}
    Therefore, with a fixed $t>0$,
    \begin{align}
     \mathbb{P}\Big(\sup_{f\in\mathcal{F}_{D\!N\!N}}\vert\mathcal{R}^\rho(f)-\mathcal{R}^\rho_n(f) \vert&\ge t+2\max\{b,h\}\epsilon \Big)\overset{\mbox{\footnotesize(a)}}{\leq} \mathbb{P}\Big(\exists j\in\{1,\ldots,\mathcal{N}\} : \vert\mathcal{R}^\rho(f_j)-\mathcal{R}^\rho_n(f_j)\vert\ge t\Big)\notag\\
        &\overset{\mbox{\footnotesize(b)}}{\leq}\mathcal{N}_n(\epsilon,\mathcal{F}_{D\!N\!N},\Vert\cdot\Vert_\infty)\max_{j\in\{1,\ldots,\mathcal{N}\}}\mathbb{P}\Big( \vert\mathcal{R}^\rho(f_j)-\mathcal{R}^\rho_n(f_j)\vert\ge t\Big)\notag\\
&=\mathcal{N}_n(\epsilon,\mathcal{F}_{D\!N\!N},\Vert\cdot\Vert_\infty)\max_{j\in\{1,\ldots,\mathcal{N}\}} \mathbb{P}\Big( \vert \sum_{i=1}^n L(f_j(\bX_i),D_i)-\mathbb{E}[L(f_j(\bX),D)])\vert\ge nt\Big)\notag\\
        &\overset{\mbox{\footnotesize(c)}}{\leq}2\mathcal{N}_n(\epsilon,\mathcal{F}_{D\!N\!N},\Vert\cdot\Vert_\infty)\exp\left(-\frac{nt^2}{2(b\Bar{D}+h\mathcal{B})^2}\right)\label{EC.3},
    \end{align}
    where (a) holds due to \eqref{EC.1}; notice that we have denoted $\mathcal{N}:=\mathcal{N}_n(\epsilon,\mathcal{F}_{D\!N\!N},\Vert\cdot\Vert_\infty)$ in (b); and (c) comes from the Hoeffding's inequality and \eqref{EC.2}.
    
    For any $\delta>0$, let $\epsilon=1/n$ and $t=\sqrt{2}(b\bar{D}+h\mathcal{B})\sqrt{\log\Big(2\mathcal{N}_n(1/n,\mathcal{F}_{D\!N\!N},\Vert\cdot\Vert_\infty)/\delta\Big)/n}$ so that the right-hand side of \eqref{EC.3} equals to $\delta$, we have 
    \begin{align*}
     &\mathbb{P}\Big(\sup_{f\in\mathcal{F}_{D\!N\!N}}\vert\mathcal{R}^\rho(f)-\mathcal{R}^\rho_n(f) \vert\ge t+2(b+h)\epsilon \Big)\le \delta.
    \end{align*}
    In other words, with probability at least $1-\delta$,
        \begin{align*}
     &\sup_{f\in\mathcal{F}_{D\!N\!N}}\vert\mathcal{R}^\rho(f)-\mathcal{R}^\rho_n(f) \vert
     \leq t+2(b+h)\epsilon\\
     &\overset{\mbox{\footnotesize(a)}}{\leq} \frac{\sqrt{2}(b\bar{D}+h\mathcal{B})\sqrt{\log 2\mathcal{N}_n(1/n,\mathcal{F}_{D\!N\!N},\Vert\cdot\Vert_\infty)}}{\sqrt{n}}+\frac{\sqrt{2}(b\bar{D}+h\mathcal{B})\sqrt{\log (1/\delta)}}{\sqrt{n}}+\frac{2(b+h)}{n},
    \end{align*}
    where we have used the inequality $\sqrt{c+d}\leq \sqrt{c}+\sqrt{d}$ for any $c,d\geq0$ in (a). Furthermore, we know from Theorem \ref{thm:covering_number} that 
    \begin{align*}
    \sup_{f\in\mathcal{F}_{D\!N\!N}}\vert\mathcal{R}^\rho(f)-\mathcal{R}^\rho_n(f) \vert
     &\leq \frac{\sqrt{2}(b\bar{D}+h\mathcal{B})}{\sqrt{n}}\bigg(C\sqrt{\log \bigg(\frac{en^2\Bar{D}}{\text{Pdim}(\mathcal{F}_{D\!N\!N})}\bigg)^{\text{Pdim}(\mathcal{F}_{D\!N\!N})}}+\sqrt{\log (1/\delta)}\bigg)+\frac{2(b+h)}{n}\\
     &\leq \frac{\sqrt{2}(b\bar{D}+h\mathcal{B})}{\sqrt{n}}\bigg(C\sqrt{\text{Pdim}(\mathcal{F}_{D\!N\!N})\log (en^2\Bar{D})}+\sqrt{\log (1/\delta)}\bigg)+\frac{2(b+h)}{n},
    \end{align*}
    for $n\geq \text{Pdim}(\mathcal{F}_{D\!N\!N})$, where $\text{Pdim}(\mathcal{F}_{D\!N\!N})$ is the pseudo-dimension of $\mathcal{F}_{D\!N\!N}$ stated in Definition \ref{def:pseudo_dimension}; and $C$ is a generic constant that may have different values from time to time.

For the piecewise linear neural networks such as the ReLU network we consider in the current paper, Theorems 3 and 6 of \cite{bartlett2019nearly} show that there exist two universal constants $c$ and $C$ so that
\begin{equation}\label{bartlett_pdim}
    c\cdot \mathcal{S}\mathcal{D} \log(\mathcal{S}/\mathcal{D})\leq \text{Pdim}(\mathcal{F}_{D\!N\!N})\leq C\cdot \mathcal{S}\mathcal{D}\log(\mathcal{S}),
\end{equation}
where $\mathcal{S}$ and $\mathcal{D}$ are the size and depth of the given $\mathcal{F}_{D\!N\!N}$, respectively; see Section \ref{subsec:DNN_model} for more description. Then, with a probability of at least $1-\delta$, we have
\begin{align*}
    \sup_{f\in\mathcal{F}_{D\!N\!N}}\vert\mathcal{R}^\rho(f)-\mathcal{R}^\rho_n(f) \vert
     &\leq \frac{\sqrt{2}(b\bar{D}+h\mathcal{B})}{\sqrt{n}}\bigg(C\sqrt{\mathcal{S}\mathcal{D}\log(\mathcal{S})\log(n)}+\sqrt{\log (1/\delta)}\bigg)+\frac{2(b+h)}{n}\\
     &\overset{\mbox{\footnotesize(a)}}{\leq} \frac{\sqrt{2}(b\bar{D}+h\mathcal{B})}{\sqrt{n}}\bigg(C\sqrt{\mathcal{S}\mathcal{D}\log(\mathcal{S})\log(n)}+\sqrt{\log (1/\delta)}\bigg),
\end{align*}
for $n\geq C_1\cdot \mathcal{S}\mathcal{D}\log(\mathcal{S})$, where both $C$ and $C_1$ are universal constants; and (a) is valid since $\sqrt{\log(n)/n}$ decays slower than $1/n$ and $C$ represents a universal constant. 
\Halmos


{\color{black}\proof{Proof of Lemma \ref{lemma:linear_functions_nn}.}
We observe the simple identity that $x=\max\{x,0\}-\max\{-x,0\}=:\sigma(x)-\sigma(-x)$ for any $x\in\mathbb{R}$, where $\sigma(\cdot)$ denotes the ReLU activation function. This observation leads to the result in Lemma \ref{lemma:linear_functions_nn}.
\Halmos
}

\proof{Proof of Proposition \ref{prop:approximation_error}.}
    {\color{black}For the approximation of functions in $\mathcal{H}^\beta([0,1]^p,B_0)$ using ReLU neural networks with the width $\mathcal{W}$ and depth $\mathcal{D}$ specified as below
    \begin{equation*}
        \mathcal{W}=38(\lfloor\beta\rfloor+1)^2p^{\lfloor\beta\rfloor+1}N\lceil\log_2(8N)\rceil \text{ and }
        \mathcal{D}=21(\lfloor\beta\rfloor+1)^2M\lceil\log_2(8M)\rceil,
    \end{equation*}
    we know from Theorem 3.3 in \cite{jiao2023deep} that there exists a ReLU network function $f\in \mathcal{F}_{\mathcal{D},\mathcal{W},\mathcal{U},\mathcal{S},\mathcal{B}_{N,M}}$ so that
    \begin{equation*}
        |f(\bx)-f_{\rho}(\bx)|\leq 18B_0 (\lfloor\beta\rfloor+1)^2p^{\lfloor\beta\rfloor+(\beta\vee 1)/2}(NM)^{-2\beta/p},
    \end{equation*}
    for $\bx\in[0,1]^p$ except for a small set $\Omega$ with Lebesgue measure $\delta Kp$ and $\delta$ can be arbitrarily small, where we have added a subscript to $\mathcal{B}_{N,M}$ to indicate the dependence of the uniform upper bound of the network function on the free parameters in the network depth and width. This is reasonable because a larger network may result in a higher upper bound for the network function. To ensure that the output is uniformly bounded within $[-\mathcal{B},\mathcal{B}]$, we apply an additional truncation for $f(\bx)$. In fact, we observe that the function 
    $y=-\max\{\mathcal{B}-\max\{x,0\},0\}+\max\{\mathcal{B}-\max\{-x,0\},0\}$ maps any $x\in\mathbb{R}$ to a value $y\in[-\mathcal{B},\mathcal{B}]$, which corresponds a two-hidden-layer network with two neurons in each hidden layer. Therefore, we assert that there exists a ReLU network function $\tilde{f}\in \mathcal{F}_{D\!N\!N}$, stacked by $f$ and this two-hidden-layer network, such that the network width remains $\mathcal{W}$, the network depth becomes $\mathcal{D}+2$, and $\|\tilde{f}\|_\infty\leq\mathcal{B}$. Meanwhile, since Assumption \ref{assumption_2} states that $\|f_\rho\|_{\infty}\leq B_0\leq \mathcal{B}$, we can still have
    \begin{equation*}
        |\tilde{f}(\bx)-f_{\rho}(\bx)|\leq 18B_0 (\lfloor\beta\rfloor+1)^2p^{\lfloor\beta\rfloor+(\beta\vee 1)/2}(NM)^{-2\beta/p},
    \end{equation*}
    for $\bx\in[0,1]^p$ except a small set $\Omega$. Therefore, we get}
    \begin{align}\label{expectation_approximation_bound}
        \mathbb{E}|\Tilde{f}(\bx)-f_{\rho}(\bx)|&\leq 18B_0 (\lfloor\beta\rfloor+1)^2p^{\lfloor\beta\rfloor+(\beta\vee 1)/2}(NM)^{-2\beta/p}+\mathbb{P}(\Omega)\cdot \underset{\bx\in\Omega}{\sup}\{|\Tilde{f}(\bx)-f_{\rho}(\bx)|\}.
    \end{align}
    By Assumption \ref{assumption_2}(i), the marginal distribution of $\bX$ is absolutely continuous with respect to the Lebesgue measure, which means that $\underset{\delta\rightarrow0}{\lim\inf} \mathbb{P}(\Omega)=0$.

    {\color{black}Furthermore, since both $\|\Tilde{f}\|_{\infty}$ and $\|f_\rho\|_{\infty}$ are bounded, taking limit infimum with respect to $\delta$ on both sides of \eqref{expectation_approximation_bound} leads to}
    \begin{equation*}
        \mathbb{E}|\Tilde{f}(\bx)-f_{\rho}(\bx)|\leq 18B_0 (\lfloor\beta\rfloor+1)^2p^{\lfloor\beta\rfloor+(\beta\vee 1)/2}(NM)^{-2\beta/p}.
    \end{equation*}
    As the newsvendor loss function  is $\max\{b,h\}$-Lipschitz, it holds that
    \begin{equation*}
        \mathcal{R}^{\rho}(\Tilde{f})-\mathcal{R}^{\rho}(f_{\rho})\leq \max\{b,h\}\mathbb{E}|\Tilde{f}(\bx)-f_{\rho}(\bx)|\leq 18\max\{b,h\}B_0 (\lfloor\beta\rfloor+1)^2p^{\lfloor\beta\rfloor+(\beta\vee 1)/2}(NM)^{-2\beta/p}.
    \end{equation*}
    The conclusion in Proposition \ref{prop:approximation_error} is then valid since 
    \begin{align*}
        \underset{f\in\mathcal{F}_{D\!N\!N}}{\inf}\left[\mathcal{R}^{\rho}(f)-\mathcal{R}^{\rho}(f_{\rho})\right]\leq \mathcal{R}^{\rho}(\Tilde{f})-\mathcal{R}^{\rho}(f_{\rho}),
    \end{align*}
    and we finish the proof.
\Halmos

In fact, the approximation rate in Proposition \ref{prop:approximation_error} can be further accelerated once the excess risk has a desirable local quadratic structure.



\begin{proposition} \label{prop:approximation_error_new}
     Let $\mathcal{F}_{D\!N\!N}$ be the class of ReLU neural networks with the width and depth, respectively, specified as
    \begin{equation*}
        \mathcal{W}=38(\lfloor\beta\rfloor+1)^2p^{\lfloor\beta\rfloor+1}N\lceil\log_2(8N)\rceil \text{ and }
        \mathcal{D}=21(\lfloor\beta\rfloor+1)^2M\lceil\log_2(8M)\rceil,
    \end{equation*}
    where $\lceil a\rceil$ denotes the smallest integer that is no less than $a$.
    Under Assumption \ref{assumption_2} and \ref{assumption_3}, 
    for any $M,N\in\mathbb{N}_+$, we have
    \begin{equation*}
        \underset{f\in\mathcal{F}_{D\!N\!N}}{\inf}\left[\mathcal{R}^{\rho}(f)-\mathcal{R}^{\rho}(f_{\rho})\right]\leq 324 \Tilde{c}\max\{b,h\}^2B_0^2 (\lfloor\beta\rfloor+1)^4p^{2\lfloor\beta\rfloor+\beta\vee 1}(NM)^{-4\beta/p},
    \end{equation*}
    where $a\vee b=\max\{a,b\}$ and $\Tilde{c}=\max\{\max\{b,h\}/\delta_\rho^0, c_\rho^0\}$. 
\end{proposition}

\proof{Proof of Proposition \ref{prop:approximation_error_new}.}
The key observation is that 
\begin{equation}\label{approx_local_quadratic}
    \underset{f\in\mathcal{F}_{D\!N\!N}}{\inf}\left[\mathcal{R}^{\rho}(f)-\mathcal{R}^{\rho}(f_{\rho})\right]\leq \max\{\frac{\max\{b,h\}}{\delta_\rho^0}, c_\rho^0\}\inf_{f\in\mathcal{F}_{D\!N\!N}}\|f-f_{\rho}\|^2_{L^2(\nu)},
\end{equation}
and the remaining proof can be similarly conducted as in the proof of Proposition \ref{prop:approximation_error}.

To prove \eqref{approx_local_quadratic}, we first notice that
\begin{equation*}
    \mathcal{R}^\rho(f)-\mathcal{R}^{\rho}(f_{\rho})=\mathbb{E}[L(f(\bX),D)-L(f_\rho(\bX),D)]\leq \max\{b,h\}\mathbb{E}[|f(\bX)-f_\rho(\bX)|],
\end{equation*}
where $L(f(\bx),d):=b(d-f(\bx))^++h(f(\bx)-d)^+$ is $\max\{b,h\}$-Lipschitz. Therefore, if $\Vert f-f_\rho \Vert_{L^\infty(\mathcal{X}^0)}>\delta^0_\rho$, we have
\begin{equation*}
    \mathcal{R}^\rho(f)-\mathcal{R}^{\rho}(f_{\rho})\leq \max\{b,h\}\mathbb{E}[|f(\bX)-f_\rho(\bX)|]\leq \max\{b,h\}\mathbb{E}\Big[ \frac{|f(\bX)-f_\rho(\bX)|^2}{\delta_\rho^0}\Big]
\end{equation*}
On the other hand, if $\Vert f-f_\rho \Vert_{L^\infty(\mathcal{X}^0)}\leq\delta^0_\rho$, Assumption \ref{assumption_3} tells that $\mathcal{R}^\rho(f) -\mathcal{R}^\rho(f_\rho)\leq c^0_\rho \Vert f-f_\rho \Vert^2_{L^2(\nu)}.$ Combining these two cases, we have \eqref{approx_local_quadratic}.
\Halmos

\proof{Proof of Theorem \ref{thm:excess_risk}.}
    Combining Lemma \ref{lemma:risk_decomposition}, Theorem \ref{thm:stochastic_error}, and Proposition \ref{prop:approximation_error_new}, the following inequality holds for at least probability $1-\delta$,
    \begin{align}\label{first_error_term}
         \mathcal{R}^{\rho}(\hat{f}_{D\!N\!N})-\mathcal{R}^{\rho}(f_{\rho})&\leq 2\sqrt{2}(b\bar{D}+h\mathcal{B})\bigg(C_1\sqrt{\frac{\mathcal{S}\mathcal{D}\log(\mathcal{S})\log(n)}{n}}+\sqrt{\frac{\log(1/\delta)}{n}}\bigg)\\
         &\qquad+18\max\{b,h\}B_0 (\lfloor\beta\rfloor+1)^2p^{\lfloor\beta\rfloor+(\beta\vee 1)/2}(NM)^{-2\beta/p},\label{second_error_term}
    \end{align}
    where $n\geq C\cdot \mathcal{S}\mathcal{D}\log(\mathcal{S})$ for a large enough $C>0$ and $M,N\in\mathbb{N}_+$. To balance these two error terms, we choose proper tuning parameters $M$ and $N$ to optimize the convergence rate in the sample size. 
    {\color{black}We assume that $M=\lfloor n^{t_1}\rfloor$ and $N=\lfloor n^{t_2}\rfloor$ with $t_1,t_2\geq 0$ to be determined later. Then, we have $$\mathcal{W}=38(\lfloor\beta\rfloor+1)^2p^{\lfloor\beta\rfloor+1}N\lceil\log_2(8N)\rceil=O(n^{t_2}),$$
    \begin{equation*}
        \mathcal{D}=21(\lfloor\beta\rfloor+1)^2M\lceil\log_2(8M)\rceil=O(n^{t_1}),
    \end{equation*}
    and
    \begin{align*}
        \mathcal{S}&=w_1(w_0+1)+\sum_{i=1}^{\mathcal{D}-1}w_{i+1}(w_i+1)+w_{\mathcal{D}+1}(w_{\mathcal{D}}+1)\\
        &\leq \mathcal{W}(p+1)+(\mathcal{W}^2+\mathcal{W})(\mathcal{D}-1)+\mathcal{W}+1\\
        &=O(\mathcal{W}^2\mathcal{D})=O(n^{2t_2+t_1}),
    \end{align*}
    }where we have deliberately neglected those $\log(n)$ factors here since we are mainly interested in the convergence order in $n$, and \cite{schmidt2020nonparametric} also pointed out that they are more likely artifacts in the proofs related to the high-dimensional statistics. We then check the orders of the two error terms (neglecting those $\log(n)$ factors) in the excess risk bound:
    {\color{black}\begin{equation*}
        \eqref{first_error_term}=O(\sqrt{\frac{\mathcal{S}\mathcal{D}
        }{n}})=O(n^{t_1+t_2-1/2}) \quad\text{ and }\quad \eqref{second_error_term}=O((MN)^{-2\beta/p})=O(n^{-2(t_1+t_2)\beta/p}).
    \end{equation*}
    Let $t_1+t_2-1/2=-2t\beta/p$, we obtain $t_1+t_2=\frac{p}{2p+4\beta}$. Plugging it back into the relevant expressions above and letting $C$ be a generic constant whose value may change in different scenarios, we get
    \begin{align*}
        \mathcal{S}\mathcal{D}\leq C\mathcal{W}^2\mathcal{D}^2&\leq Ct_1^2t_2^2(\lfloor\beta\rfloor+1)^8 p^{2(\lfloor\beta\rfloor+1)}n^{\frac{p}{p+2\beta}}\log_2(n)\\
        &\leq C\left(\frac{p}{2p+4\beta}\right)^2(\lfloor\beta\rfloor+1)^8 p^{2(\lfloor\beta\rfloor+1)}n^{\frac{p}{p+2\beta}}\log_2(n)\\
        &\leq C(\lfloor\beta\rfloor+1)^8 p^{2(\lfloor\beta\rfloor+1)}n^{\frac{p}{p+2\beta}}\log_2(n).
    \end{align*}
    }Therefore, with probability $1-\delta$,
    \begin{align*}
        \mathcal{R}^{\rho}(\hat{f}_{D\!N\!N})-\mathcal{R}^{\rho}(f_{\rho})&\leq 2\sqrt{2}(b\bar{D}+h\mathcal{B})\bigg(C(\lfloor\beta\rfloor+1)^4p^{\lfloor\beta\rfloor+1}(\log(n))^2n^{-\frac{\beta}{2\beta+p}}+\sqrt{\frac{\log(1/\delta)}{n}}\bigg)\\
         &\qquad+18\max\{b,h\}B_0 (\lfloor\beta\rfloor+1)^2p^{\lfloor\beta\rfloor+1}n^{-\frac{\beta}{2\beta+p}}\\
         &\leq 2\sqrt{2}(b\bar{D}+h\mathcal{B})\bigg(C(\lfloor\beta\rfloor+1)^4p^{\lfloor\beta\rfloor+1}(\log(n))^2n^{-\frac{\beta}{2\beta+p}}+\sqrt{\frac{\log(1/\delta)}{n}}\bigg),
    \end{align*}
    for $n$ large enough. Finally, we notice that the stochastic error bound is valid only when $n\geq C\cdot \mathcal{S}\mathcal{D}\log(\mathcal{S})$. Under the current setting of $\mathcal{S}$ and $\mathcal{D}$, we know that $\mathcal{S}\mathcal{D}\log(\mathcal{S})$ grows slower than $n$. This condition then holds when $n$ is large enough.
\Halmos

Using similar arguments, the rate can be slightly improved to $O(n^{-\frac{2\beta}{4\beta+p}})$ with the addition of Assumption \ref{assumption_3}.

\begin{theorem}\label{thm:excess_risk_e_companion}
    {\color{black}Let the network width and depth be defined according to \eqref{choice_of_depth_width} with $NM=\lfloor n^{\frac{p}{2p+4\beta}}\rfloor$.
    Under Assumptions \ref{assumption_1}, \ref{assumption_2}, and \ref{assumption_3} and for all $n\geq C\cdot \mathcal{S}\mathcal{D}\log(\mathcal{S})$, where $C>0$ is sufficiently large,} 
    with probability at least $1-\delta$ over the random draw of $\bS_n$, we have
    \begin{equation*}
        \mathcal{R}^{\rho}(\hat{f}_{D\!N\!N})-\mathcal{R}^{\rho}(f_{\rho})\leq 2\sqrt{2}(b\bar{D}+h\mathcal{B})\bigg(C_1(\lfloor\beta\rfloor+1)^4p^{3\lfloor\beta\rfloor+1}(\log(n))^2n^{-\frac{2\beta}{4\beta+p}}+\sqrt{\frac{\log(1/\delta)}{n}}\bigg),
    \end{equation*}
    where $C_1>0$ is an independent universal constant.
\end{theorem}

\subsection{Section \ref{subsec:lower_bound}: Sharper Bound in Expectation}\label{subsec_append:proof_expected_bound}

To prove Theorem \ref{thm:upper_expected_bound}, we need to resort to a different risk decomposition to that in Lemma \ref{lemma:risk_decomposition}. Nevertheless, the general rationale and procedure are similar to the proofs in Theorem \ref{thm:stochastic_error}, Proposition \ref{prop:approximation_error}, and Theorem \ref{thm:excess_risk}.

\begin{lemma}\label{lemma:risk_decomposition2}
    The expected excess risk of the ERM $\hat{f}_{D\!N\!N}$ satisfies
    \begin{equation*}
        \mathbb{E}\left[\mathcal{R}^{\rho}(\hat{f}_{D\!N\!N})-\mathcal{R}^{\rho}(f_{\rho})\right]\leq \underbrace{\mathbb{E}\left[\mathcal{R}^{\rho}(\hat{f}_{D\!N\!N})-2\mathcal{R}_n^{\rho}(\hat{f}_{D\!N\!N})+\mathcal{R}^{\rho}(f_\rho)\right]}_\text{stochastic error}+\underbrace{2\underset{f\in\mathcal{F}_{D\!N\!N}}{\inf}\left\{\mathcal{R}^{\rho}(f)-\mathcal{R}^{\rho}(f_{\rho})\right\}}_\text{approximation error}.
    \end{equation*}
\end{lemma}

\proof{Proof of Lemma \ref{lemma:risk_decomposition2}.}
For any $f\in\mathcal{F}_{D\!N\!N}$,
\begin{align*}
	\mathbb{E}[\mathcal{R}^\rho(\hat{f}_{D\!N\!N})-\mathcal{R}^\rho(f_\rho)]
	&\overset{\mbox{\footnotesize(a)}}{\leq}\mathbb{E}[\mathcal{R}^\rho(\hat{f}_{D\!N\!N})-\mathcal{R}^\rho(f_\rho)]+2\mathbb{E}[\mathcal{R}^\rho_n(f)-\mathcal{R}^\rho_n(\hat{f}_{D\!N\!N})]\\
	&=\mathbb{E}[\mathcal{R}^\rho(\hat{f}_{D\!N\!N})-2\mathcal{R}^\rho_n(\hat{f}_{D\!N\!N})+\mathcal{R}^\rho(f_\rho)]+2\mathbb{E}[\mathcal{R}^\rho_{n}(f)-\mathcal{R}^\rho(f_\rho)]\\
 &=\mathbb{E}[\mathcal{R}^\rho(\hat{f}_{D\!N\!N})-2\mathcal{R}^\rho_n(\hat{f}_{D\!N\!N})+\mathcal{R}^\rho(f_\rho)]+2(\mathcal{R}^\rho(f)-\mathcal{R}^\rho(f_\rho)),
\end{align*}
where (a) holds due to the definition of the empirical risk minimizer $\hat{f}_{D\!N\!N}$ so that $\mathcal{R}^\rho_{n}(\hat{f}_{D\!N\!N})\le\mathcal{R}^\rho_{n}(f)$ for any $f\in\mathcal{F}_{D\!N\!N}$; and (b) is true because $(\bX_i,D_i)$ for $i=1,\ldots,n$ are i.i.d.. Since the inequality holds for any $f\in\mathcal{F}_{D\!N\!N}$, we have
\begin{align*}
	\mathbb{E}[\mathcal{R}^\rho(\hat{f}_{D\!N\!N})-\mathcal{R}^\rho(f_\rho)]&\le\mathbb{E}\Big[\mathcal{R}^\rho(\hat{f}_{D\!N\!N})-2\mathcal{R}^\rho_n(\hat{f}_{D\!N\!N})+\mathcal{R}^\rho(f_\rho)\Big]+2\inf_{f\in\mathcal{F}_{D\!N\!N}}\{\mathcal{R}^\rho(f)-\mathcal{R}^\rho(f_\rho)\}.
\end{align*}
This completes the proof. 
\Halmos

Although the risk decomposition in Lemma \ref{lemma:risk_decomposition2} differs from Lemma \ref{lemma:risk_decomposition}, we choose to use the same terminologies with a slight abuse of notation. Still, the stochastic error is mainly incurred by random realizations, and the approximation error is caused by the representation of the target function using a DNN.

\begin{theorem}\label{thm:stochastic_error2}
    Under Assumption \ref{assumption_1}, for $n\geq C\cdot \mathcal{S}\mathcal{D}\log(\mathcal{S})$ for a large enough $C>0$, the ERM $\hat{f}_{D\!N\!N}$ satisfies 
    \begin{equation*}
        \mathbb{E}\Big[\mathcal{R}^{\rho}(\hat{f}_{D\!N\!N})-2\mathcal{R}_n^{\rho}(\hat{f}_{D\!N\!N})+\mathcal{R}^{\rho}(f_\rho)\Big] \leq \frac{C\log(n)\mathcal{S}\mathcal{D}\log(\mathcal{S})}{n},
    \end{equation*}
    where $C>0$ is a universal constant.
\end{theorem}

\proof{Proof of Theorem \ref{thm:stochastic_error2}.}
We first notice that the expectation operator in the stochastic error $\mathbb{E}\left[\mathcal{R}^{\rho}(\hat{f}_{D\!N\!N})-2\mathcal{R}_n^{\rho}(\hat{f}_{D\!N\!N})+\mathcal{R}^{\rho}(f_\rho)\right]$ is taken with respect to the random samples $\bS=\{(\bX_i,D_i)\}_{i=1}^n$, based on which the ERM $\hat{f}_{D\!N\!N}$ is obtained within the function class $\mathcal{F}_{D\!N\!N}$. Let $\bS^\prime=\{(\bX^\prime_i,D^\prime_i)\}_{i=1}^n$ be another set of samples independent of $\bS$. For any $f\in\mathcal{F}_{D\!N\!N}$ and $\bX_i\in\mathcal{X}$, define
\begin{equation}\label{definition_functiong}
	g(f,\bX_i):=\mathbb{E}_{D|\bX}\Big[L(f(\bX_i),D_i)- L(f_\rho(\bX_i),D_i)\Big|\bX_i\Big],
\end{equation}
where $L(f(\bx),d):=b(d-f(\bx))^++h(f(\bx)-d)^+$. By definition, we can rewrite the stochastic error in terms of the newly defined function $g$ as below:
\begin{align}
	&\mathbb{E}\Big[\mathcal{R}^\rho(\hat{f}_{D\!N\!N})-2\mathcal{R}^\rho_n(\hat{f}_{D\!N\!N})+\mathcal{R}^\rho(f_\rho)\Big]\notag\\
    =&\mathbb{E}_{\bS} \Big[ \mathbb{E}_{\bS'}[L(\hat{f}_{D\!N\!N}(\bX'),D')]-\frac{2}{n}\sum_{i=1}^n [L(\hat{f}_{D\!N\!N}(\bX_i),D_i)]+L(f_{\rho}(\bX),D)\Big]\notag\\
    \overset{\mbox{\footnotesize(a)}}{=}& \mathbb{E}_{\bS} \Big[ \mathbb{E}_{\bS'}[L(\hat{f}_{D\!N\!N}(\bX'),D')]-\mathbb{E}_{\bS'}[L(f_{\rho}(\bX'),D')]-\frac{2}{n}\sum_{i=1}^n\Big( L(\hat{f}_{D\!N\!N}(\bX_i),D_i)-L(f_{\rho}(\bX_i),D_i)\Big)\Big]\notag\\
    =&\mathbb{E}_{\bS}\Big[\frac{1}{n}\sum_{i=1}^{n}\Big(\mathbb{E}_{\bS^\prime}[g(\hat{f}_{D\!N\!N},\bX'_i)]-2g(\hat{f}_{D\!N\!N},\bX_i)\Big)\Big], \label{sto_g}
\end{align}
where in (a), we have used the facts that $\mathcal{R}^\rho(f_\rho)$ is not sample-dependent (so that $\mathcal{R}^\rho(f_{\rho})=\mathbb{E}_{\bS}[\mathcal{R}^\rho(f_{\rho})]=\mathbb{E}_{\bS}\mathbb{E}_{\bS'}[\mathcal{R}^\rho(f_{\rho})]$) and $\bS_i:=(\bX_i,D_i)$ are i.i.d..

To bound \eqref{sto_g}, we notice the fact that for any random variable $W$, it holds that $\mathbb{E} [W]\le \mathbb{E} [\max\{W,0\}]=\int_{0}^\infty \mathbb{P}(W>t) dt$. Therefore, it suffices to derive the upper bound for its tail probability.
In fact, since $\hat{f}_{D\!N\!N}\in\mathcal{F}_{D\!N\!N}$, we have
\begin{align}\notag
	&\mathbb{P}\left(\frac{1}{n}\sum_{i=1}^{n}\Big(\mathbb{E}_{\bS^\prime}[g(\hat{f}_{D\!N\!N},\bX'_i)]-2g(\hat{f}_{D\!N\!N},\bX_i)\Big)>t \right)\\ \notag
	\le&\mathbb{P}\left(\exists f\in\mathcal{F}_{D\!N\!N}: \frac{1}{n}\sum_{i=1}^{n}\Big(\mathbb{E}_{\bS^\prime}[g( f,\bX^\prime_i)]-2g(f,\bX_i)\Big)>t \right)\\ 
	\overset{\mbox{\footnotesize(a)}}{=}&\mathbb{P}\left(\exists f\in\mathcal{F}_{D\!N\!N}:\mathbb{E}[g(f,\bX)]-\frac{2}{n}\sum_{i=1}^{n}\big[g(f,\bX_i)\big]>t\right),\label{sto_1}
\end{align}
where (a) holds because $\bS$ and $\bS'$ are i.i.d..
Since $L(y,y')$ is a $\max\{b,h\}$-Lipschitz function, we know from \eqref{definition_functiong} that
\begin{equation*}
    g(f, \bX_i)\le\max\{b,h\}\Vert f-f_\rho\Vert_\infty\le\max\{b,h\}(\mathcal{B}+\Bar{D}),
\end{equation*}
where the last inequality is valid since $f_\rho$ is a conditional quantile function, which has the same range as $D$, i.e., $0\leq f_\rho\leq\Bar{D}$.
On the other hand, given any $\bX_i$, the conditional quantile function $f_{\rho}$ should minimize $\mathbb{E}_{D|\bX}[L(f(\bX),D)]$ by the definitions in \eqref{feature_newsvendor} and \eqref{true_conditional_quantile}, implying that 
$g(f, X_i)\ge0$.
Consequently, we have $g(f,\cdot)\in[0,\max\{b,h\}(\mathcal{B}+\Bar{D})]$ for all $f\in\mathcal{F}_{D\!N\!N}$. Then, one can exactly follow the proof of Theorem 11.6 in \cite{gyorfi2002distribution_appendix} to derive the inequality below, whose tedious proof is omitted here but can be provided upon request:
\begin{equation}\label{inequality_sup_f}
    \mathbb{P}\bigg(\sup_{f\in\mathcal{F}_{D\!N\!N}}\frac{\mathbb{E}[g(f,\bX)]-\frac{1}{n}\sum_{i=1}^n[g(f,\bX_i)]}{s+\mathbb{E}[g(f,\bX)]+\frac{1}{n}\sum_{i=1}^n[g(f,\bX_i)]}>\epsilon\bigg)\leq 4\mathcal{N}_n\Big(\frac{s\epsilon}{16},\mathcal{G},\|\cdot\|_\infty\Big)\exp\Big(-\frac{\epsilon^2sn}{15\max\{b,h\}(\mathcal{B}+\Bar{D})}\Big),
\end{equation}
where the function class $\mathcal{G}$ is induced by $\mathcal{F}_{D\!N\!N}$ via $\mathcal{G}:=\{g(f,\bx): f\in\mathcal{F}_{D\!N\!N}\}$.
We then have
\begin{align}
    &\mathbb{P}\left(\exists f\in\mathcal{F}_{D\!N\!N}:\mathbb{E}[g(f,\bX)]-\frac{2}{n}\sum_{i=1}^{n}\big[g(f,\bX_i)\big]>t\right)\notag\\
    \overset{\mbox{\footnotesize(a)}}{\leq}& 
    \mathbb{P}\left(\underset{f\in\mathcal{F}_{D\!N\!N}}{\sup} \Big\{\frac{\mathbb{E}[g(f,\bX)]-\frac{1}{n}\sum_{i=1}^{n}[g(f,\bX_i)]}{2t+\mathbb{E}[g(f,\bX)]+\frac{1}{n}\sum_{i=1}^{n}[g(f,\bX_i)]}\Big\}>\frac{1}{3}\right)\notag\\
    \overset{\mbox{\footnotesize(b)}}{\leq}& 4\mathcal{N}_n\Big(\frac{t}{24},\mathcal{G},\|\cdot\|_\infty\Big)\exp\Big(-\frac{2tn}{135\max\{b,h\}(\mathcal{B}+\Bar{D})}\Big)\label{sto_2},
\end{align}
where (a) holds since when an even A implies an event B, $\mathbb{P}(A)\leq\mathbb{P}(B)$; and (b) is valid by choosing $\epsilon=1/3$ and $s=2t>0$ in \eqref{inequality_sup_f}.

Combining (\ref{sto_1}) and (\ref{sto_2}), for any truncating sequence $a_n>1/n$,  we have
\begin{align*}
	&\mathbb{E}_{\bS}\Big[\frac{1}{n}\sum_{i=1}^{n}\Big(\mathbb{E}_{\bS^\prime}[g(\hat{f}_{D\!N\!N},\bX'_i)]-2g(\hat{f}_{D\!N\!N},\bX_i)\Big)\Big]\\
	\le & \int_0^\infty \mathbb{P}\bigg( \frac{1}{n}\sum_{i=1}^{n}\Big(\mathbb{E}_{\bS^\prime}[g(\hat{f}_{D\!N\!N},\bX^\prime_i)]-2g(\hat{f}_{D\!N\!N},\bX_i)\Big)>t\bigg) dt\\
	\le & \int_0^{a_n} 1 dt+ \int_{a_n}^\infty4\mathcal{N}_{n}\Big(\frac{t}{24},\mathcal{G},\Vert\cdot\Vert_\infty\Big)\exp\Big(-\frac{2tn}{135\max\{b,h\}(\mathcal{B}+\Bar{D})}\Big) dt\\
	\overset{\mbox{\footnotesize(a)}}{\leq}& a_n + 4\mathcal{N}_{n}\Big(\frac{1}{24n},\mathcal{G},\Vert\cdot\Vert_\infty\Big) \int_{a_n}^\infty\exp\Big(-\frac{2tn}{135\max\{b,h\}(\mathcal{B}+\Bar{D})}\Big) dt\\
	=& a_n+ 4\mathcal{N}_{n}\Big(\frac{1}{24n},\mathcal{G},\Vert\cdot\Vert_\infty\Big)\exp\Big(-\frac{2a_nn}{135\max\{b,h\}(\mathcal{B}+\Bar{D})}\Big)\frac{135\max\{b,h\}(\mathcal{B}+\Bar{D})}{2n},
\end{align*}
where (a) holds because $\mathcal{N}_{n}\big(\delta,\mathcal{G},\Vert\cdot\Vert_\infty\big)$ is decreasing in $\delta$ and $a_n>\frac{1}{n}$. We can then choose 
 $a_n:=\log \Big(4\mathcal{N}_n(1/(24n),\mathcal{G},\Vert\cdot\Vert_\infty)\Big) \frac{135\max\{b,h\}(\mathcal{B}+\Bar{D})}{2n}$ to obtain
\begin{align*} \notag
	\mathbb{E}_S\Big(\frac{1}{n}\sum_{i=1}^{n}\bigg[\mathbb{E}_{S^\prime}\big\{g(\hat{f}_{D\!N\!N},X^\prime_i)\big\}-2g(\hat{f}_{D\!N\!N},X_i)\bigg]\Big)&\le a_n+\frac{135\max\{b,h\}(\mathcal{B}+\Bar{D})}{2n}\\
    &= \frac{135\max\{b,h\}(\mathcal{B}+\Bar{D})\log \Big(4e \mathcal{N}_n(\frac{1}{24n},\mathcal{G},\Vert\cdot\Vert_\infty)\Big)}{2n}.
\end{align*}

For any $f_1,f_2\in\mathcal{F}_{D\!N\!N}$, by the definition of $g$ in \eqref{definition_functiong} and its Lipschitz property, one can easily see that $\Vert g(f_1,\cdot)-g(f_2,\cdot)\Vert_\infty\le \max\{b,h\}\Vert f_1-f_2\Vert_\infty$ and therefore, $\mathcal{N}_n(1/(24n),\mathcal{G},\Vert\cdot\Vert_\infty)\le \mathcal{N}_n(1/(24\max\{b,h\}n),\mathcal{F}_{D\!N\!N},\Vert\cdot\Vert_\infty)$. By Theorem \ref{thm:covering_number} and \eqref{bartlett_pdim}, we then have
\begin{align*} \notag
	&\mathbb{E}_{\bS}\Big[\frac{1}{n}\sum_{i=1}^{n}\Big(\mathbb{E}_{\bS^\prime}[g(\hat{f}_{D\!N\!N},\bX'_i)]-2g(\hat{f}_{D\!N\!N},\bX_i)\Big)\Big]\\ \notag
	&\qquad\qquad\qquad\qquad\le \frac{135\max\{b,h\}(\mathcal{B}+\Bar{D})\log \Big(4e \mathcal{N}_n(\frac{1}{24\max\{b,h\}n},\mathcal{F}_{D\!N\!N},\Vert\cdot\Vert_\infty)\Big)}{2n}\\ \notag
    &\qquad\qquad\qquad\qquad \le \frac{135\max\{b,h\}(\mathcal{B}+\Bar{D})\Big[\log (4e) +\text{Pdim}(\mathcal{F}_{D\!N\!N})\log\big(24en^2\bar{D}\max\{b,h\}\big)\Big]}{n}\\\notag
    &\qquad\qquad\qquad\qquad \le \frac{135\max\{b,h\}(\mathcal{B}+\Bar{D})\Big[\log (4e) +C\mathcal{S}\mathcal{D}\log(\mathcal{S})\log\big(24en^2\bar{D}\max\{b,h\}\big)\Big]}{n}\\
    &\qquad\qquad\qquad\qquad \le \frac{C\log(n)\mathcal{S}\mathcal{D}\log(\mathcal{S})}{n},
\end{align*}
for $n\geq \text{Pdim}(\mathcal{F}_{D\!N\!N})$, where $\text{Pdim}(\mathcal{F}_{D\!N\!N})$ is the pseudo-dimension of $\mathcal{F}_{D\!N\!N}$, $\mathcal{S}$ and $\mathcal{D}$ are the size and depth of the given $\mathcal{F}_{D\!N\!N}$, and $C$ is a universal constant.

To conclude, we have established that for $n\ge\text{Pdim}(\mathcal{F}_{D\!N\!N})$,
\begin{align}\notag
	&\mathbb{E}\Big[\mathcal{R}^\rho(\hat{f}_{D\!N\!N})-2\mathcal{R}^\rho_n(\hat{f}_{D\!N\!N})+\mathcal{R}^\rho(f_\rho)\Big]\le \frac{C\log(n)\mathcal{S}\mathcal{D}\log(\mathcal{S})}{n},
\end{align}
where $C>0$ is a universal constant. This completes the proof. 
\Halmos

The upper bound for the approximation error is the same as that in Proposition \ref{prop:approximation_error_new} since they share the same definition. Therefore, we can prove Theorem \ref{thm:upper_expected_bound} as follows.

\proof{Proof of Theorem \ref{thm:upper_expected_bound}.}
    Under Assumption \ref{assumption_3}, we know from Proposition \ref{prop:approximation_error_new} that
    \begin{equation*}
        \underset{f\in\mathcal{F}_{D\!N\!N}}{\inf}\left[\mathcal{R}^{\rho}(f)-\mathcal{R}^{\rho}(f_{\rho})\right]\leq 324 \Tilde{c}\max\{b,h\}^2B_0^2 (\lfloor\beta\rfloor+1)^4p^{2\lfloor\beta\rfloor+\beta\vee 1}(NM)^{-4\beta/p},
    \end{equation*}
    where $\Tilde{c}=\max\{\max\{b,h\}/\delta_\rho^0, c_\rho^0\}$, $N$ and $M$ are two positive integers, and $\mathcal{W}=38(\lfloor\beta\rfloor+1)^2p^{\lfloor\beta\rfloor+1}N\lceil\log_2(8N)\rceil \text{ and }
        \mathcal{D}=21(\lfloor\beta\rfloor+1)^2M\lceil\log_2(8M)\rceil$. Therefore, Lemma \ref{lemma:risk_decomposition2} and Theorem \ref{thm:stochastic_error2} imply that
    \begin{equation*}
        \mathbb{E}\left[\mathcal{R}^{\rho}(\hat{f}_{D\!N\!N})-\mathcal{R}^{\rho}(f_{\rho})\right]\leq \frac{C\log(n)\mathcal{S}\mathcal{D}\log(\mathcal{S})}{n}+648 \Tilde{c}\max\{b,h\}^2B_0^2 (\lfloor\beta\rfloor+1)^4p^{2\lfloor\beta\rfloor+\beta\vee 1}(NM)^{-4\beta/p}.
    \end{equation*}
    {\color{black}We want to find the optimal convergence rate in the sample size $n$ by tuning parameters $N$ and $M$. 
    Following a similar argument as in the proof of Theorem \ref{thm:excess_risk}, we can choose $NM=\lfloor n^{\frac{p}{2p+4\beta}}\rfloor$ to reach the optimal rate. Under this setting, we also have
    \begin{equation*}
        \mathcal{S}\mathcal{D}\leq C n^{\frac{p}{p+2\beta}}(\log(n))^2.
    \end{equation*}
    }Therefore, 
    \begin{align*}
         \mathbb{E}\left[\mathcal{R}^{\rho}(\hat{f}_{D\!N\!N})-\mathcal{R}^{\rho}(f_{\rho})\right]&\leq \frac{C(\log(n))^4 n^{\frac{p}{p+2\beta}}}{n}+Cn^{-\frac{2\beta}{p+2\beta}}\\
         &\leq C(\log(n))^4n^{-\frac{2\beta}{p+2\beta}},
    \end{align*}
    where $C$ is a universal constant depending on the model parameters.
\Halmos

\proof{Proof of Theorem \ref{thm:lower_expected_bound}.}
 The proof is conducted in the following way. We first construct a finite subset of the target function space, then reduce the problem to multiple hypothesis testing, and finally establish the lower bound by applying Fano's inequality \citep{scarlett2019introductory}.
    We can rewrite the demand-feature model as a nonparametric quantile regression problem
     $D=f_\rho(\bX)+\epsilon$, where $D \in \mathbb{R}$ is the response (demand), $\bX \in\mathcal{X}$ is a $p$-dimensional vector of predictors (features), and $\epsilon$ is an unobservable random variable that may depend on $\bX$. Moreover, it satisfies that $\mathbb{E}\{\vert\epsilon\vert\mid \bX\}<\infty$ and its $\rho$-th conditional quantile given $\bX$ is zero. According to Assumption \ref{assumption_2}, the function $f_\rho:\mathcal{X}\to \mathbb{R}$  is H\"older continuous of order $\beta>0$, namely,  $f_\rho\in\mathcal{H}^{\beta}([0,1]^p,B_0)$. Indeed, this consists of our target function space.
	
    The first step can be conducted using the Varshamov-Gilbert lemma (see, e.g., Lemma D.2 of \citealt{lu2021machine_appendix}). Specifically, for any positive integer $m$ satisfying $m^{p}\ge8$, there exists a subset $\mathcal{V}=\{v^{(0)},\cdots,v^{(2^{m^{p}/8})}\}$ of $m^{p}$-dimensional hypercube $\{0,1\}^{m^{p}}$ such that $v^{(0)}=(0,\ldots,0)$ and the $\ell_1$ distance between every two elements is larger than $m^{p}/8$, i.e.,
	$$\sum_{i=1}^{m^{p}} \Vert v^{(j)}_i-v^{(k)}_i\Vert_{1}\ge \frac{m^{p}}{8},\quad{\rm for\ all\ }0\le j\not=k\le2^{m^{p}/8}.$$
	In a similar vein to the proof of Theorem D.1 in \citet{lu2021machine_appendix}, we consider a simple $C^\infty$ bump function supported on $[0,1]^d$ defined by
	$$g(\bx)=\prod_{i=1}^{p}h(x_i), \quad \text{ for }\bx=(x_1,\ldots,x_p),$$
	where $h:\mathbb{R}\to\mathbb{R}$ is a positive integrable function in $C^\infty(\mathbb{R})$. We then construct multiple hypothesises on the regular grid $(\bx^{(j)}),j\in[m]^{p}$ by
	$$f_k(\bx)=\sum_{j\in[m]^{p}}v^{(k)}_j\frac{\omega}{m^{\beta+p/2}} g(m(\bx-\bx^{(j)})),\quad k=1,2,\ldots,2^{m^{p}/8},$$
	where $\omega$ is a constant to be determined later. One can easily check that $f_k$ is $\beta$-H\"older smooth and
	\begin{align}\label{lower_b1}
		\frac{C_1\cdot\omega}{m^{2\beta}}\leq \Vert f_i-f_k \Vert^2_2=\frac{C\cdot\omega}{m^{2\beta+p}}\sum_{j\in[m]^{p}}\vert v_j^{(i)}-v_j^{(k)}\vert_1\leq \frac{C_2\cdot\omega}{m^{2\beta}} ,\quad \text{for } 1\le i\not= k\le 2^{m^{p}/8},
	\end{align}
 where $C,C_1$ and $C_2$ are universal constants.
	
    Next, we reduce the problem to multiple hypothesis testing. Let $V$ be a random index uniformly drawn from the set $\{1,\ldots,2^{m^{p}/8}\}$, and $\bS_n=\{(\bX_i,D_i)\}_{i=1}^n$ be the random samples drawn from the distribution $P_V$ associated with the model $D=f_V(\bX)+\epsilon$ for $V\in \{1,\ldots,2^{m^{p}/8}\}$. 
    We remark that although the data-generating process could be different depending on the realization of $V$, the random variable $\epsilon$ is fixed so that the distributional properties assumed in Theorem \ref{thm:lower_expected_bound} are shared among these models. 
    For any estimator $\hat{f}_n$ based on the sample $\bS_n=\{(\bX_i,D_i)\}_{i=1}^n$, we further denote $\hat{V}$ as the index corresponding to the closest $f_j$ to $\hat{f}_n$ in terms of the metric $\Vert\cdot\Vert_{2}^2$, i.e., $\hat{V}=\underset{v=1,\ldots,2^{m^{p}/8}}{\argmin} \Vert f_v-\hat{f}_n\Vert_{2}^2$. Using the triangle inequality and \eqref{lower_b1}, if $\Vert f_V-\hat{f}_n\Vert_{2}^2<\varepsilon/2$ with $\varepsilon:=C_1\cdot \omega/m^{2\beta}$, we can claim that $\hat{V}=V$ because for any $j\neq V$,
    \begin{equation*}
        \Vert f_j-\hat{f}_n\Vert_{2}^2\geq  \Vert f_j-f_V\Vert_{2}^2 - \Vert f_V-\hat{f}_n\Vert_{2}^2\geq \varepsilon-\varepsilon/2=\varepsilon/2.
    \end{equation*} 
    Therefore, we have $\mathbb{P}(\Vert f_V-\hat{f}_n\Vert_{2}^2< \varepsilon/2)\leq \mathbb{P}(\hat{V}= V)$, implying that
$$\mathbb{P}(\Vert f_V-\hat{f}_n\Vert_{2}^2\ge \varepsilon/2)\ge \mathbb{P}(\hat{V}\not= V).$$
We then know that
\begin{align*}
	\sup_{f_\rho\in \mathcal{H}^{\beta}([0,1]^p,B_0)} \mathbb{E} \Vert \hat{f}_n-f_\rho\Vert_{2}^2 &\overset{\mbox{\footnotesize(a)}}{\geq}\sup_{f_\rho\in \mathcal{H}^{\beta}([0,1]^p,B_0)} \varepsilon\times \mathbb{P}(\Vert \hat{f}_n-f_\rho\Vert_{2}^2\ge\varepsilon)\\
	&\overset{\mbox{\footnotesize(b)}}{\geq} \varepsilon\times \max_{v=1,\ldots,2^{m^{p}/8}} \mathbb{P}(\Vert \hat{f}_n-f_v\Vert_{2}^2\ge \varepsilon)\\
	&\ge \varepsilon\times \max_{v=1,\ldots,2^{m^{p}/8}} \mathbb{P}(\hat{V}\not= v)\\
	&\overset{\mbox{\footnotesize(c)}}{\geq} \varepsilon\times \frac{1}{2^{m^{p}/8}} \sum_{v=1}^{2^{m^{p}/8}}\mathbb{P}(\hat{V}\not= v)\\
	&\overset{\mbox{\footnotesize(d)}}{\geq} \varepsilon\times \left(1-\frac{I(V;\bS_n)+\log2}{\log(2^{m^{p}/8})}\right)\\
	& = \frac{C_1\cdot\omega}{m^{2\beta}} \times \left(1-\frac{I(V;\bS_n)+\log2}{{m^{p}\log(2)/8}}\right),
\end{align*}
where (a) follows from taking supremum after an application of the Markov inequality; (b) is valid because $\{f_{v}\}_{v=1}^{2^{m^{p}/8}}\subset \mathcal{H}^{\beta}([0,1]^p,B_0)$; (c) is due to the fact that the maximum value of a real sequence is larger than its average; and (d) follows from the Fano's inequality (Theorem 1 in \citealt{scarlett2019introductory}), where $I(V; \bS_n)$ denotes the mutual information of $V$ and $\bS_n$ representing how much information $V$ reveals about $\bS_n$ \citep{scarlett2019introductory}. 

To obtain the lower bound, we then need to bound the mutual information $I(V; \bS_n)$. Recall that $\bS_n=\{(\bX_i,D_i)\}_{i=1}^n$ contains $n$ i.i.d. random samples, then we can tensorize the mutual information $I(V;\bS_n)$ through a sum of $n$ separate mutual information terms. Indeed, Lemma 2 in \citet{scarlett2019introductory} tells that
$$I(V;\bS_n)\le \sum_{i=1}^n I(V;(\bX_i,D_i))=n\times I(V;(\bX,D)).$$
By the KL divergence-based bounds (see, e.g., Lemma 4 in \citealt{scarlett2019introductory}), we further have
\begin{align*}
	I(V;(\bX,D))\le \max_{v,v^\prime\in\{1,\ldots,2^{m^{p}/8}\}} D_{\rm KL}(P_{v}\mid\mid P_{v^\prime}),
\end{align*}
where $D_{\rm KL}(P_{v}\mid\mid P_{v^\prime})$ denotes the KL divergence between two distributions $P_{v}$ and $P_{v^\prime}$, and $P_v$ stands for the joint distribution of $(\bX,D)$ under the model $D=f_v(\bX)+\epsilon$. For convenience, we respectively let $p_v$ and $p_{v^\prime}$ be the density function of the joint distribution $P_{v}$ and $P_{v^\prime}$, and denote $p_\epsilon$ as the conditional density of $\epsilon$ given $\bX$. By definition of the $\chi^2$ divergence, we get
\begin{align*}
	D_{\rm KL}(P_{v}\mid\mid P_{v^\prime})&\overset{\mbox{\footnotesize(a)}}{\leq} D_{\rm \chi^2}(P_{v}\mid\mid P_{v^\prime})\\
 &=\mathbb{E}_{(\bX,D)\sim P_{v^\prime}}\left[\frac{\left({p_{v}(\bX,D)}-{p_{v^\prime}(\bX,D)}\right)^2}{p_{v'}^2(\bX,D)}\right]\\
 &=\mathbb{E}_{(\bX,D)\sim P_{v^\prime}}\left[\frac{\left({p_{v}(D\mid \bX)}-{p_{v^\prime}(D\mid \bX)}\right)^2}{p^2_{v^\prime}(D\mid \bX)}\right]\\
 &=\mathbb{E}_{(\bX,D)\sim P_{v^\prime}}\left[\frac{\left(p_{\epsilon}(D- f_v(\bX))-p_{\epsilon}(D- f_{v'}(\bX))\right)^2}{p^2_{v^\prime}(D\mid \bX)}\right]\\
	&\overset{\mbox{\footnotesize(b)}}{\leq} \frac{L_\epsilon}{\kappa^2}\mathbb{E}_{(\bX,D)\sim P_{v^\prime}}\left\vert f_{v}(\bX)-f_{v^\prime}(\bX)\right\vert^2\\
 &\overset{\mbox{\footnotesize(c)}}{\leq} \frac{L_\epsilon C_X}{\kappa^2}\|f_v-f_{v'}\|_2^2\\
    &\overset{\mbox{\footnotesize(d)}}{\leq} \frac{L_\epsilon C_X}{\kappa}\times \left[\frac{C_2\cdot\omega}{m^{2\beta}}\right]\\
	&\le \frac{C_2\cdot\omega}{m^{2\beta}},
\end{align*}
where (a) is a standard relationship between KL and $\chi^2$ divergence measures (see, e.g., Lemma 6 of \citealt{scarlett2019introductory}); (b) holds because the density function of the conditional distribution $D|\bX$ is uniformly bounded below by $\kappa$, which is $L_\epsilon$-Lipschitz continuous (and therefore, $p_\epsilon$ is also $L_\epsilon$-Lipschitz continuous); (c) comes from the fact that the marginal density distribution of $\bX$ is bounded above by $C_X$; (d) is valid by \eqref{lower_b1}.

Combining the obtained bounds, for any estimator $\hat{f}_n$ based on the sample $\bS_n=\{(\bX_i,D_i)\}_{i=1}^n$, we have
\begin{align*}
	\sup_{f_\rho\in \mathcal{H}^{\beta}([0,1]^p,B_0)} \mathbb{E} \Vert \hat{f}_n-f_\rho\Vert_{2}^2 &\ge\frac{C_1\cdot\omega}{m^{2\beta}} \times \left(1-\frac{C_2\cdot\omega\cdot n}{m^{p+2\beta}}\right),
\end{align*}
where we have left those finite numbers to the universal constants. 
By choosing $m=\lfloor n^{1/(p+2\beta)}\rfloor$ and a proper $\omega$, we can reach the following conclusion:
\begin{align*}
	\inf_{\hat{f}_n}\sup_{f_\rho\in \mathcal{H}^{\beta}([0,1]^p,B_0)} \mathbb{E} \Vert \hat{f}_n-f_\rho\Vert_{2}^2 &\ge C\times n^{-2\beta/(p+2\beta)}.
\end{align*}
{\color{black}Finally, we can apply Theorem \ref{thm:self_calibration} in Section \ref{sec:misspecified_models} to get the lower bound for excess risk:}
\begin{align*}
	\inf_{\hat{f}_n}\sup_{f_\rho\in \mathcal{H}^{\beta}([0,1]^p,B_0)} \mathbb{E} \Big[ \mathcal{R}^\rho(\hat{f}_n)-\mathcal{R}^\rho(f_\rho)\Big] &\ge C\times n^{-2\beta/(p+2\beta)}. \Halmos
\end{align*}


\section{Proofs in Section \ref{sec:extensions}}

    For the completeness of the proof, we present the definition of sequential covering number in the following \citep{kuznetsov2015learning}.

\begin{definition}\label{definition:sequential_covering}
A $\mathcal{Z}$-valued complete binary tree $z$ is a sequence $(z_1,\ldots,z_T)$ of $T$ mappings $z_t:\{\pm 1\}^{t-1}\to\mathcal{Z},t\in[1,T]$. A path in the tree is $\sigma=(\sigma_1,\ldots,\sigma_{T-1}).$ Simplify, we write $z_t(\bf{\sigma})$ for $z_t(\sigma_1,\ldots,\sigma_{t-1})$. A set $V$ of $\mathbb{R}$-valued trees of the depth $T$ is a sequential $\alpha$ cover(with respect to $\bf{q}$-weighted $\ell_p$ norm) of a function class $\mathcal{G}$ on a tree $z$ of the depth $T$ if for all $g\in\mathcal{G}$ and all $\sigma\in\{\pm 1\}^T$, there is $v\in V$ such that 
$$\left( \sum_{t=1}^T\vert \bf{v}_t({\bf \sigma})-g({\bf z_t}({\bf \sigma}))\vert^p\right)^{\frac{1}{p}}\le \Vert {\bf q}\Vert_q^{-1}\alpha,$$
    where $\Vert\cdot\Vert_q$ is the dual norm and ${\bf q}=(q_1,\ldots,q_T)$ is an arbitrary sequence.  The (sequential) covering number $\mathcal{N}_p(\alpha,\mathcal{G},\bf{z})$ of a function class $\mathcal{G}$ on a given tree $\bf{z}$ is defined to be the size of the minimal sequential cover.
\end{definition}

    \proof{Proof of Theorem \ref{thm:excess_risk_dependent}.}
    We write $L_f(z)= L(f(\bx),d):=b(d-f(\bx))^++h(f(\bx)-d)^+\leq b\Bar{D}+h\mathcal{B}$ for any $z=(\bx,d)\in\mathbb{R}^{d+1}$ and $f\in\mathcal{F}_{D\!N\!N}$. We define the function class $\mathcal{F}_L=\{L_f(\cdot):f\in\mathcal{F}_{D\!N\!N}\}$.

    To prove the theorem, we essentially follow the idea of risk decomposition in Lemma \ref{lemma:risk_decomposition}, but the treatment of the stochastic error is slightly different.
    Recall that the stochastic error with dependent data is defined as
    $$\sup_{f\in\mathcal{F}_{D\!N\!N}}\vert \mathcal{R}^\rho(f)-\mathcal{R}_T^\rho(f)\vert .$$

    By definition, we know that $\mathcal{F}_{D\!N\!N}$ contains networks with the parameter $\phi$, depth $\mathcal{D}$, width $\mathcal{W}$, size $\mathcal{S}$, and number of neurons $\mathcal{U}$ so that $\|f_{\phi}\|_{\infty}\leq\mathcal{B}$. 
    Also, it is obvious that the newsvendor loss function $L(y,y')$ is $\max\{b,h\}$-Lipschitz in both its two arguments. Then, by Corollary 2 in \cite{kuznetsov2015learning}, for any $\delta>0$, with probability at least $1-\delta$, for all $\epsilon>0$,
 \begin{align}\label{ineq1}
  \sup_{f\in\mathcal{F}_{D\!N\!N}}\vert \mathcal{R}^\rho(f)-\mathcal{R}_T^\rho(f)\vert\le 2\epsilon+\Delta
  + (b\bar{D}+h\mathcal{B}) \frac{\Vert \bf{q}\Vert_2}{\sqrt{T}}\sqrt{\frac{2\log \mathbb{E}_{S_T}[\mathcal{N}_1(\epsilon,\mathcal{F}_L,S_T)]+\log(1/\delta)}{T}},
\end{align}
where $\mathcal{N}_1(\epsilon,\mathcal{F}_L, S_T)$ is the (sequential) covering number of the function class $\mathcal{F}_L$ with radius $\epsilon$ based on samples $S_T$,  ${\bf q}=(q_1,\ldots,q_T)$ is an arbitrary sequence, and $\Delta$ is the discrepancy defined by
$$\Delta=\sup_{L\in\mathcal{F}_L}\Big( \mathbb{E}[L(Z_{T+1})\mid Z_1,\ldots,Z_T]-\sum_{t=1}^T q_t\mathbb{E}[L(Z_{t})\mid Z_1,\ldots,Z_{t-1}]\Big).$$

Recall that for any $z=(x,d)\in\mathbb{R}^{d+1}$ and any $f_1,f_2\in\mathcal{F}_{D\!N\!N}$, $L_{f_1}(z)=L_{f_1}(x,d)$ and $L_{f_2}(z)=L_{f_2}(x,d)$, and $\vert L_{f_1}(z)-L_{f_2}(z)\vert\le \max\{b,h\} \vert f_1(x)-f_2(x)\vert.$ For any $\epsilon>0$, let $\mathcal{F}^*_{D\!N\!N}$ be a $(\epsilon/\max\{b,h\})$-covering of $\mathcal{F}_{D\!N\!N}$ with number $\mathcal{N}_T(\epsilon/\max\{b,h\},\mathcal{F}_{D\!N\!N},\Vert\cdot\Vert_\infty)$. Then it is easy to check that $\mathcal{F}^*_{L}=\{L_f: f\in\mathcal{F}^*_{D\!N\!N}\}$ is a $\epsilon$-covering of $\mathcal{F}_L$ with ${\bf q}=(1,\ldots,1)$ (see Definition \ref{definition:sequential_covering}). Consequently, we have $\mathcal{N}_1(\epsilon,\mathcal{F}_L,S_T)\le \mathcal{N}_T(\epsilon/\max\{b,h\},\mathcal{F}_{D\!N\!N},\Vert\cdot\Vert_\infty)$ for any $S_T$. We choose $\epsilon=1/T$, 
 by Theorem \ref{thm:covering_number} and Theorem 3 in \cite{bartlett2019nearly}, we know that for $T\ge \text{Pdim}(\mathcal{F}_{D\!N\!N})$,
 \begin{align*}
 \log\mathbb{E}_{S_T}[\mathcal{N}_1(\epsilon,\mathcal{F}_L,S_T)]&\le \log\mathcal{N}_T(\epsilon/\max\{b,h\},\mathcal{F}_{D\!N\!N},\Vert\cdot\Vert_\infty)\\
 &\le \text{Pdim}(\mathcal{F}_{D\!N\!N})\log\Big( \frac{eT\mathcal{B}\max\{b,h\}}{\epsilon\text{Pdim}(\mathcal{F}_{D\!N\!N})}\Big)\\
 &\le \text{Pdim}(\mathcal{F}_{D\!N\!N})\log(eT^2\mathcal{B}\max\{b,h\})\\
 &\le C \mathcal{B}\max\{b,h\}\mathcal{D}\mathcal{S}\log(\mathcal{S})\log(T),
 \end{align*}
where $C>0$ is some universal constant. Together with (\ref{ineq1}), for a sufficient large $T$ and any $\delta>0$, we have, with probability at least $1-\delta$,
\begin{align}\notag
  \sup_{f\in\mathcal{F}_{D\!N\!N}}\vert \mathcal{R}^\rho(f)-\mathcal{R}_T^\rho(f)\vert &\le C [\max\{b,h\}(\bar{D}\mathcal{B})]^2 \sqrt{\frac{2\mathcal{D}\mathcal{S}\log(\mathcal{S})\log(T)}{T}} \\ \label{ineq2}
  & \qquad\qquad\qquad\qquad+(b\bar{D}+h\mathcal{B})\sqrt{\frac{\log(1/\delta)}{T}} +\Delta,
\end{align}
where $C>0$ is some universal constant. Then by Proposition \ref{prop:approximation_error_new} and the risk decomposition in Lemma \ref{lemma:risk_decomposition}, the following inequality holds with at least probability $1-\delta$,
\begin{align*}\notag
   \mathcal{R}^\rho(\hat{f}_{D\!N\!N})-\mathcal{R}^\rho(f_\rho) &\le C [\max\{b,h\}(\bar{D}\mathcal{B})]^2 \sqrt{\frac{2\mathcal{D}\mathcal{S}\log(\mathcal{S})\log(T)}{T}} +2(b\bar{D}+h\mathcal{B})\sqrt{\frac{\log(1/\delta)}{T}} +2\Delta\\
   & \qquad\qquad\qquad +18\max\{b,h\}B_0 (\lfloor\beta\rfloor+1)^2p^{\lfloor\beta\rfloor+(\beta\vee 1)/2}(NM)^{-2\beta/p},
\end{align*}
    where $T\ge C\cdot \mathcal{S}\mathcal{D}\log(\mathcal{S})$ for a large enough $C>0$ and $M,N\in\mathbb{N}_+$. Using the same arguments as in the proof of Theorem \ref{thm:excess_risk}, we can reach the conclusion stated in Theorem \ref{thm:excess_risk_dependent}.
\Halmos


\proof{Proof of Theorem \ref{thm:excess_risk_GLM}.}
    The excess risk can again be decomposed into a stochastic error and an approximation error as stated in Lemma \ref{lemma:risk_decomposition}. Specifically, the stochastic error bound remains the same as in Theorem \ref{thm:stochastic_error} since it is independent of the target function by definition in Lemma \ref{lemma:risk_decomposition}. Therefore, we only need to check the approximation error under this GLM setting.
    
    We first notice that the target function can be represented in the following composite form:
    \begin{equation*}
        f_{\rho}(\bx)=g(\btheta^\top \bx)=h_1\circ h_0(\bx),
    \end{equation*}
    where
    \begin{equation*}
        h_1(y)=g(y) \text{ for } y\in\mathbb{R}\quad \text{ and }\quad h_0(\bx)=\theta^\top\bx\text{ for } \bx\in\mathbb{R}^p.
    \end{equation*}
    Such a composite structure is indeed similar to the DNN structure \eqref{network_functional_composition} so that we can first use two sub-networks $\Tilde{h}_1$ and $\Tilde{h}_0$ to approximate $h_1$ and $h_0$, respectively. Next, by stacking these two sub-networks, we can construct a new network in approximating $f_\rho$. 
    Specifically, the domain of the univariate function $h_1$ is $[a,c]$, 
    and to apply Theorem 3.3 in \cite{jiao2023deep} with a domain being $[0,1]$, we first need an additional invertible linear layer $A(\cdot): [a,c]\rightarrow [0,1]$ as the first layer in the sub-network $\Tilde{h}_1$ to finish the transformation of the domain. Indeed, this invertible layer can be easily established via $A(y)=\sigma(\frac{1}{c-a}y-\frac{a}{c-a})$ where $\sigma(\cdot)$ is the ReLU activation function. Then, by Theorem 3.3 in \cite{jiao2023deep},
    there exists an $\Tilde{h}_1\in\mathcal{F}_{D\!N\!N}$ with the width $\mathcal{W}=38(\lfloor\beta\rfloor+1)^2 N\lceil\log_2(8N)\rceil$ and depth $\mathcal{D}=21(\lfloor\beta\rfloor+1)^2M\lceil\log_2(8M)\rceil+1$  for $M,N\in\mathbb{N}_+$, such that 
    \begin{equation}\label{approximation_univariate}
        |\Tilde{h}_1(y)-h_1(y)|\leq 18B_0 (\lfloor\beta\rfloor+1)^2(NM)^{-2\beta},
    \end{equation}
    for $y\in \mathbb{R}$ except a negligible set with an arbitrarily small Lebesgue measure.

    On the other hand, for the construction of sub-network $\Tilde{h}_0$, we know from Lemma \ref{lemma:linear_functions_nn} that there exists a one-layer network that can exactly represent a linear function. In other words, there is an $\Tilde{h}_0\in\mathcal{F}_{D\!N\!N}$ with the width of each layer being $(p,2,1)$, such that
    \begin{equation}\label{approximation_linear}
        |\Tilde{h}_0(\bx)-h_0(\bx)|=0, 
    \end{equation}
    for any $\bx\in\mathbb{R}^p$.

    Since the mappings $A(\cdot)$ and $\Tilde{h}_0$ are both continuous, and the probability measure of $\bX$ is absolutely continuous with respect to the Lebesgue measure by Assumption \ref{assumption_2}(i), combining \eqref{approximation_univariate} and \eqref{approximation_linear}, we then have  
    \begin{equation*}
        \mathbb{E}|\Tilde{h}_1\circ\Tilde{h}_0(\bx)-h_1\circ h_0(\bx)|\leq 18B_0 (\lfloor\beta\rfloor+1)^2(NM)^{-2\beta}.
    \end{equation*}
    Therefore, the approximation error can be bounded as follows:
    \begin{align*}
        \underset{f\in\mathcal{F}_{D\!N\!N}}{\inf}\left[\mathcal{R}^{\rho}(f)-\mathcal{R}^{\rho}(f_{\rho})\right]&\leq \mathcal{R}^{\rho}(\Tilde{h}_1\circ \Tilde{h}_0)-\mathcal{R}^{\rho}(h_1\circ h_0)\\
        &\overset{\mbox{\footnotesize(a)}}{\leq} \max\{b,h\}\mathbb{E}|\Tilde{h}_1\circ \Tilde{h}_0(\bx)-h_1\circ h_0(\bx)|\\
        &\leq 18\max\{b,h\}B_0 (\lfloor\beta\rfloor+1)^2(NM)^{-2\beta},
    \end{align*}
    where $(a)$ holds because the newsvendor loss function is $\max\{b,h\}$-Lipschitz. Together with Lemma \ref{lemma:risk_decomposition} and Theorem \ref{thm:stochastic_error}, the following inequality holds for at least probability $1-\delta$,
    \begin{align*}
        \mathcal{R}^{\rho}(\hat{f}_{D\!N\!N})-\mathcal{R}^{\rho}(f_{\rho})&\leq 2\sqrt{2}(b\bar{D}+h\mathcal{B})\bigg(C_3\sqrt{\frac{\mathcal{S}\mathcal{D}\log(\mathcal{S})\log(n)}{n}}+\sqrt{\frac{\log(1/\delta)}{n}}\bigg)\\
         &\qquad+18\max\{b,h\}B_0 (\lfloor\beta\rfloor+1)^2(NM)^{-2\beta}.
    \end{align*}
    {\color{black}Using similar arguments as in the proof of Theorem \ref{thm:excess_risk}, we choose $MN=\lfloor n^{\frac{1}{4\beta+2}}\rfloor$ so that 
    the convergence order in the sample size is optimized. Under this setting, we get
    \begin{equation*}
        \mathcal{R}^{\rho}(\hat{f}_{D\!N\!N})-\mathcal{R}^{\rho}(f_{\rho})\leq 2\sqrt{2}(b\bar{D}+h\mathcal{B})\bigg(C_3(\lfloor\beta\rfloor+1)^4(\log(n))^2n^{-\frac{\beta}{2\beta+1}}+\sqrt{\frac{\log(1/\delta)}{n}}\bigg),
    \end{equation*}
    for $n$ large enough. The stacked network consisting of $\Tilde{h}_1$ and $\Tilde{h}_0$ then has the width $\tilde{\mathcal{W}}=\max\{\mathcal{W}, 2p\}$ and depth $\tilde{\mathcal{D}}=\mathcal{D}+3$.}
    \Halmos


\section{Model Misspecification Issues}\label{sec:misspecified_models}

{\color{black}
\cite{ban2019big} have made significant contributions by thoroughly exploring the importance of incorporating feature information into decision-making for the newsvendor problem, both theoretically and numerically. In particular, by adaptively integrating both the basic and nonlinearly transformed features into their decision models, they effectively address the model misspecification issues of the linear and kernel policies proposed in their work.
In this section, we echo their insights by demonstrating that misspecified models can result in non-vanishing errors. In contrast, as shown in Section \ref{sec:theoretical_guarantees}, the DNN method's excess risk can converge to zero for a large class of target functions, provided there is sufficient data.} {\color{black} 
In practical implementations, a DNN with large enough depth and width should suffice to obtain a good decision with a small approximation error. In addition, there is no need to specify parametric forms of approximation rules, making it convenient in various practical scenarios. These advantages help mitigate the misspecification issues commonly found in learning problems. In our context, we consider the following definition of model misspecification.}

\begin{definition}[Model Misspecification]
Let $\{\mathcal{F}_n\}_{n=1}^{\infty}$ be the collection of hypothesis spaces and $f_\rho$ be the \emph{unknown} target function. The model is considered misspecified if $f_\rho\notin \underset{n\rightarrow\infty}{\lim}\overset{n}{\underset{i=1}{\cup}}\mathcal{F}_n.$ 
\end{definition}

{\color{black}Notice that we allow the hypothesis set $\mathcal{F}_n$ to expand with an increasing sample size to align with the formulation in Theorem \ref{thm:excess_risk}. 
}

    \begin{assumption}[{\sc Local Minimal Separation}]\label{assumtion:local_minimum_seperation}
        There exist $\gamma\in(0,\Bar{D}]$ and $\kappa>0$ such that for any $|\zeta|\leq \gamma$,
        \begin{equation*}
            \vert F(f_{\rho}(\bx)+\zeta\vert \bx)-F(f_{\rho}(\bx)\vert \bx) \vert \geq \kappa |\zeta|
        \end{equation*}
        for any $\bx\in \mathcal{X}$ up to a zero-probability set under the probability measure of $\bX$, where $F(\cdot\vert\bx)$ is the conditional CDF of $D$ given $\bX=\bx$.
    \end{assumption}

    Assumption \ref{assumtion:local_minimum_seperation} requires that the conditional density has a positive lower bound around the quantile function, which also warrants the strong convexity of the cost function in this region. 

    \begin{theorem}[{\sc Excess Risk w.r.t. Decision Bias}]\label{thm:self_calibration}
        Under Assumptions \ref{assumption_1}, \ref{assumption_2}, and \ref{assumtion:local_minimum_seperation}, for any $f: \mathcal{X}\rightarrow\mathbb{R}$, we have
        \begin{equation*}
            \mathcal{R}^{\rho}(f)-\mathcal{R}^{\rho}(f_{\rho})\geq \frac{\gamma\kappa(b+h)}{2(\Bar{D}+\gamma)} \Vert f -f_\rho\Vert^2_{L^2(\bX)}.
        \end{equation*}
    \end{theorem}

    {\color{black} Theorem \ref{thm:self_calibration} characterizes the lower bound for the excess risk caused by any suboptimal decision $f$ in terms of its $L^2$-norm difference from the true optimal decision $f_{\rho}$. It is straightforward that if the decision bias is larger, this excess risk bound also gets larger. Generally speaking, we can arrive at the universal insight that the excess risk should increase with the deviation of $f$ from $f_{\rho}$, even though we may not have an analytical expression for the latter. Therefore, parametric (misspecified) models will result in non-vanishing excess risk even with an infinite amount of data.} 

     \proof{Proof of Theorem \ref{thm:self_calibration}.} 
    We follow the idea of \cite{steinwart2011estimating_appendix} to construct a \emph{self-calibration} function to finish the proof (see the proof of Theorem 2.7 therein). For any \emph{fixed} $\bx\in\bX$, the self-calibration function is defined as
    \begin{equation*}
        \varphi_{L}(\epsilon,\bx):=\inf_{t\in\mathbb{R}: |t-f_{\rho}(\bx)|\geq \epsilon} \mathbb{E}_{D|\bx}[L(t,D)]-\mathbb{E}_{D|\bx}[L(f_{\rho}(\bx),D)],
    \end{equation*}
    where $\mathbb{E}_{D|\bx}[\cdot]$ denotes the expectation taken with respect to the conditional distribution of $D$ given $\bX=\bx$ and $L(t,d):=b(d-t)^++h(t-d)^+$. By definition, the self-calibration function captures the sensitivity of $L$-risk around the minimizer $f_{\rho}(\bx)$, and it is named according to the property that
    \begin{equation}\label{self_calibration}
        \varphi_{L}(|t-f_{\rho}(\bx)|,\bx)\leq \mathbb{E}_{D|\bx}[L(t,D)]-\mathbb{E}_{D|\bx}[L(f_{\rho}(\bx),D)], \quad \text{for any }t\in\mathbb{R},
    \end{equation}
    if we let $\epsilon:=|t-f_{\rho}(\bx)|$. Since the newsvendor loss function $L$ is convex, so does the map $t\rightarrow \mathbb{E}_{D|\bx}[L(t,D)]-\mathbb{E}_{D|\bx}[L(f_{\rho}(\bx),D)],$ we then have
    \begin{equation}\label{self_calibration_bound}
        \varphi_{L}(\epsilon,\bx)=\min\{\mathbb{E}_{D|\bx}[L(f_{\rho}(\bx)+\epsilon,D)], \mathbb{E}_{D|\bx}[L(f_{\rho}(\bx)-\epsilon,D)]\}-\mathbb{E}_{D|\bx}[L(f_{\rho}(\bx),D)].
    \end{equation}
    We then examine each term individually. 
    First, for the first term, if $\epsilon\in[0,\gamma]$, 
    \begin{align*}
        \mathbb{E}_{D|\bx}[L(f_{\rho}(\bx)+\epsilon,D)]-\mathbb{E}_{D|\bx}[L(f_{\rho}(\bx),D)]&\overset{\mbox{\footnotesize(a)}}{=} (b+h)\int_{f_{\rho}(\bx)}^{f_{\rho}(\bx)+\epsilon}[F(y|\bx)-\rho]dy\\
        &=(b+h)\int_{0}^{\epsilon}[F(y+f_{\rho}(\bx)|\bx)-F(f_{\rho}(\bx)|\bx)]dy\\
        &\overset{\mbox{\footnotesize(b)}}{\geq} (b+h)\int_{0}^{\epsilon}\kappa y dy\\
        &=\frac{1}{2}(b+h)\kappa\epsilon^2,
    \end{align*}
    where equation $(a)$ follows from a standard integration-by-parts argument (see, e.g., Lemma A-1 of \citealt{besbes2023big_appendix}); and inequality $(b)$ holds due to Assumption \ref{assumtion:local_minimum_seperation}. On the other hand, if $\epsilon\in[\gamma, \Bar{D}]$,
    \begin{align*}
        \mathbb{E}_{D|\bx}[L(f_{\rho}(\bx)+\epsilon,D)]-\mathbb{E}_{D|\bx}[L(f_{\rho}(\bx),D)]
        &=(b+h)\int_{0}^{\epsilon}[F(y+f_{\rho}(\bx)|\bx)-F(f_{\rho}(\bx)|\bx)]dy\\
        &\overset{\mbox{\footnotesize(a)}}{\geq} (b+h)\int_{0}^{\gamma}\kappa y dy + (b+h)\int_{\gamma}^{\epsilon}\kappa \gamma dy\\
        &=\frac{1}{2}(b+h)\kappa\gamma^2+(b+h)\kappa\gamma(\epsilon-\gamma)\\
        &=(b+h)\kappa(\gamma\epsilon-\frac{1}{2}\gamma^2),
    \end{align*}
    where the inequality $(a)$ applies Assumption \ref{assumtion:local_minimum_seperation} and also since $F(y|\bx)$ is non-decreasing so that $F(y|\bx)\geq F(f_{\rho}(\bx)+\gamma|\bx)$. By Lemma \ref{lemma:quadratic_lower_function} presented below, we know that 
    \begin{equation*}
        \mathbb{E}_{D|\bx}[L(f_{\rho}(\bx)+\epsilon,D)]-\mathbb{E}_{D|\bx}[L(f_{\rho}(\bx),D)]\geq \frac{\gamma\kappa(b+h)}{2(\Bar{D}+\gamma)}\epsilon^2, \quad \text{for any }\epsilon\in[0,\Bar{D}].
    \end{equation*}
    For the second term, we can also analogously show that
    \begin{equation*}
        \mathbb{E}_{D|\bx}[L(f_{\rho}(\bx)-\epsilon,D)]-\mathbb{E}_{D|\bx}[L(f_{\rho}(\bx),D)]\geq \frac{\gamma\kappa(b+h)}{2(\Bar{D}+\gamma)}\epsilon^2, \quad \text{for any }\epsilon\in[0,\Bar{D}].
    \end{equation*}

    As a result, for any function $f:\mathcal{X}\rightarrow \mathbb{R}$ and a fixed $\bx\in\mathcal{X}$, by letting $t:=f(\bx)$ and $\epsilon:=|f(\bx)-f_{\rho}(\bx)|$, the self-calibration property \eqref{self_calibration} and \eqref{self_calibration_bound} give
    \begin{align*}
        \mathbb{E}_{D|\bx}[L(f(\bx),D)]-\mathbb{E}_{D|\bx}[L(f_{\rho}(\bx),D)]&\geq \varphi_{L}(|f(\bx)-f_{\rho}(\bx)|,\bx)\\
        &\geq \frac{\gamma\kappa(b+h)}{2(\Bar{D}+\gamma)}|f(\bx)-f_{\rho}(\bx)|^2.
    \end{align*}
    Taking expectation with respect to the distribution of $\bX$ on both sides, we obtain the desired inequality in Theorem \ref{thm:self_calibration}.
    \Halmos

    \begin{lemma}\label{lemma:quadratic_lower_function}
        For $\epsilon\in[0,\Bar{D}]$ and the function $h: [0,\Bar{D}]\rightarrow [0,\infty)$ defined by
        \begin{equation*}
            h(\epsilon):=\begin{cases}
            \frac{1}{2}(b+h)\kappa\epsilon^2, \quad &\text{if }\epsilon\in[0,\gamma],\\
            (b+h)\kappa(\gamma\epsilon-\frac{1}{2}\gamma^2), \quad &\text{if }\epsilon\in[\gamma,\Bar{D}],
            \end{cases}
        \end{equation*}
        we have
        \begin{equation*}
            h(\epsilon)\geq \frac{\gamma\kappa(b+h)}{2(\Bar{D}+\gamma)}\epsilon^2, \quad\text{ for all }\epsilon\in[0,\Bar{D}].
        \end{equation*}
    \end{lemma}
    \proof{Proof of Lemma \ref{lemma:quadratic_lower_function}.}
    When $\epsilon\in[0,\gamma]$, because $\gamma\leq \Bar{D}$, it holds that
    \begin{equation*}
        h(\epsilon)=\frac{1}{2}(b+h)\kappa\epsilon^2 = \frac{\gamma}{\gamma+\gamma}(b+h)\kappa\epsilon^2\geq \frac{\gamma\kappa(b+h)}{\Bar{D}+\gamma}\epsilon^2\geq \frac{\gamma\kappa(b+h)}{2(\Bar{D}+\gamma)}\epsilon^2.
    \end{equation*}
    When $\epsilon\in[\gamma,\Bar{D}]$, we define 
    \begin{equation*}
        g(\epsilon):=h(\epsilon)-\frac{\gamma\kappa(b+h)}{2(\Bar{D}+\gamma)}\epsilon^2= (b+h)\kappa(\gamma\epsilon-\frac{1}{2}\gamma^2)-\frac{\gamma\kappa(b+h)}{2(\Bar{D}+\gamma)}\epsilon^2,
    \end{equation*}
    whose first-order derivative is
    \begin{equation*}
        g'(\epsilon)=(b+h)\kappa\gamma-\frac{\gamma\kappa(b+h)}{\Bar{D}+\gamma}\epsilon=(b+h)\kappa\gamma\frac{\Bar{D}+\gamma-\epsilon}{\Bar{D}+\gamma}\geq 0.
    \end{equation*}
    Therefore, $g(\epsilon)$ is non-decreasing on $[\gamma,\Bar{D}]$, implying
    \begin{equation*}
        g(\epsilon)\geq g(\gamma)=\frac{1}{2}(b+h)\kappa\gamma^2\frac{\Bar{D}}{\Bar{D}+\gamma}\geq 0,
    \end{equation*}
    from which we can claim the result in Lemma \ref{lemma:quadratic_lower_function}.
    \Halmos


{\color{black}
\section{Supplementary Materials for Section \ref{sec:simulation_implementation}}

In this section, we present additional numerical experiments and results to supplement the analysis in Section \ref{sec:simulation_implementation}. First, we examine the claim from Section \ref{subsec:numerical_depth_width} that the DNN method can achieve comparably good performance across a range of network structures, demonstrating its flexibility. Beyond Figure \ref{fig:excess_risk_surface}, this phenomenon holds robustly across different sample sizes (see Figures \ref{fig:excess_risk_surface_samples_256} and \ref{fig:excess_risk_surface_samples_512}), various critical levels used (see Figures \ref{fig:excess_risk_surface_samples_1024_0.5} and \ref{fig:excess_risk_surface_samples_1024_0.75}), as well as different model setups (as illustrated in the subplots of each figure).

\begin{figure}[!ht]
	\centering
 \includegraphics[scale=0.436]{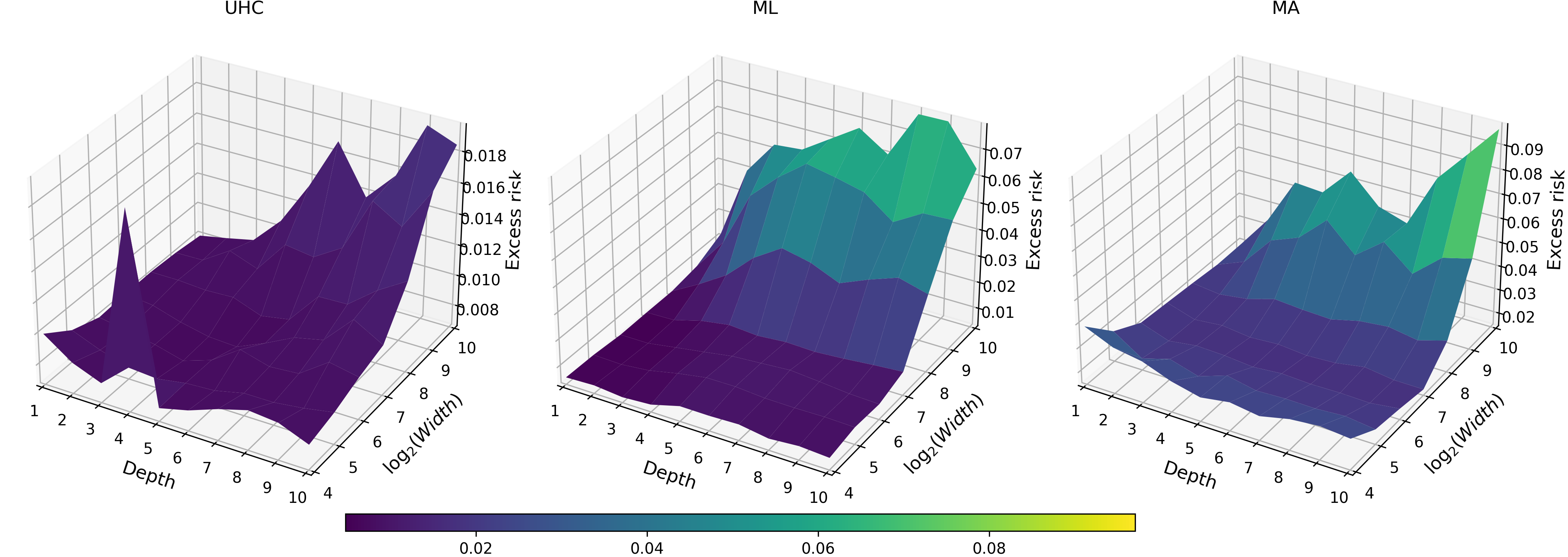}	
	\caption{Excess risk against network depth and width when $n=256$ and $\rho=0.25$.}
		\label{fig:excess_risk_surface_samples_256}
\end{figure}

\begin{figure}[!ht]
	\centering
 \includegraphics[scale=0.436]{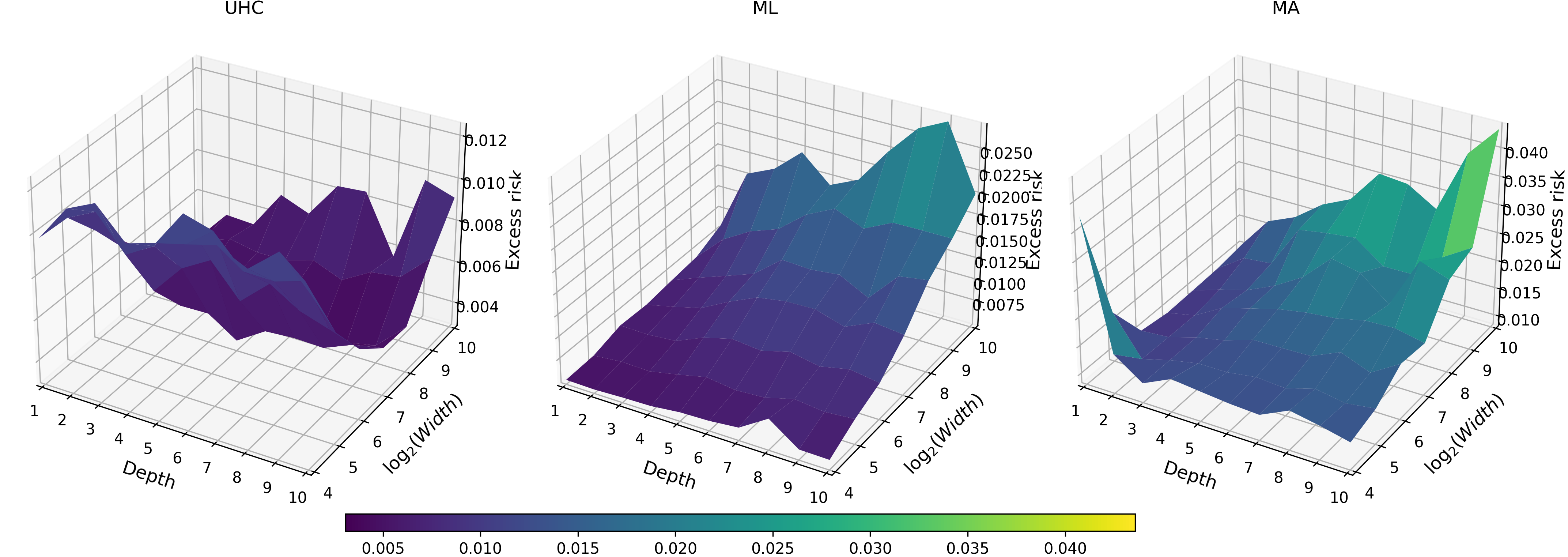}	
	\caption{Excess risk against network depth and width when $n=512$ and $\rho=0.25$.}
		\label{fig:excess_risk_surface_samples_512}
\end{figure}

\begin{figure}[!ht]
	\centering
 \includegraphics[scale=0.436]{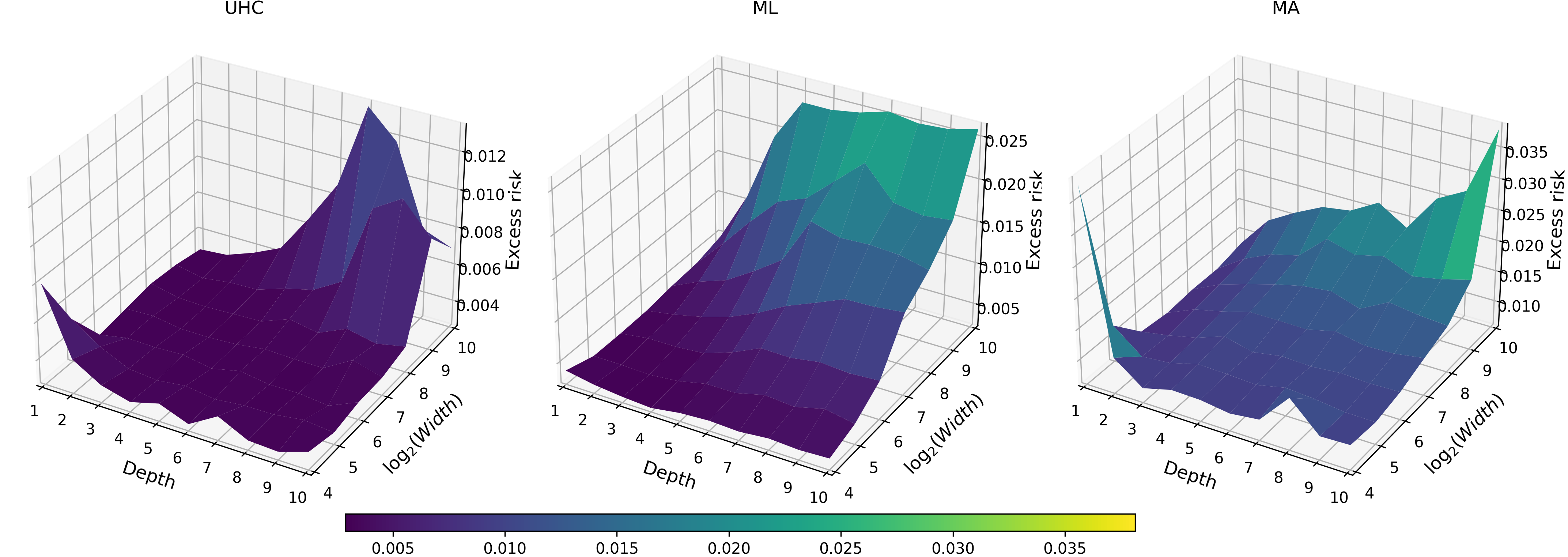}	
	\caption{Excess risk against network depth and width when $n=1024$ and $\rho=0.5$.}
		\label{fig:excess_risk_surface_samples_1024_0.5}
\end{figure}

\begin{figure}[!ht]
	\centering
 \includegraphics[scale=0.436]{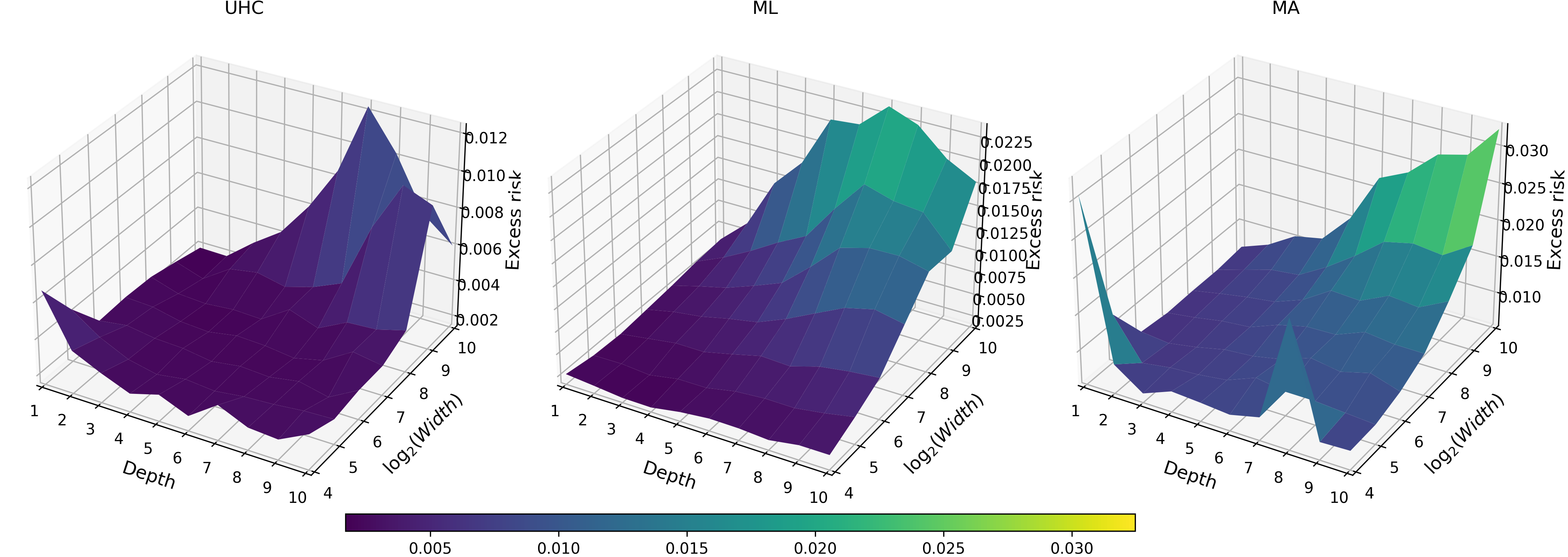}	
	\caption{Excess risk against network depth and width when $n=1024$ and $\rho=0.75$.}
		\label{fig:excess_risk_surface_samples_1024_0.75}
\end{figure}

As an interesting and independent trial, we further explore the impact of different activation functions on the performance of the DNN method. Specifically, we consider the Sigmoid and Tanh activation functions, defined respectively as
\begin{align*}
    \sigma_{Sigmoid}(x):=\frac{1}{1+e^{-x}} \quad \text{ and }\quad \sigma_{Tanh}(x):=\frac{e^x-e^{-x}}{e^x+e^{-x}}.
\end{align*}
Compared to these alternatives, the ReLU activation function is simpler and often avoids vanishing gradients, even when inputs are large. Moreover, as noted in Section \ref{subsec:DNN_model}, ReLU networks have gained widespread popularity in deep learning applications due to their computational efficiency, ease of training, and strong generalization performance (see, e.g., \citealt{nair2010rectified,krizhevsky2012imagenet}).

For a brief comparison, we present the numerical performance of DNNs using these two activation functions in Figures \ref{fig:excess_risk_Sigmoid_surface_samples_1024_0.75} and \ref{fig:excess_risk_Tanh_surface_samples_1024_0.75},  with other settings the same as those in Figure \ref{fig:excess_risk_surface}. We have the following observations. First, for the relatively simple tasks considered, the DNN method performs well with any of the activation functions once finely tuned, as the minimum loss across depth-width combinations remains consistently low. Second, Sigmoid and Tanh networks seem to be more sensitive to the network design, as their performance can deteriorate rapidly, evidenced by the wider range of errors across different figures. In particular, they perform relatively badly when the network gets deeper, which may be attributed to the vanishing gradient problem, as previously noted. Third, since both Sigmoid and Tanh activation functions introduce (more) nonlinearities, networks using these activation functions often require less depth to achieve a good approximation. 

\begin{figure}[!ht]
	\centering
 \includegraphics[scale=0.436]{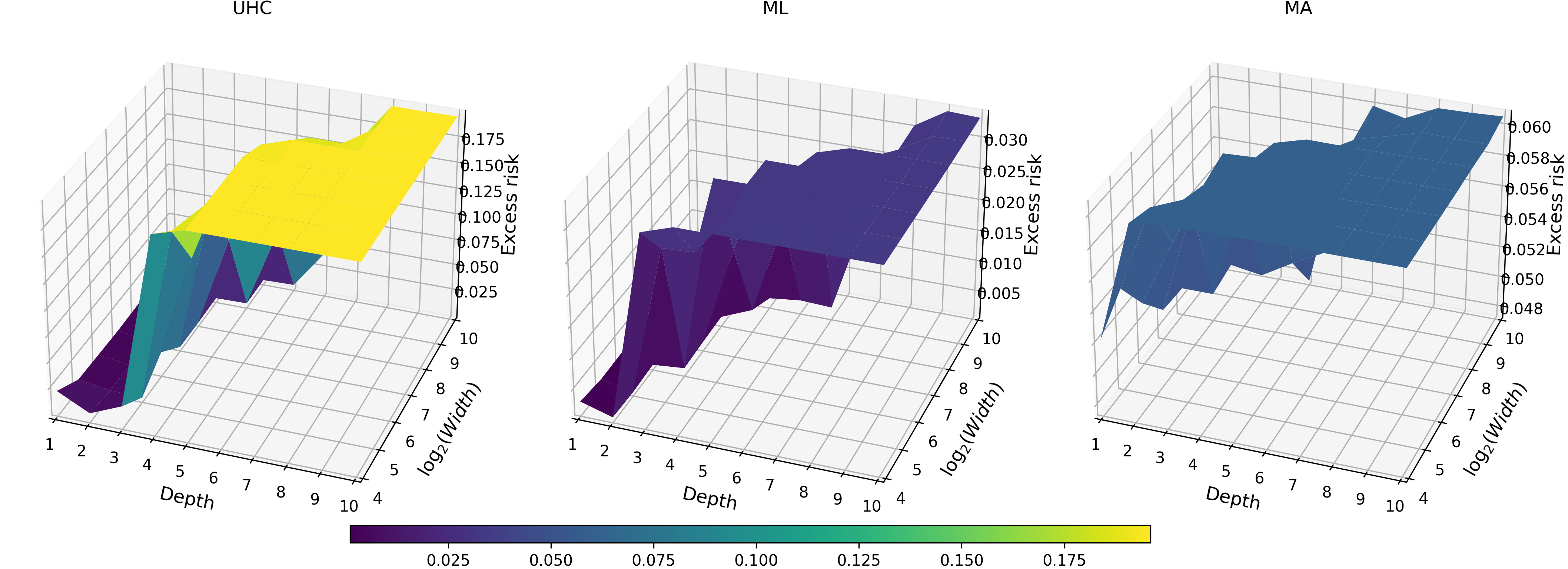}	
	\caption{Excess risk against network depth and width when using Sigmoid activation function ($n=1024, \rho=0.25$).}
		\label{fig:excess_risk_Sigmoid_surface_samples_1024_0.75}
\end{figure}

\begin{figure}[!ht]
	\centering
 \includegraphics[scale=0.436]{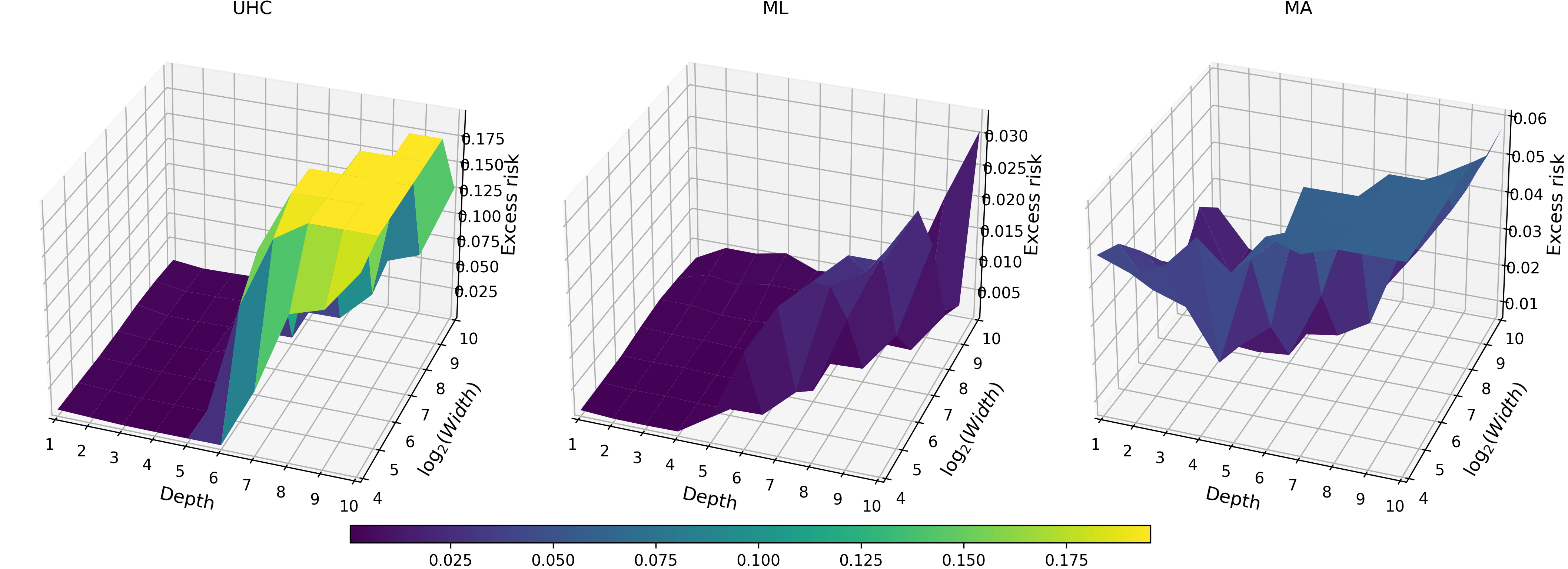}	
	\caption{Excess risk against network depth and width when using Tanh activation function ($n=1024, \rho=0.25$).}
		\label{fig:excess_risk_Tanh_surface_samples_1024_0.75}
\end{figure}

}


{\color{black}
\section{Supplementary Materials for Section \ref{sec:case_study}}

\subsection{Justification for the Selected Network Architecture}\label{sec_supp:justifincation_network_design}

Here, we provide the rationale for selecting a three-layer network with 512 neurons per layer for our numerical experiments in Section \ref{sec:case_study}. This choice is not made arbitrarily but is guided by the principles outlined in Section \ref{sec:simulation_implementation}. Specifically, as suggested in Section \ref{subsec:numerical_depth_width}, we first select a sufficiently wide network by setting the width to 512 and then progressively evaluate the performance of networks with varying depths. The final choice is based on achieving strong performance while avoiding unnecessary complexity in the number of layers. Following this principle, we conduct experiments on both the preliminary and full datasets, using a critical level of $\rho=0.7$ as a simple illustration. As shown in Figure \ref{fig:testing_error_varying_depth}, the three-layer network achieves the best performance on both the preliminary and full datasets among different depths.
Additionally, the network's moderate depth makes implementation quite convenient.
Based on these observations, we select the three-layer network with 512 neurons in each layer as the benchmark for the DNN method in all numerical experiments in Section \ref{sec:case_study}.

\begin{figure}[!ht]
	\centering
 \includegraphics[scale=0.436]{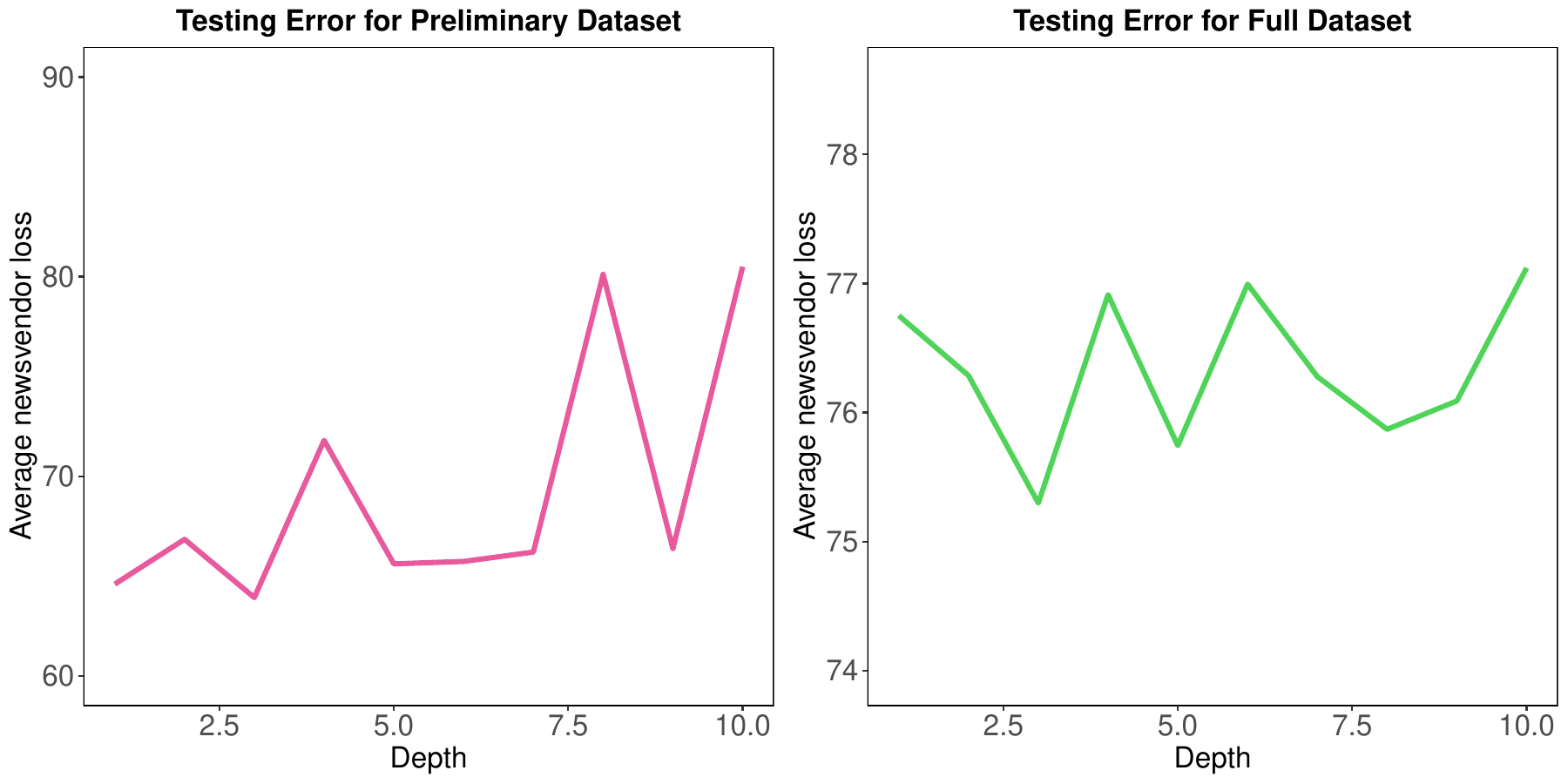}	
	\caption{Average newsvendor loss against network depth when width=512 and $\rho=0.7$.}
		\label{fig:testing_error_varying_depth}
\end{figure}
}

{\color{black}\subsection{Supplementary Experiment Outcomes}\label{subsec_appendix:more_numerical_outcomes}

Based on Tables \ref{table:sub-dataset:absolute} and \ref{table:full-dataset:absolute}, Figures \ref{fig:relative_loss_small_dataset} and \ref{fig:relative_loss_full_dataset} present the relative average newsvendor loss of alternative methods compared to the DNN, for the small and full datasets, respectively.

\begin{figure}[!ht]
	\centering
 \includegraphics[scale=0.466]{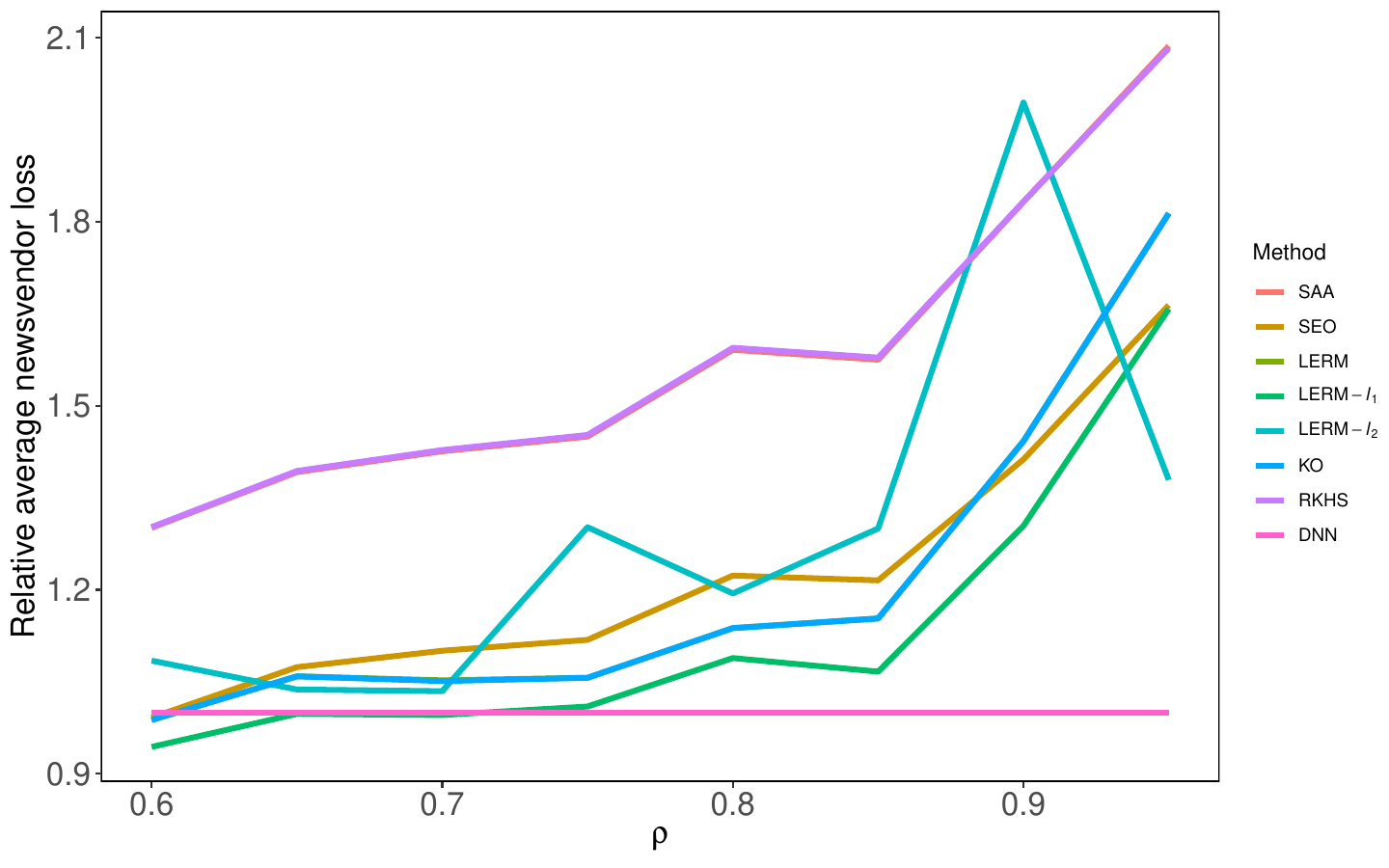}	
	\caption{Relative average newsvendor loss (compared to DNN) for different methods on the small dataset.}
		\label{fig:relative_loss_small_dataset}
\end{figure}

\begin{figure}[!ht]
	\centering
 \includegraphics[scale=0.466]{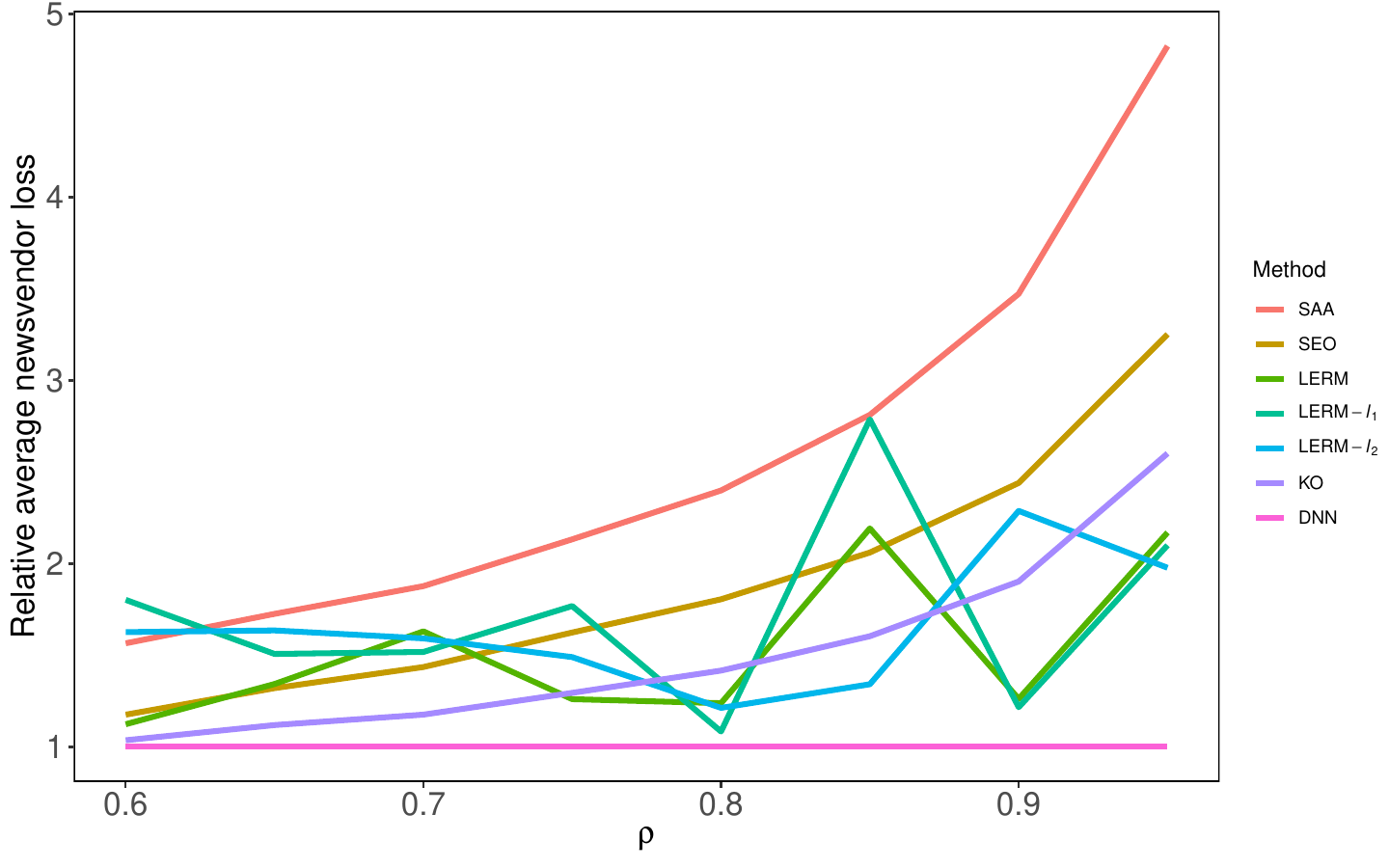}	
	\caption{Relative average newsvendor loss (compared to DNN) for different methods on the full dataset. The performance of RKHS is not shown because it significantly increases the scale of the plot.}
		\label{fig:relative_loss_full_dataset}
\end{figure}

Next, we present the training and execution times for each method on the full dataset (see Tables \ref{table:training-execution_times_full}). All numerical experiments were conducted on a computer equipped with a 12th Gen Intel(R) Core(TM) i9-12900 CPU (2.40 GHz) and 128 GB of RAM. In particular, unlike in \cite{ban2019big}, the KO method in our setting requires a few seconds to implement, likely due to the significantly higher dimensionality and larger size of the training data, which increases the time needed to compute the distances for the kernel weights.

\begin{table}[!ht]
				\caption{\color{black}Training and Execution times (in seconds) for different methods on the full dataset. Each cell displays the training time on the first line, followed by the execution time on the second line. A dash (``-'') indicates not applicable.}
			\centering\color{black}
		{\def\arraystretch{1.66}  
			\begin{tabular*}{\textwidth}{@{\extracolsep{\fill}}lcccccccc}
				\hline
				\hline
				$\rho$   & 0.6 & 0.65 & 0.7 &0.75  & 0.8  & 0.85 & 0.9 & 0.95 \\ \hline
				\parbox[c]{1cm}{SAA}  	&\parbox[c]{1.46cm}{\centering\scriptsize -\\$\approx0$}	&\parbox[c]{1.46cm}{\centering\scriptsize -\\$\approx0$}	&\parbox[c]{1.46cm}{\centering\scriptsize -\\$\approx0$}	&\parbox[c]{1.476cm}{\centering\scriptsize -\\$\approx0$}	&\parbox[c]{1.476cm}{\centering\scriptsize -\\$\approx0$}	&\parbox[c]{1.476cm}{\centering\scriptsize -\\$\approx0$}	&\parbox[c]{1.476cm}{\centering\scriptsize -\\$\approx0$} & \parbox[c]{1.46cm}{\centering\scriptsize -\\$\approx0$}
                     \\ \hline
                    \parbox[c]{1cm}{SEO}  	&\parbox[c]{1.476cm}{\centering\scriptsize 0.14\\0.002}	&\parbox[c]{1.476cm}{\centering\scriptsize 0.09\\0.002}	&\parbox[c]{1.476cm}{\centering\scriptsize 0.08\\0.003}	&\parbox[c]{1.476cm}{\centering\scriptsize 0.08\\0.002}	&\parbox[c]{1.476cm}{\centering\scriptsize 0.25\\0.002}	&\parbox[c]{1.66cm}{\centering\scriptsize 0.08\\0.002} & \parbox[c]{1.476cm}{\centering\scriptsize 0.08\\0.003}& \parbox[c]{1.476cm}{\centering\scriptsize 	0.34\\0.003}
                \\ \hline 
				\parbox[c]{0.86cm}{LERM}  	&\parbox[c]{1.476cm}{\centering\scriptsize 6.10\\0.001}	&\parbox[c]{1.476cm}{\centering\scriptsize 6.17\\0.002}	&\parbox[c]{1.66cm}{\centering\scriptsize 6.23\\0.001}	&\parbox[c]{1.476cm}{\centering\scriptsize 6.20\\0.001}	&\parbox[c]{1.476cm}{\centering\scriptsize 6.10\\0.001}	&\parbox[c]{1.66cm}{\centering\scriptsize 6.29\\0.001} & \parbox[c]{1.476cm}{\centering\scriptsize 6.17\\0.001}& \parbox[c]{1.476cm}{\centering\scriptsize 	6.07\\0.001}
                 \\ \hline  
    			{\small LERM-$\ell_1$} 	&\parbox[c]{1.66cm}{\centering\scriptsize 6.21\\0.001}	&\parbox[c]{1.66cm}{\centering\scriptsize 6.07\\0.001}	&\parbox[c]{1.476cm}{\centering\scriptsize 6.13\\0.001}	&\parbox[c]{1.476cm}{\centering\scriptsize 6.17\\0.001}	&\parbox[c]{1.476cm}{\centering\scriptsize 6.14\\0.001}	&\parbox[c]{1.66cm}{\centering\scriptsize 6.22\\0.001} & \parbox[c]{1.476cm}{\centering\scriptsize 6.20\\0.001}& \parbox[c]{1.476cm}{\centering\scriptsize 	6.13\\0.001}\\ \hline  
    			{\small LERM-$\ell_2$} 	&\parbox[c]{1.66cm}{\centering\scriptsize 8.85\\ 0.001}	&\parbox[c]{1.66cm}{\centering\scriptsize 13.96\\0.001}	&\parbox[c]{1.66cm}{\centering\scriptsize 7.35\\0.001}	&\parbox[c]{1.66cm}{\centering\scriptsize 6.27\\0.001}	&\parbox[c]{1.476cm}{\centering\scriptsize 6.31\\0.001}	&\parbox[c]{1.476cm}{\centering\scriptsize 6.36\\0.001} & \parbox[c]{1.66cm}{\centering\scriptsize 6.42\\0.001}& \parbox[c]{1.476cm}{\centering\scriptsize 	6.44\\0.001}
 \\ \hline 
				\parbox[c]{1cm}{KO}  	&\parbox[c]{1.476cm}{\centering\scriptsize -\\ 6.91}	&\parbox[c]{1.476cm}{\centering\scriptsize -\\7.15}	&\parbox[c]{1.476cm}{\centering\scriptsize -\\5.07}	&\parbox[c]{1.476cm}{\centering\scriptsize -\\10.52}	&\parbox[c]{1.476cm}{\centering\scriptsize -\\6.62}	&\parbox[c]{1.476cm}{\centering\scriptsize -\\8.74} & \parbox[c]{1.476cm}{\centering\scriptsize -\\6.96}& \parbox[c]{1.476cm}{\centering\scriptsize 	-\\11.76}
 \\ \hline  
                \parbox[c]{1cm}{RKHS}  &\parbox[c]{1.476cm}{\centering\scriptsize 368.26\\0.10}	&\parbox[c]{1.476cm}{\centering\scriptsize 367.55\\0.14}	&\parbox[c]{1.66cm}{\centering\scriptsize 433.51\\0.12}	&\parbox[c]{1.476cm}{\centering\scriptsize 430.40\\0.12}	&\parbox[c]{1.476cm}{\centering\scriptsize 430.60\\0.15}	&\parbox[c]{1.476cm}{\centering\scriptsize 387.65\\0.12} & \parbox[c]{1.476cm}{\centering\scriptsize 364.89\\0.11}& \parbox[c]{1.476cm}{\centering\scriptsize 	362.15\\0.10}
 \\ \hline 
				\parbox[c]{1cm}{DNN}  &\parbox[c]{1.476cm}{\centering\scriptsize 13.96\\0.01}	&\parbox[c]{1.476cm}{\centering\scriptsize 13.83\\0.01}	&\parbox[c]{1.476cm}{\centering\scriptsize 14.41\\0.01}	&\parbox[c]{1.476cm}{\centering\scriptsize 14.29\\0.01}	&\parbox[c]{1.476cm}{\centering\scriptsize 14.43\\0.01}	&\parbox[c]{1.476cm}{\centering\scriptsize 14.38\\0.01} & \parbox[c]{1.476cm}{\centering\scriptsize 13.93\\0.01}& \parbox[c]{1.476cm}{\centering\scriptsize 	13.50\\0.01} \\
				\hline
				\hline
			\end{tabular*}
			\label{table:training-execution_times_full}
		}
\end{table}

}

{\color{black}\section{Further Discussions}\label{sec_supp:further_discussion}

Although our main focus is on the single-product newsvendor problem for ease of exposition, our results naturally extend to the multi-product case. There are also many avenues for future research.  We also anticipate that similar theoretical results would apply to other end-to-end stochastic optimization problems where the objective function is Lipschitz continuous and convex in the decision variable. {\color{black} However, the same approach may not apply to the optimal pricing problem, where a mapping from prices to demand (or revenue) is usually learned from data, and then one needs to optimize the network function with respect to (some of) its inputs. In such cases, this additional optimization step requires other techniques, such as the gradient descent method \citep{chen2022using}.} For some other examples, there is currently a theoretical gap between the expected risk bound in Theorem \ref{thm:upper_expected_bound} and the high-probability risk bound in Theorem \ref{thm:excess_risk}, where the convergence rate of the former is faster than that of the latter. As shown in the proof in Section \ref{subsec_append:proof_expected_bound}, the faster rate in Theorem \ref{thm:upper_expected_bound} arises from a different yet tricky decomposition of the expected excess risk in Lemma \ref{lemma:risk_decomposition2}, compared to the decomposition in Lemma \ref{lemma:risk_decomposition}. While the approximation error is similar in both cases, the definition of the stochastic error differs. To this end, a variant of Bernstein inequality \citep{gyorfi2002distribution_appendix} is used to show that the stochastic error corresponding to the expected excess risk converges more rapidly in terms of the sample size, resulting in a faster convergence rate for the excess risk. 
While the high-probability bound offers strong probabilistic guarantees for the excess risk, the expectation bound only characterizes its behavior in the expected level. Therefore, improving the convergence rate of the high-probability bound, if possible, would likely require more advanced techniques for analyzing the stochastic error, such as the use of a more refined concentration inequality.
In addition, it would be interesting to integrate various constraints on the DNN solution, addressing practical considerations such as capacity limitations \citep{elmachtoub2023estimate_appendix}. While this could potentially be tackled by implementing a penalized loss function that accounts for the relevant constraint, a thorough examination of the corresponding theoretical and numerical performance is still required. Furthermore, we anticipate that the DNN method holds significant potential for solving other data-driven OM problems.

Finally, while we have established theoretical guarantees for using the DNN method in the considered newsvendor problem, there remains ample room for future research. To mention a few, first, we have mainly focused on analyzing worst-case error bounds so that the theory holds uniformly across a wide range of cases. However, maintaining such generality may suffer from the problem that the DNN method can perform better than what the worst-case bounds indicate in practical examples. For instance,  in Section \ref{subsec:understanding_convergence_rates}, we observe that the convergence rate of the excess risk in the sample size differs from the minimax optimal rate. Going beyond worst-case analysis, one can resort to other theoretical characterizations on helping explain why the DNN method generalizes so well. For instance, \cite{belkin2021fit_appendix} introduces the interpolation paradigm, showing that over-parameterized models can perfectly interpolate training data while still generalizing well. Another frequently employed method involves precise analyses of specific simple models, such as liner and ridge regressions \citep{bartlett2020benign_appendix,hastie2022surprises_appendix,tsigler2023benign_appendix}. For a comprehensive review of various approaches explaining DNN's generalization, readers can refer to \cite{bartlett2021deep_appendix} and \cite{theisen2023beyond_appendix}. Second, one can possibly explore the optimization side of the DNN method in OM applications. As discussed in Section \ref{sec:simulation_implementation}, providing a uniform theoretical analysis for the optimization error is highly challenging. Nevertheless, the success of DNNs is partially attributed to their optimization properties, which could be investigated both numerically and theoretically in simple setups. For example, \cite{neyshabur2015search_appendix} finds that the DNN  optimization algorithm introduces a bias that favors simpler models, which contributes to the phenomenon of implicit regularization. Along similar lines, \cite{soudry2018implicit_appendix} and \cite{ji2019implicit_appendix} analyze the implicit bias of gradient descent on linearly separable and non-separable data, respectively, as well as its effects on model generalization. In particular, \cite{arora2019fine_appendix} explores the connection between the optimization landscape and generalization performance in a two-layer network setting. 

}

\end{document}